\documentclass{article}

\usepackage[a4paper,hmargin=3cm,vmargin=3cm]{geometry}

\usepackage{amsmath}
\usepackage{amsfonts}
\usepackage{amssymb}
\usepackage{epsfig}
\usepackage{theorem}
\usepackage{comment}
\usepackage[dvips]{color}
\usepackage{color}
\usepackage{graphicx}

\newtheorem{theorem}{Theorem}[section]
\newtheorem{lemma}[theorem]{Lemma}
\newtheorem{corollary}[theorem]{Corollary}
\newtheorem{proposition}[theorem]{Proposition}

{\theorembodyfont{\rmfamily} 

\newtheorem{convention}[theorem]{Convention}
\newtheorem{construction}[theorem]{Construction}
\newtheorem{remark}[theorem]{Remark}
\newtheorem{example}[theorem]{Example}

\newtheorem{definition}[theorem]{Definition}

\newtheorem{assumption}[theorem]{Assumption}

}

\def\C{{\mathbb C}}

\def\Z{{\mathbb Z}}

\def\Q{{\mathbb Q}}

\newcommand\al{\alpha}

\newcommand\ga{\gamma}
\newcommand\la{\lambda}
\newcommand\be{\beta}
\newcommand\eps{\varepsilon}
\newcommand\Cbar{\overline{\C}}
\newcommand\Qbar{\overline{\Q}}
\newcommand\eq{\,=\,}

\newcommand\bsk{\bigskip}
\newcommand\msk{\medskip}
\newcommand\ssk{\smallskip}

\def\fleche{\longrightarrow}
\def\pusto{\varnothing}

\title{Minimum Degree of the Difference of Two Polynomials over~$\Q$, 
       and Weighted Plane Trees.} 
\author{
Fedor Pakovich\thanks{Department of Mathematics, Faculty of Natural
Sciences, Ben-Gurion University of the Negev, P.O.B. 653, Beer Sheva, 
Israel; {\tt pakovich@math.bgu.ac.il}.}\hspace{1.5mm}
and 
Alexander K. Zvonkin\thanks{LaBRI, Universit\'e Bordeaux I, 351 cours de la 
Lib\'eration, F-33405 Talence Cedex France; 
{\tt zvonkin@labri.fr}.}
}
\date{\today}

\begin{document}

\maketitle

\begin{abstract}
A {\em weighted bicolored plane tree}\/ (or just {\em tree}\/ for short) 
is a bicolored plane tree whose edges are endowed with positive integral 
weights. The degree of a vertex is defined as the sum of the weights of 
the edges incident to this vertex. Using the theory of {\em dessins 
d'enfants}\/, which studies the action of the absolute Galois group on 
graphs embedded into Riemann surfaces, we show that a weighted plane tree 
is a graphical representation of a pair of coprime polynomials 
$P,Q\in\C\,[x]$ such that: (a)~$\deg P = \deg Q$, and $P$ and $Q$ have 
the same leading coefficient; (b)~the multiplicities of the roots of~$P$ 
(respectively, of~$Q$) are equal to the degrees of the black (respectively, 
white) vertices of the corresponding tree; (c)~the degree of the difference 
$P-Q$ attains the minimum which is possible for the given multiplicities 
of the roots of $P$~and~$Q$. Moreover, if a tree in question is uniquely 
determined by the set of its black and white vertex degrees (we call such 
trees {\em unitrees}), then the corresponding polynomials are defined 
over $\Q$.

The pairs of polynomials $P,Q$ such that the degree of the difference 
$P-Q$ attains the minimum, and especially those defined over $\Q$, 
are related to some important questions of number theory. Dozens of 
papers, from 1965~\cite{BCHS-65} to 2010~\cite{BeuSte-10}, were dedicated 
to their study. The main result of this paper is a complete classification 
of the {\em unitrees}\/ which provides us with the most massive class of 
such pairs defined over~$\Q$. We also study combinatorial invariants 
of the Galois action on trees, as well as on the corresponding polynomial 
pairs, which permit us to find yet more examples defined over~$\Q$.
In a subsequent paper we compute the polynomials $P,Q$ corresponding 
to all the unitrees.
\end{abstract}

\setcounter{tocdepth}{3}
\tableofcontents

\section{Introduction}

In 1965, Birch, Chowla, Hall, and Schinzel \cite{BCHS-65} asked a question 
which soon became famous:\label{BCHS}
\begin{quote}
Let $A$ and $B$ be two coprime polynomials with complex coefficients; 
what is the possible minimum degree of the difference $R=A^3-B^2$\,? 
\end{quote}
It is reasonable to suppose that $A^3$ and $B^2$ have the same degree and 
the same leading coefficients. Let us take $\deg A = 2k$, $\deg B = 3k$, 
so that $\deg A^3 = \deg B^2 = 6k$. Then the following was conjectured 
in~\cite{BCHS-65}:\label{init-problem}

\begin{enumerate}
\item   For $R=A^3-B^2$ one always has\, $\deg R \ge k+1$. 
\item   This bound is sharp: that is, it is attained for infinitely many
        values of $k$. 
\end{enumerate}

The first conjecture was proved the same year by Davenport 
\cite{Davenport-65}. The second one turned out to be much more difficult
and remained open for 16 years: in 1981 Stothers \cite{Stothers-81}
showed that the bound is in fact attained not only for infinitely many
values of $k$ but for all of them. 

The above problem may be generalized in various ways. The following one
was considered in 1995 by Zannier \cite{Zannier-95}. Let $\al,\be\vdash n$
be two partitions of~$n$, 
$$\al=(\al_1,\ldots,\al_p), \quad \be=(\be_1,\ldots,\be_q), \quad 
\sum_{i=1}^p\al_i=\sum_{j=1}^q\be_j=n,
$$
and let $P$ and $Q$ be two coprime polynomials of degree $n$ having 
the following factorization pattern:
\begin{eqnarray}\label{eq:P-and-Q}
P(x) \eq \prod_{i=1}^p \, (x-a_i)^{\al_i}\,, \qquad
Q(x) \eq \prod_{j=1}^q \, (x-b_j)^{\be_j}\,.
\end{eqnarray}
In these expressions we consider the multiplicities $\al_i$ and $\be_j$, 
$i=1,2,\ldots,p$, $j=1,2,\ldots,q$ as being given, while the roots 
$a_i$ and $b_j$ are not fixed, though they must all be distinct. The problem 
is to find the minimum possible degree of the difference $R = P-Q$. 
In his paper, Zannier proved the following. Let $d={\rm gcd}\,(\al,\be)$ 
denote the greatest common divisor of the numbers 
$\al_1,\ldots,\al_p,\be_1,\ldots,\be_q$. If
\begin{eqnarray}\label{qwert}
p+q\,\le\,\frac{n}{d}+1
\end{eqnarray}
then
\begin{eqnarray}\label{eq:main-bound}
\deg R \,\ge\, (n+1)-(p+q),
\end{eqnarray}
and this bound is always attained. If, on the other hand, 
${\displaystyle p+q>\frac{n}{d}+1}$, then a weaker bound
\begin{eqnarray}\label{no-qwert}
\deg R \,\ge\, \frac{(d-1)\,n}{d},
\end{eqnarray}
is valid, and it is also attained. 

\begin{definition}[Davenport--Zannier triple]\label{def:DZ-triple}
Let $P,Q,R\in\C\,[x]$ be coprime polynomials with factorization pattern 
(\ref{eq:P-and-Q}), $\deg P=\deg Q=n$, while the degree of the polynomial 
$R=P-Q$ equals $(n+1)-(p+q)$. Then the triple $(P,Q,R)$ is called a 
{\em Davenport--Zannier triple}, or, in a more concise way, a 
{\em DZ-triple}.
\end{definition}

The main subject of this paper is a study of DZ-triples
{\em defined over}\/ $\Q$, that is, the triples $P,Q,R\in\Q\,[x]$.

\msk

The paper is organized as follows. 

A preliminary work is carried out in Section~\ref{sec:dessins}.
First, we show that bound (\ref{eq:main-bound}) follows from the 
Riemann--Hurwitz formula for the function $f=P/R$. Then we reduce the 
problem about polynomials to a problem about {\em weighted bicolored 
plane trees}. A weighted bicolored plane tree is a {\em plane tree}\/
(``plane'' means that the cyclic order of branches around each vertex is 
fixed) whose vertices are colored in black and white in such a way that 
the ends of each edge have opposite colors, and whose {\em edges are endowed 
with positive integral weights}. The {\it degree of a vertex} is defined 
as the sum of the weights of the edges incident to this vertex. The sum 
of the weights of all the edges is called the {\em total weight}\/ of the 
tree. We show that a DZ-triple with prescribed factorization pattern
(\ref{eq:P-and-Q}) exists if and only if there exists a weighted bicolored 
plane tree of the total weight $n = \deg P = \deg Q$ having $p$ black 
vertices of degrees $\al_1,\ldots,\al_p$ and $q$ white vertices of 
degrees $\be_1,\ldots,\be_q$. As a corollary, we give a spectacularly 
simple proof of Stothers's 1981 result for the squares and cubes, namely, 
the attainability of the lower bound $\deg (A^3-B^2)\ge k+1$ where 
$\deg A=2k$ and $\deg B=3k$, see Example~\ref{ex:stothers}. The results 
of  Section~\ref{sec:dessins} as well as the framework of the whole paper
are based on the theory of {\em dessins d'enfants}\/ (see, e.\,g., 
Chapter~2 of \cite{LanZvo-04}, or a collection of papers~\cite{Schneps-94}, 
or a recent book~\cite{GirGon-12}). This theory establishes a correspondence
between simple combinatorial objects, graphs drawn on two-dimensional
surfaces, and a vast world of Riemann surfaces, algebraic curves, 
number fields, Galois theory, etc.

In Section~\ref{sec:existence} we prove the existence theorem for 
weighted bicolored plane trees. Namely, we show that a necessary and
sufficient condition for the existence of a tree with the above 
characteristics is inequality (\ref{qwert}). The attainability of
bound \eqref{eq:main-bound} is deduced from this result. 
In Section 4 we establish bound (\ref{no-qwert}) and its attainability 
in the case when inequality \eqref{qwert} is not satisfied. Although 
bounds \eqref{eq:main-bound} and (\ref{no-qwert}) and their attainability 
were proved by U.\,Zannier, we reprove these results here for the sake 
of completeness, and also in order to show how the pictorial language 
clarifies and simplifies the exposition.

In Sections~\ref{sec:unitrees} and \ref{sec:galois} we study DZ-triples 
defined over $\Q$. This case is the most interesting one since by 
specializing $x$ to a rational value one may obtain an important 
information concerning differences of integers with given factorization 
patterns. This subject is actively studied in number-theoretic works: see, 
for example, a recent paper by Beukers and Stewart~\cite{BeuSte-10} (2010) 
and the bibliography therein. Our approach here is based on the following 
corollary of the  theory of  dessins d'enfants which gives a sufficient 
(though not necessary) condition for a DZ-triple to be defined over\/ $\Q$: 
the triple is defined over $\Q$ if there exists {\it exactly one}\/ weighted 
bicolored plane tree with the degrees of black vertices equal to 
$\al_1,\ldots,\al_p$, and the degrees of white vertices equal to 
$\be_1,\ldots,\be_q$. We will call such trees {\it unitrees}.

In Section~\ref{sec:unitrees} we prove the main result of the paper, 
namely, a complete classification of unitrees: see Theorem~\ref{th:main}. 
The formulation of this theorem is rather long so we do not enunciate 
it in the Introduction. We mention only that the class of unitrees 
consists of ten infinite series of trees and ten sporadic trees which 
do not belong to the above series.

While the results of Section~\ref{sec:unitrees} may be considered as 
conclusive, Section~\ref{sec:galois} represents only first steps in a
far-ranging programme of study of the Galois action on weighted plane 
trees and of combinatorial invariants of this action. This approach
permits to find yet more DZ-triples defined over $\Q$ and also to study
DZ-triples over other number fields. Note that essentially all previously 
found examples of DZ-triples over~$\Q$ correspond either to unitrees 
or to the trees constructed in Section~\ref{sec:galois}. We have also
found quite a few new examples of DZ-triples.

Finally, in Section~\ref{sec:further} we mention some further possible
developments of the subject.

\msk

This paper deals only with the combinatorial aspect of the whole 
construction. The computation of the corresponding DZ-triples is 
postponed to a separate publication (see \cite{PakZvo-13}) since 
the techniques used for this purpose are very different form the 
ones used in this paper. In particular, a great deal of symbolic 
computations as well as certain polynomial identities are required. 
For an individual unitree, the computation of the corresponding 
DZ-triple is a difficult task but the verification of the result 
is easy. Indeed, when polynomials $P,Q,R$ (together with their 
appropriate factorizations) are given, it is immediate to observe 
that their coefficients are rational, and the only thing to verify 
is that $R$ is indeed equal to $P-Q$.
The situation becomes significantly more complicated for infinite 
series of trees since in this case the proof may become rather 
elaborate.

\section{From polynomials through Belyi functions to
               weighted trees}
\label{sec:dessins}

\subsection{Function $f=P/R$ and its critical values}

Let $\al,\be\vdash n$ be two partitions of $n$, 
$\al=(\al_1,\ldots,\al_p)$, $\be=(\be_1,\ldots,\be_q)$, 
$\sum_{i=1}^p\al_i=\sum_{j=1}^q\be_j=n$, and let $P,Q\in\C[x]$ be two 
polynomials of degree $n$ with the factorizations
\begin{eqnarray}\label{eq:P,Q}
P(x) \eq \prod_{i=1}^p \, (x-a_i)^{\al_i}\,, \qquad
Q(x) \eq \prod_{j=1}^q \, (x-b_j)^{\be_j}\,.
\end{eqnarray}
We suppose all $a_i,b_j$, $i=1,\ldots,p$, $j=1,\ldots,q$\, to be distinct.
Let the difference $R=P-Q$ have the following factorization:
\begin{eqnarray}\label{eq:R}
R(x) \eq \prod_{k=1}^r\,(x-c_k)^{\ga_k}\,, \qquad 
\deg R \eq \sum_{k=1}^r\,\ga_k\,.
\end{eqnarray}
Our goal is to minimize $\deg R$; obviously, 
\begin{eqnarray}\label{eq:R>r}
\deg R \ge r.
\end{eqnarray}
Consider the following rational function of degree $n$:
$$
f \eq \frac{P}{R}\,;
$$
note that
$$
f-1 \eq \frac{Q}{R}\,.
$$

\begin{definition}[Critical value]\label{def:crit-value}
A point $y\in\Cbar=\C\cup\{\infty\}$ is called {\em critical value}\/
of a rational function $f$ if the equation $f(x)=y$ has multiple roots.
\end{definition}
 
The expressions written above for the function $f=P/Q$ provide us with 
at least three critical values of $f$: 
\begin{itemize}
\item   $y=0$, provided that not all $\al_i$ are equal to 1; 
\item   $y=1$, provided that not all $\be_j$ are equal to 1; and 
\item   $y=\infty$, if only we do not consider the trivial case 
        $\deg R = \deg P - 1$; if $\deg R < \deg P - 1$ then
        $f$ has a multiple pole at infinity.
\end{itemize} 
Denote $y_1,\ldots,y_m$ the other critical values of $f$, if there are any, 
and let $n_l$ be the number of preimages of~$y_l$, $l=1,\ldots,m$; by the 
definition of a critical value, $n_l<n$.

\begin{lemma}[Number of roots of $R$]\label{lem:r}
The number $r$ of distinct roots of the polynomial $R$ is
\begin{eqnarray}\label{eq:r}
r \eq (n+1) - (p+q) + \sum_{l=1}^m\,(n-n_l).
\end{eqnarray}
\end{lemma}

In fact, equality (\ref{eq:r}) is a particular case of the Riemann--Hurwitz 
formula, but for the sake of completeness we give its proof here. 

\begin{figure}[htbp]
\begin{center}
\epsfig{file=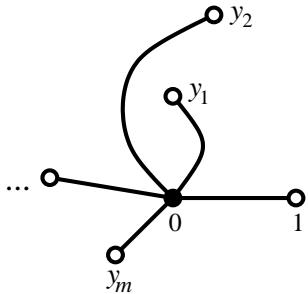,width=4cm}
\caption{\small Star-tree whose vertices are critical values of $f$.}
\label{fig:star}
\end{center}
\end{figure}

\paragraph{Proof.} Let us draw a {\em star-tree}\/ with the center at 0 
and with its rays going to the critical values $1,y_1,\ldots,y_m$, see 
Figure~\ref{fig:star}. 
Considered as a map on the sphere, this tree has $m+2$ vertices, $m+1$ 
edges, and a single outer face with its ``center'' at $\infty$.

Now let us take the preimage of this tree under~$f$. 
We will get a graph drawn on the preimage sphere which has $n\,(m+1)$ 
edges since each edge is ``repeated'' $n$ times in the preimage. 
Its vertices are the preimages of the points $0,1,y_1,\ldots,y_m$, 
so their number is equal to $p+q+\sum_{l=1}^m n_l$. 

What occurs to the faces?

If we puncture at $\infty$ the single open face in the image sphere, 
we get a punctured disk without any ramification points inside.
The only possible unramified covering of a punctured open disk is a
disjoint collection of punctured disks; their number is equal to
the number of poles of $f$, namely, $r+1$ ($r$ roots of~$R$ and $\infty$).
Inserting a point into each puncture we get $r+1$ simply connected
open faces in the preimage sphere. The fact that they are simply
connected implies that the graph drawn on the preimage sphere is
connected. Thus, the preimage of our star-tree is a plane map.
What remains is to apply Euler's formula:
$$
\left( p+q+\sum_{l=1}^m n_l \right) - n\,(m+1) + (r+1) \eq 2,
$$
which leads to \eqref{eq:r}. 
\hfill $\Box$

\bsk

Notice that in order to prove Lemma \ref{lem:r}, instead of the tree of 
Figure~\ref{fig:star}, we could take any other plane map with vertices 
at the critical values (see e.g. the proof of Proposition \ref{prop:thom} below).


\begin{corollary}[Lower bound] We have
\begin{eqnarray}\label{eq:lower-bound}
\deg R \,\ge\, (n+1)-(p+q)\,.
\end{eqnarray}
\end{corollary}
The proof follows from (\ref{eq:r}) and (\ref{eq:R>r}).
\hfill$\Box$

\bsk

Note that $\deg R$ cannot be negative; therefore, when $p+q>n+1$ the 
latter bound cannot be attained. In this case one can attain the bound
$\deg R\ge 0$, that is, the polynomial $R$ can be made equal to a constant.
This situation
is studied in more detail 
in Section~\ref{sec:weak-bound}.

Equation \eqref{eq:r} provides us with guidelines of how to get the
minimum degree of~$R$.

\begin{proposition}[Bound \eqref{eq:lower-bound} attainability]
\label{prop:lower-bound}
Bound \eqref{eq:lower-bound} is attained if and only if the following
conditions are satisfied:
\begin{itemize}
\item   $p+q\le n+1$.
\item   The number $m$ of the critical values of $f$ other than 
        $0,1,\infty$, is equal to zero, so that the sum 
        $\sum_{l=1}^m\,(n-n_l)$ in the right-hand side of~\eqref{eq:r}
        is eliminated altogether. The tree of\/ {\rm Figure~\ref{fig:star}} 
        is then reduced to merely the segment $[0,1]$.
\item   All the roots of $R$ are simple, that is, $\ga_1=\ldots=\ga_r=1$, 
        so that $\deg R=r$.
        Another formulation of the same condition: the partition
        $\ga\vdash n$, $\ga=(\ga_0,\ga_1,\ga_2,\ldots,\ga_r)$, 
        which corresponds to the multiplicities of the poles, has
        the form of a hook\/: 
        $\ga=(n-r,\underbrace{1,1,\ldots,1}_{r\ \mathrm{times}})=(n-r,1^r)$.
\end{itemize}
\end{proposition}

The conditions which imply the existence of such a function~$f$ will be 
obtained in Section~\ref{sec:existence}.

\subsection{Dessins d'enfants and Belyi functions}

Considering rational functions with only three critical values brings us
into the framework of the theory of {\em dessins d'enfants}. Here we
give a brief summary of this theory (only in a planar setting); the 
missing details, proofs, and bibliography can be found, for example, 
in~\cite{LanZvo-04}, Chapter 2.

\begin{definition}[Belyi function]
A rational function $f:\Cbar\to\Cbar$ is called {\em Belyi function}\/
if it does not have critical values outside the set $\{0,1,\infty\}$.
\end{definition}

For such a function, the tree considered in the proof of 
Lemma~\ref{lem:r} is reduced to the segment $[0,1]$. Let us take 
this segment, color the point 0 in black and the point 1 in white, and 
consider the preimage $D=f^{-1}([0,1])$; we will call this preimage a 
{\em dessin}.

\begin{proposition}[Dessin]
The dessin $D=f^{-1}([0,1])$ is a connected graph drawn on the sphere,
and its edges do not intersect outside the vertices. Therefore, $D$ may
also be considered as a plane map. This map has a bipartite structure:
black vertices are preimages of\/ $0$, and white vertices are preimages
of\/~$1$.
\end{proposition}

The degrees of the black vertices are equal to the multiplicities of
the roots of the equation $f(x)=0$, and the degrees of the white ones
are equal to the multiplicities of the roots of the equation $f(x)=1$.
The sum of the degrees in both cases is equal to $n=\deg f$, which is
also the number of edges.

The map $D$ being bipartite, the number of edges surrounding each face
is even. It is convenient, in defining the face degrees, to divide this
number by two.

\begin{definition}[Face degree]\label{def:face-degree}
We say that {\em an edge is incident to a face}\/ if, while remaining
inside this face and making a circuit of it in the positive (trigonometric) 
direction, we follow the edge from its black end toward the white one. 
Thus, only half of the edges surrounding a face are incident to it. 
Moreover, each edge is incident to exactly one face. The {\em degree 
of a face}\/ is equal to the number of edges incident to it.
\end{definition}

According to this definition and to the remarks preceding it, every 
edge is incident to one black vertex, to one white vertex, and to one 
face. The sum of the face degrees is equal to $n=\deg f$.

\begin{proposition}[Faces and poles]
Inside each face there is a single pole of~$f$, and the multiplicity
of this pole is equal to the degree of the face.
\end{proposition}

\begin{definition}[Passport of a dessin]
The triple $\pi=(\al,\be,\ga)$ of partitions $\al,\be,\ga\vdash n$ which 
correspond to the degrees of the black vertices, of the white vertices, 
and of the faces of a dessin, is called a {\em passport}\/ of the dessin.
\end{definition}

\begin{definition}[Combinatorial orbit]\label{def:comb-orbit}
A set of the dessins having the same passport is called a 
{\em combinatorial orbit}\/ corresponding to this passport.
\end{definition}

The construction which associates a map to a Belyi function works also
in the opposite direction. Two bicolored plane maps are {\em isomorphic}\/ 
if there exists an orientation preserving homeomorphism of the sphere 
which transforms one map into the other, respecting the colors of the
vertices. Let $M$ be a bicolored map on the sphere. Then, the sphere 
may be endowed with a complex structure, thus becoming the Riemann 
complex sphere, and a representative of the isomorphism class of $M$ 
can be drawn as a dessin~$D$ obtained via a Belyi function. The following 
statement is a particular case of the classical Riemann's existence 
theorem: 

\begin{proposition}[Existence of Belyi functions]
\label{prop:belyi-exist}
For every bicolored plane map $M$ there exists a dessin $D$ isomorphic 
to $M$, that is, $D=f^{-1}([0,1])$ where $f$ is a Belyi function. 
The function $f=f(x)$ is unique up to a linear fractional transformation 
of the variable $x$.
\end{proposition}

Of course, when we draw a map we do not respect the specific geometric
form of the corresponding dessin. We are content with the fact that such 
a dessin exists.



Now Proposition~\ref{prop:lower-bound} may be reformulated in purely 
combinatorial terms:

\begin{proposition}[Bound \eqref{eq:lower-bound} attainability]
\label{prop:attain}
The lower bound\/ \eqref{eq:lower-bound} is attained if and 
only if there exists a bicolored plane map with the passport 
$\pi=(\al,\be,\ga)$ in which the partitions \linebreak[4] 
$\al=(\al_1,\ldots,\al_p)$ and $\be=(\be_1,\ldots,\be_q)$ are given, 
and $\ga$ has the form $\ga=(n-r,1^r)$ where $1$~is repeated 
$r=(n+1)-(p+q)$ times.
\end{proposition}

In geometric terms, all the faces of our map except the outer one
must be of degree~1. Recall that the {\em number}\/ of faces, which 
is equal to $r+1$, is prescribed by Euler's formula.


\begin{example}[Cubes and squares: a solution]\label{ex:stothers}
Let us look once again at the problem posed by Birch et al.\! in 
\cite{BCHS-65} (see page~\pageref{BCHS}). In order to show that 
if $\deg A=2k$, $\deg B=3k$, and $R=A^3-B^2$, then the lower bound 
$\deg R\ge k+1$ is attained, we must construct a map with the following 
properties: all its black vertices are of degree~3; all its white 
vertices are of degree~2; and all its finite faces are of degree 1. 

\bsk

In order to simplify our pictures we sometimes use the following 
convention.

\begin{convention}[White vertices of degree 2]\label{con:without-white}
When all the white vertices are of degree 2, it is convenient, in order
to simplify a graphical representation of such maps, to draw only
black vertices and to omit the white ones, considering them as being
implicit. In such a picture, a line connecting two black vertices 
contains an invisible white vertex in its middle, and is thus not an 
edge but a union of two edges.
\end{convention}

The construction of the maps we need to solve the above problem 
about $\min\,\deg(A^3-B^2)$ is trivial: first we draw a tree with all 
internal nodes being of degree~3, and then attach loops to its leaves: 
see Figure~\ref{fig:stothers}.

\begin{figure}[htbp]
\begin{center}
\epsfig{file=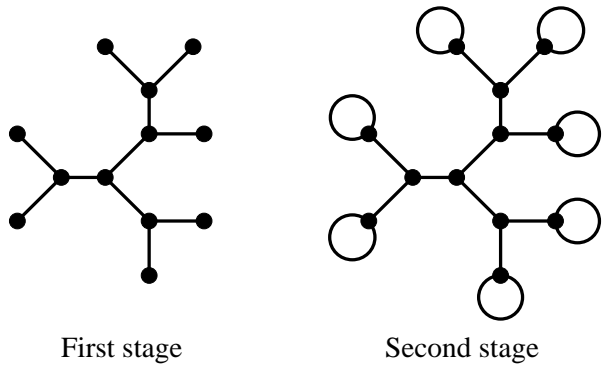,width=7.8cm}
\caption{\small This map solves the problem which remained open for 
16 years: there exist polynomials $A$~and~$B$, $\deg A = 2k$,
$\deg B = 3k$, such that $\deg\,(A^3-B^2)=k+1$.}
\label{fig:stothers}
\end{center}
\end{figure}

We see in this example a remarkable efficiency of the pictorial
representation of problems concerning polynomials. If this
representation was known in~1965, the proof of the conjecture would
have taken 16~minutes instead of 16~years.
\end{example}

\subsection{Number fields}

As it was told in Proposition~\ref{prop:belyi-exist}, a Belyi function 
$f(x)$ corresponding to a dessin is defined up to a linear fractional
transformation of $x$. In this family of equivalent Belyi functions
it is always possible to find one whose coefficients are algebraic
numbers. If we act simultaneously on all the coefficients of such a 
function by an element 
of the {\em absolute Galois group}\/ ${\rm Gal}\,(\Qbar|\Q)$, that is,
by an automorphism $\sigma$ of the field~$\Qbar$ of algebraic numbers, 
or, in other words, if we replace all the coefficients $a_i$ of $f$ by 
their algebraically conjugate numbers $\sigma(a_i)$, we obtain once again 
a Belyi function. Furthermore, one can prove that in such a way the action 
of ${\rm Gal}\,(\Qbar|\Q)$ on Belyi functions descends to an action on 
dessins. There exist many combinatorial invariants of this action, the 
first and the simplest of them being the passport of the dessin. Thus, a 
combinatorial orbit (see Definition~\ref{def:comb-orbit}) may constitute 
a single Galois orbit, or may further split into a union of several 
Galois orbits. Every combinatorial orbit is finite, and therefore every 
Galois orbit is also finite.

One of the most important notions concerning the Galois action on dessins
is that of the field of moduli. 

\begin{construction}[Field of moduli]
Let $D$ be a dessin, and let $\Gamma_D\le {\rm Gal}\,(\Qbar|\Q)$ be 
its stabilizer. Since the orbit $D$ is finite, the group $\Gamma_D$ is a 
subgroup of finite index in ${\rm Gal}\,(\Qbar|\Q)$. Let $H\le\Gamma_D$ 
be the maximal normal subgroup of ${\rm Gal}\,(\Qbar|\Q)$ contained in 
$\Gamma_D$. According to 
the Galois correspondence between subgroups of ${\rm Gal}\,(\Qbar|\Q)$
and algebraic extensions of $\Q$, there exists a number field $K$
corresponding to $H$. This field is called the  {\em field of moduli}\/
of the dessin $D$. By construction, this field is unique: a dessin
cannot have two different fields of moduli.
\end{construction}

Below we list some properties of the fields of moduli.
Let ${\cal D}=\{D_1,\ldots,D_m\}$ be an orbit of the Galois action on 
dessins.
\begin{itemize}
\item   The field of moduli $K$ is the same for all the elements of the orbit.
\item   The degree of $K$ as an extension of $\Q$ is equal to $m=|{\cal D}|$.
\item   The coefficients of Belyi functions corresponding to the
        dessins $D\in{\cal D}$, if they are chosen as algebraic
        numbers, always belong to a finite extension $L$ of~$K$.
\item   The action of the group ${\rm Gal}\,(\Qbar|\Q)$ on the orbit 
        ${\cal D}$ coincides with the action of ${\rm Gal}\,(K|\Q)$.
\item   The action of ${\rm Gal}\,(L|K)$ on Belyi functions may change 
        a position of a dessin $D\in{\cal D}$\/ on the complex sphere but 
        does not change its combinatorial structure; in other words, as a
        map, the dessin in question remains the same.
\end{itemize}

In the absolute majority of cases the situation is much simpler: the
field of moduli of an orbit is the smallest number field to which the 
coefficients of the corresponding Belyi functions belong. However, in 
some specially constructed examples we need a larger field~$L\supset K$ 
to be able to find Belyi functions. There exists a simple {\em sufficient 
condition}\/ which ensures that the coefficients do belong to $K$, 
see \cite{Couveignes-94}: this condition is the existence of a 
{\em bachelor}.

\begin{definition}[Bachelor]\label{def:bachelor}
A {\em bachelor}\/ is a black vertex (a white vertex; a face) such that 
there is no other black vertex (no other white vertex; no other face)
of the same degree.
\end{definition}

\begin{remark}[Positioning of bachelors]\label{rem:bachelor}
If a dessin contains several bachelors then up to three of them can 
be placed at rational points, that is, at points in $\Q\cup\{\infty\}$, 
and this will not prevent the Belyi function for the dessin in question 
to be defined over the field of moduli.
\end{remark}

For the dessins we study in this paper a bachelor always exists: it is 
the outer face (since all the other faces are of degree~1).
Recalling that the degree of $K$ is equal to the size of the orbit
we may conclude the following:

\ssk
\begin{center}
\framebox{\shortstack{If a combinatorial orbit consists of a single
element, \\ then it is also a Galois orbit, and its moduli field is $\Q$.}}
\end{center}
\ssk

Summarizing what was stated above we may affirm the following:

\begin{proposition}[Coefficients in $\Q$]\label{prop:in-Q}
If for a given passport $\pi=(\al,\be,\ga)$, where the partition~$\ga$
is of the form $\ga=(n-r,1^r)$, there exists a unique bicolored 
plane map, then there exists a corresponding Belyi function with rational 
coefficients, and therefore there also exists a DZ-triple with rational 
coefficients.
\end{proposition}

Note that Proposition~\ref{prop:attain}, which concerns the existence, 
is of the ``if and only if'' type, while Proposition~\ref{prop:in-Q}, 
which concerns the definability over $\Q$, provides only an ``if''-type 
condition.

\subsection{How do the weighted trees come in}

\begin{figure}[htbp]
\begin{center}
\epsfig{file=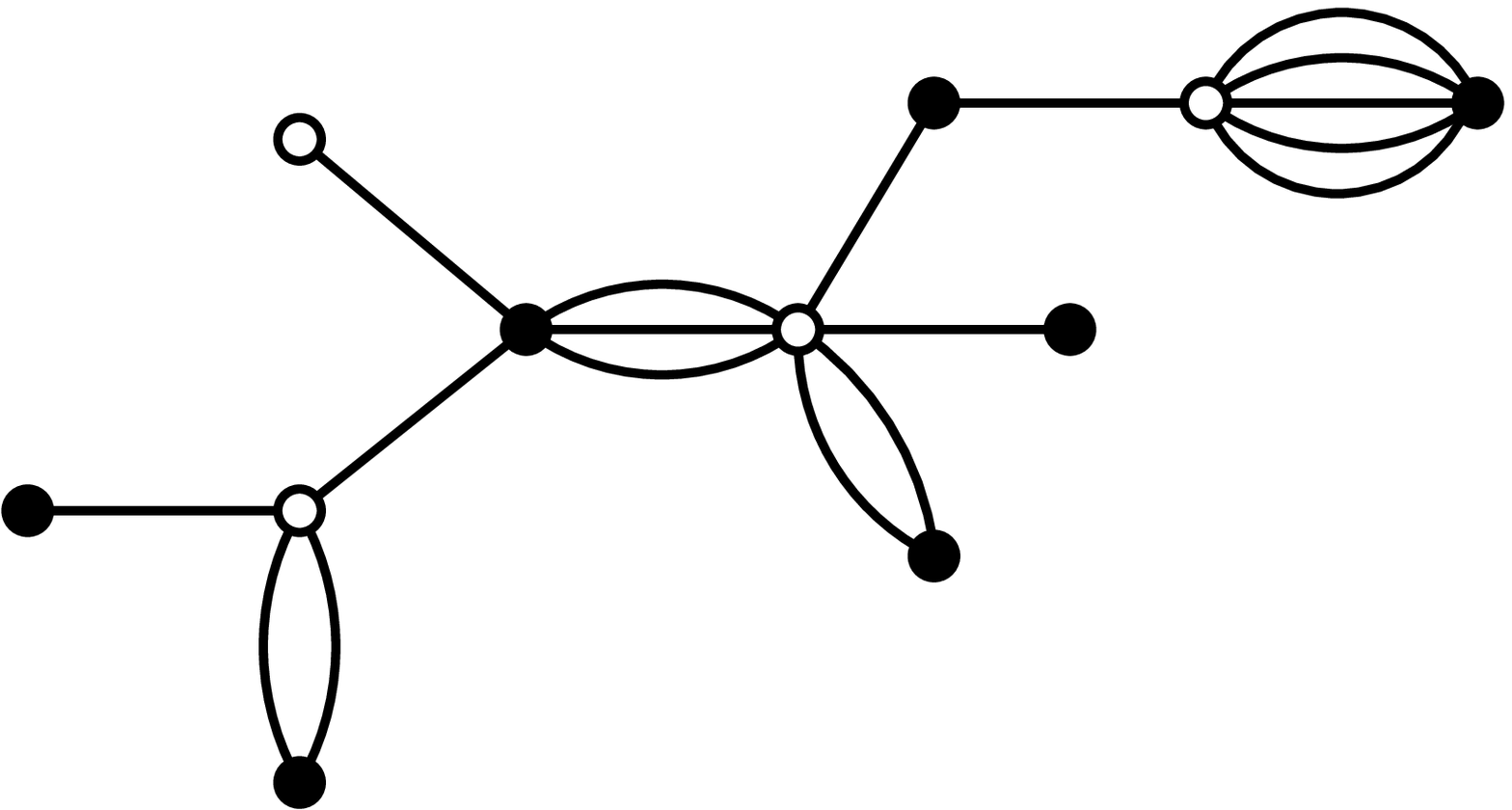,width=6.6cm}
\hspace{0.5cm}
\epsfig{file=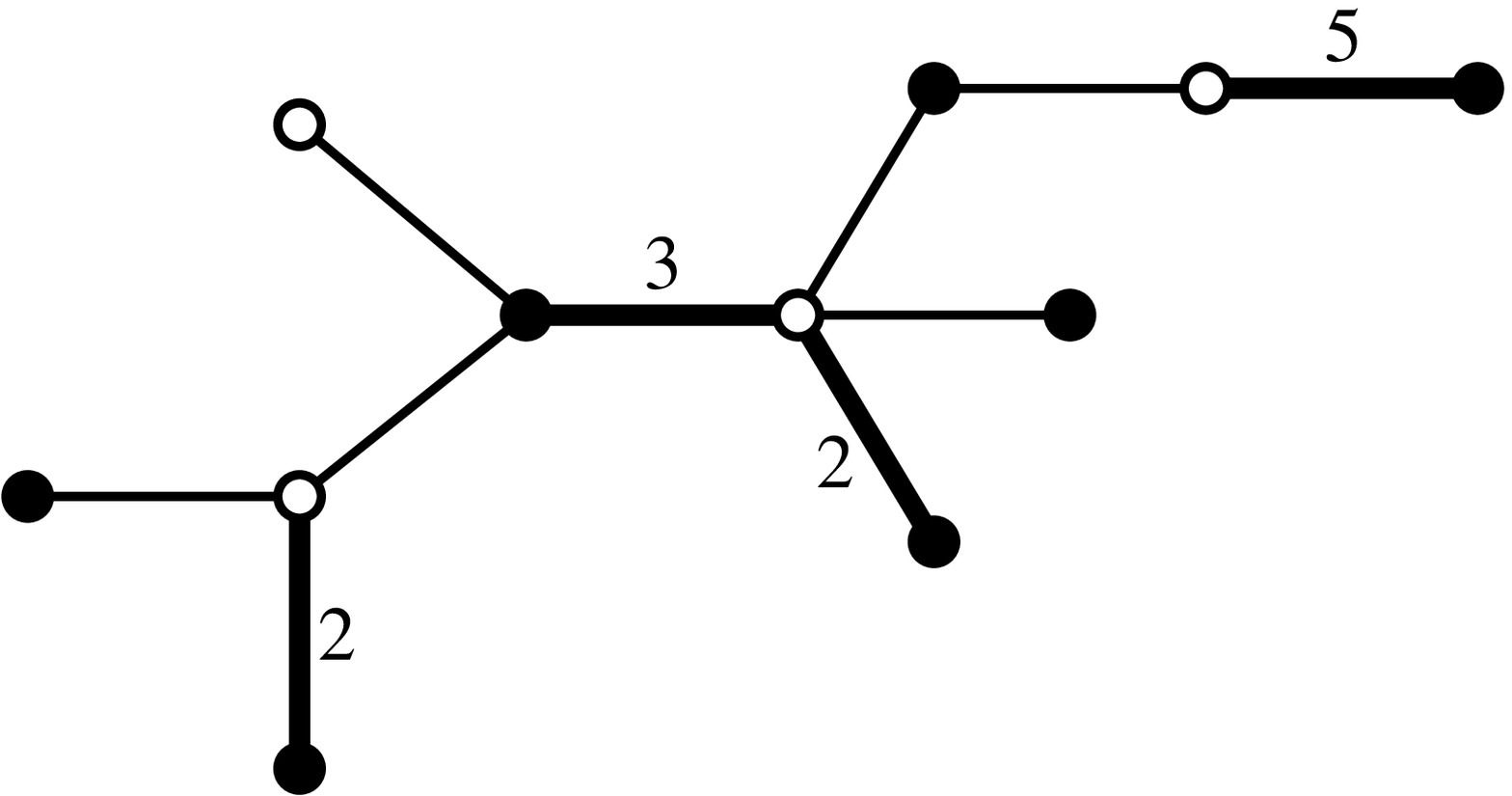,width=5.6cm}
\caption{\small The passage from a map with all its finite faces being 
of degree~1, to a weighted tree. The weights which are not explicitly 
indicated are equal to~1; the edges of the weight bigger than~1 are
drawn thick.}
\label{fig:map->tree}
\end{center}
\end{figure}

Though the weighted trees are, in our opinion, natural and interesting 
objects to be studied for their own sake, in our paper they are used as
a merely technical tool which is easy to manipulate. In 
Figure~\ref{fig:map->tree}, left, it is shown how a typical bicolored 
map whose all finite faces are of degree~1 looks like. (Recall that 
according to Definition~\ref{def:face-degree} a face of degree~1 is 
surrounded by two edges, but only one of these edges is incident to 
the face.) It is convenient to symbolically represent such a map in 
a form of a tree (see Figure~\ref{fig:map->tree}, right) by replacing 
several multiple edges which connect neighboring vertices, by a single 
edge with a weight equal to the number of these multiple edges. In this 
way, the operations of cutting and gluing subtrees, exchanging the weights 
between edges, etc., become easier to implement and to understand.

\begin{definition}[Weighted tree]\label{def:weighted}
A {\em weighted bicolored plane tree}, or a {\em weighted tree}, or just a
{\em tree}\/ for short, is a bicolored plane tree whose edges are endowed 
with positive integral {\em weights}. The sum of the weights of the edges 
of a tree is called the {\em total weight}\/ of the tree. The {\em degree}\/ 
of a vertex is the sum of the weights of the edges incident to this vertex.
The {\em weight distribution}\/ of a weighted tree is a partition 
$\mu\vdash n$, $\mu=(\mu_1,\mu_2,\ldots,\mu_m)$ where $m=p+q-1$ is 
the number of edges, and $\mu_i$, $i=1,\ldots,m$ are the weights 
of the edges. Leaving aside the weights and considering only the
underlying plane tree, we speak of a {\em topological tree}.
Weighted trees whose weight distribution is $\mu=1^n$ will 
be called {\em ordinary trees}. Ordinary trees correspond to 
{\em Shabat polynomials}\/: these are particular cases of Belyi functions, 
with a single pole at infinity.

We call a {\em leaf}\/ a vertex which has only one edge incident to it,
whatever is the weight of this edge. By abuse of language, we will
also call a leaf this edge itself.
\end{definition}

The adjective {\em plane}\/ in this definition means that our trees
are considered not as mere graphs but as plane maps. More precisely,
this means that the cyclic order of branches around each vertex of the 
tree is fixed, and changing this order will in general give a different 
tree. {\em All the trees considered in this paper will be endowed with 
the ``plane'' structure}\/; therefore, the adjective ``plane'' will often 
be omitted.

The adjective {\em plane}\/ in this definition means that our trees
are considered not as mere graphs but as plane maps. More precisely,
this means that the cyclic order of branches around each vertex of the 
tree is fixed, and changing this order will in general give a different 
tree. {\em All the trees considered in this paper will be endowed with 
the ``plane'' structure}\/; therefore, the adjective ``plane'' will often 
be omitted.

\begin{definition}[Isomorphic trees]\label{def:iso}
Two weighted trees are {\em isomorphic}\/ if the underlying bicolored 
plane maps are isomorphic. In other words, they are isomorphic if there 
exists a color-preserving bijection between the vertices of the 
trees and a bijection between the edges which respects the incidence
of edges and vertices, the cyclic order of the edges around each vertex, 
and which also respects the weights of the edges.
\end{definition}

\begin{definition}[Passport of a tree]
\label{def:pass-tree}
The pair $(\al,\be)$ of partitions $\al,\be\vdash n$ of the total
weight~$n$ of a tree, corresponding to the degrees of the black 
vertices and of the white vertices of a weighted tree, is called a 
{\em passport}\/ of this tree. 
\end{definition}

\begin{example}[Tree of Figure \ref{fig:map->tree}]
The total weight of the tree shown in Figure~\ref{fig:map->tree} 
is $n=18$; its passport is $(\al,\be)=(5^22^31^2,7^16^14^11^1)$; 
the face degree distribution is $\ga=10^11^8$, and the weight 
distribution is $\mu=5^13^12^21^6$.
\end{example} 

\begin{remark}[Weighted trees vs.\! ``weighted maps'']
We must not confuse weighted trees with ``weighted maps''. The weighted
tree on the left, and the ``weighted map'' on the right of
Figure~\ref{fig:tree-vs-graph} have the same set of black and white
vertex degrees: $(\al,\be)=(5^13^12^1,5^2)$, but the face degree
partitions of the underlying maps are different: $\ga=4^11^6$ for the 
map represented by the tree, and $\ga=3^12^11^5$ for the map on the right. 
In particular, the corresponding dessins cannot belong to the same 
Galois orbit.

Through the whole paper, we speak exclusively about weighted trees.

\begin{figure}[htbp]
\begin{center}
\epsfig{file=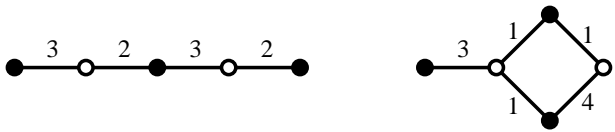,width=8cm}
\end{center}
\caption{\small The weighted tree on the left and the ``weighted map'' 
on the right have the same set of black and white vertex degrees, but 
their face degrees are different.}
\label{fig:tree-vs-graph}
\end{figure}
\end{remark}


Now Propositions~\ref{prop:lower-bound} and \ref{prop:attain} may be 
reformulated as follows:

\begin{theorem}[Lower bound attainability]\label{th:attain}
Let $\al,\be\vdash n$ be two partitions of $n$ having $p$~and~$q$ parts,
respectively.
Then the lower bound\/ \eqref{eq:lower-bound} is attained if and 
only if there exists a weighted tree with the passport $(\al,\be)$.
\end{theorem}

\section{Existence theorem}
\label{sec:existence}

In this section we study the following question: for a given pair 
of partitions $\al,\be\vdash n$, does there exist a weighted tree of the 
total weight $n$ with the passport $(\al,\be)$\,? Equivalently, does there 
exist a rational function with three critical values, and with the 
multiplicities of the preimages of these critical values being, first, 
two given partitions $\al,\be\vdash n$, and then, the third partition being 
equal to $\ga=(n-r,1^r)$\,? 

This question is a particular case of a more general problem of 
{\em realizability}\/ of ramified coverings: does there exist a ramified 
covering of a given Riemann surface with the given ``local data'' (that is, 
with given multiplicities of the preimages of ramification points)\/? 
The problem goes back to the classical paper by Hurwitz~\cite{Hurwitz-1891} 
(1891). Though many particular cases are well studied, the problem
in its full generality remains unsolved. Among numerous publications 
dedicated to the realizability we would like to mention early works by 
Husemoller~\cite{Husemoller-62} (1962) and Thom~\cite{Thom-65} (1965); 
an important paper by Edmonds, Kulkarni, and Stong~\cite{EdKuSt-84} (1984); 
and recent publications \cite{PasPet-09} (2009), \cite{CoPeZa-08} (2008), 
and \cite{Pakovich-09}~(2009).

The main result of this section is the following theorem (recall that
${\rm gcd}(\al,\be)$ denotes the greatest common divisor of the numbers
$\al_1,\ldots,\al_p,\be_1,\ldots,\be_q$):

\begin{theorem}[Realizability of a passport by a tree]\label{th:realize}
Let $\al,\be\vdash n$ be two partitions of $n$, $\al=(\al_1,\ldots,\al_p)$,
$\be=(\be_1,\ldots,\be_q)$, and let\/ ${\rm gcd}\,(\al,\be)=d$. Then
a weighted tree with the passport $(\al,\be)$ exists if and only if
\begin{eqnarray}\label{eq:realize}
p+q\,\le\,\frac{n}{d}+1\,.
\end{eqnarray}
\end{theorem}

By Theorem \ref{th:attain} the attainability of the bound 
\eqref{eq:main-bound} (coinciding with \eqref{eq:lower-bound}) follows 
from this statement. The attainabilty of the bound~\eqref{no-qwert} 
in the case when condition (\ref{eq:realize}) is not satisfied will be 
established in Section~\ref{sec:weak-bound}.

\msk

Theorem \ref{th:realize} and Theorem~\ref{th:weak} below are equivalent 
to the main result (Theorem~1) of Zannier~\cite{Zannier-95}. 
In his paper, Zannier remarks that it would be interesting to apply
the theory of dessins d'enfants to this problem in a more direct way,
and mentions a remark by G.\,Jones that such an approach might produce a
simpler proof. This is indeed the case, as we will see in this section.
Beside that, this theory enables us to find a huge class of DZ-triples 
over $\Q$ (in a way, ``almost all'' of them), see
Section~\ref{sec:unitrees}; and it also gives us a more direct access 
to Galois theory, see Section~\ref{sec:galois}. We have already had a 
first glimpse of the power of the ``dessin method'' in 
Example~\ref{ex:stothers}.

\subsection{Forests}

A {\em forest}\/ is a disjoint union of trees.

\begin{proposition}[Realizability of a passport by a forest]
\label{prop:forest-exist}
Any pair $(\al,\be)$ of partitions of~$n$ can be realized as a passport 
of a forest of weighted trees.
\end{proposition}

\paragraph{Proof.} If there are two equal parts $\al_i=\be_j$ in the 
partitions $\al$ and $\be$, we make a separate edge with the weight 
$s=\al_i=\be_j$ and proceed with the new passport $(\al',\be')$, where 
$\al'$ and $\be'$ are obtained from $\al$ and $\be$ by eliminating 
their parts $\al_i$ and $\be_j$, respectively.

\begin{figure}[htbp]
\begin{center}
\epsfig{file=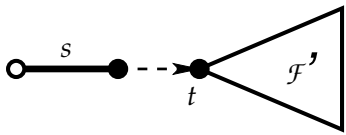,width=4.5cm}
\caption{\small Construction of a forest in the case $\al_i>\be_j$. Here 
$s=\be_j$ and $t=\al_i-\be_j$.} 
\label{fig:forest}
\end{center}
\end{figure}

If there are no equal parts, suppose, without loss of generality, that
there are two parts $\al_i>\be_j$. Then we do the following (see
Figure~\ref{fig:forest}):

\begin{itemize} 
\item[(a)]    make an edge with the weight $s=\be_j$; 
\item[(b)]    consider the new passport $(\al',\be')$ where $\be'$ is 
              obtained from $\be$\/ by eliminating the part $\be_j$, and 
              $\al'$ is obtained from $\al$ by replacing $\al_i$ with 
              $t=\al_i-\be_j$; 
\item[(c)]    construct inductively a forest ${\cal F'}$ of the total 
              weight $n-s$ corresponding to the passport $(\al',\be')$; 
              by definition, this forest must have a black vertex of 
              degree~$t$; 
\item[(d)]    glue the edge of weight $s$ to the forest ${\cal F'}$ by
              fusing two vertices, as is shown in Figure~\ref{fig:forest}, 
              and get a forest $\cal F$ corresponding to $(\al,\be)$
              (since $s+t=\al_i$).
\end{itemize}
The proposition is proved.
\hfill $\Box$

\newpage

\subsection{Stitching several trees to get one: the case ${\rm gcd}(\al,\be)=1$}

\begin{theorem}[Existence]\label{th:existence}
Suppose that ${\rm gcd}\,(\al,\be)=1$. Then the passport $(\al,\be)$
can be realized by a weighted tree if and only if $p+q\le n+1$.
\end{theorem}

\begin{figure}[htbp]
\begin{center}
\epsfig{file=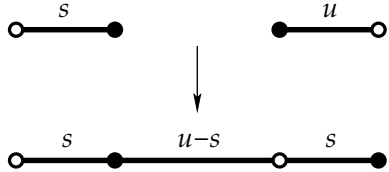,width=5cm}
\caption{\small Stitching two edges.}
\label{fig:gluing-edges}
\end{center}
\end{figure}

\begin{figure}[htbp]
\begin{center}
\epsfig{file=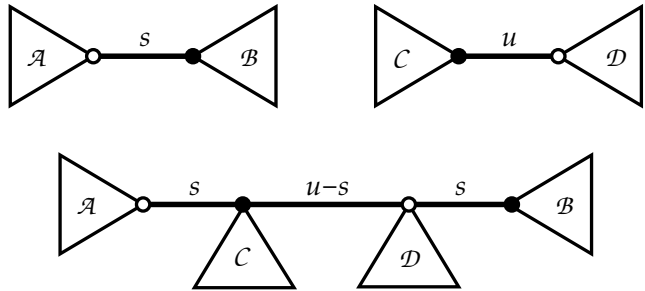,width=8.4cm}
\caption{\small Stitching two trees.}
\label{fig:gluing-trees}
\end{center}
\end{figure}

\paragraph{Proof.} According to Proposition~\ref{prop:forest-exist}
we may suppose that we already have a forest corresponding to the
passport $(\al,\be)$. Now suppose that there are two edges of weights 
$s$ and $u$, $s<u$, which belong to different trees. Then we may stitch
them together by the operation shown in Figure~\ref{fig:gluing-edges}. 
The degrees of the vertices in the new, connected figure are the same 
as in the old, disconnected one. 
Figure~\ref{fig:gluing-trees} shows that the operation works in the same 
way when there are subtrees attached to the ends of the adjoined edges.

We repeat this stitching operation until it becomes impossible to continue. 
The latter may happen in two ways. Either we have got a connected tree, 
and then we are done. Or there are no more edges with different weights 
while the forest remains disconnected. Then, taking into account that 
${\rm gcd}\,(\al,\be)=1$, we conclude that the weights of all edges are 
equal to~1, that is, we have got a forest consisting of $l>1$ ordinary 
trees. In this case, the number of vertices $p+q$ equals $n+l$ and 
therefore is strictly greater than $n+1$, which contradicts the condition 
of the theorem. 
\hfill $\Box$

\bsk

Note that the side edges in Figures \ref{fig:gluing-edges} and 
\ref{fig:gluing-trees} have the same weight $s$. We will use this 
property while performing the operation inverse to stitching, 
namely, the ripping of a connected tree in two, in the proof of 
Proposition~\ref{prop:chain-sts} (see Figure~\ref{fig:unglue}).

\subsection{Non-coprime weights}

Now suppose that ${\rm gcd}\,(\al,\be)=d>1$.

\begin{lemma}[When all weights are multiples of $d>1$]\label{lem:d>1}
The degrees of all vertices of a forest are divisible by $d>1$ if
and only if the weights of all edges are also divisible by $d$.
\end{lemma}

\paragraph{Proof.} In one direction this is evident: the degrees of the
vertices are sums of weights, and therefore, if all the weights are
multiples of $d$, then the same is true for the degrees.

In the opposite direction, if all the vertex degrees are divisible by $d$, 
then it is true, in particular, for the leaves. Cut any leaf off the tree 
to which it belongs, and the statement is reduced to the same one for 
a smaller forest.
\hfill $\Box$

\bsk

Thus, dividing by $d$ all the vertex degrees, that is, all the elements 
of the partitions $\al$ and $\be$, we return to the situation of 
Theorem~\ref{th:existence}, with the same numbers $p$ and $q$, 
and with the total weight equal to $n/d$.
This finishes the proof of Theorem~\ref{th:realize}.
\hfill$\Box$

\bsk

We hope the reader will appreciate the simplicity of the above proof: 
number theorists have been approaching this result for 30 years (1965--1995). 
Once again, the credit goes to the pictorial representation of polynomials 
with the desired properties.

\section{Weak bound}
\label{sec:weak-bound}

When condition (\ref{eq:realize}) of Theorem~\ref{th:realize} is 
satisfied then, according to Theorem~\ref{th:attain}, the main 
bound~(\ref{eq:main-bound}) is attained. If this condition is 
{\em not satisfied}, then the following holds:

\begin{theorem}[Weak bound]\label{th:weak}\,
Let\, ${\rm gcd}\,(\al,\be)=d$, and\, let\, 
${\displaystyle p+q>\frac{n}{d}+1}$.\,
Then\vspace{-2mm}
\begin{eqnarray}\label{eq:weak}
\deg R \,\ge\, \frac{(d-1)\,n}{d},
\end{eqnarray}
and this bound is attained.
\end{theorem}

\subsection{Polynomials and cacti}
\label{sec:cacti}

We will need the following proposition which was proved in 1965 
by Thom~\cite{Thom-65}, and then reproved in many other publications. 
For the reader's convenience we provide a short proof based on 
''Dessins d'enfants" theory following \cite{LanZvo-04}, Corollary~1.6.9.

\begin{proposition}[Realizability of polynomials]\label{prop:thom}
Let $\Lambda=(\la_1,\la_2,\ldots,\la_k)$ be a set of $k\ge 1$ partitions 
$\la_i\vdash n$ of number $n$. Denote by\/ $p_i$ the number of parts of\/ 
$\la_i$, $i=1,2,\ldots,k$. Let $y_1,y_2,\ldots,y_k\in\C$ be arbitrary 
complex numbers. Then a necessary and sufficient condition for the 
existence of a polynomial $T\in\C[x]$ of degree~n, whose all finite
critical values are contained in the set $\{y_1,y_2,\ldots,y_k\}$, with 
the multiplicities of the roots of the equations $T(x)=y_i$ corresponding 
to the partitions $\la_i$, $i=1,2,\ldots,k$, is the following equality:
\vspace{-2mm}
\begin{eqnarray}\label{eq:thom}
\sum_{i=1}^k\,p_i \eq (k-1)\,n+1\,, \qquad \mbox{or, equivalently,} \qquad
\sum_{i=1}^k\,(n-p_i) \eq n-1\,.
\end{eqnarray}
\end{proposition}

\begin{figure}[htbp]
\begin{center}
\epsfig{file=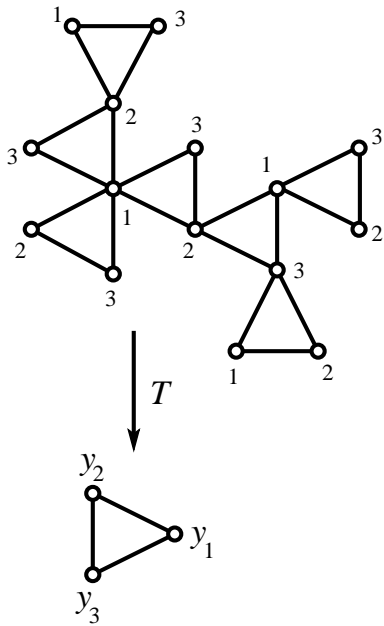,width=5cm}
\end{center}
\caption{\small A cactus. In this example there are three finite 
critical values $y_1,y_2,y_3$; therefore, $D$~is a triangle. A cactus 
is of degree~7, and therefore it contains seven triangles. Vertices
which are preimages of $y_1$, respectively of $y_2$ and of $y_3$, are 
labeled by 1, respectively by 2 and by 3. In this example we have
$\Lambda=(3^12^11^2,2^21^3,2^11^5)$. Namely, $\la_1=3^12^11^2$
shows how many triangles are glued together at vertices labeled by 1,
while partitions $\la_2$ and $\la_3$ correspond to labels 2~and~3.}
\label{fig:cactus}
\end{figure}

\paragraph{Proof.} For purely aesthetic reasons, instead of taking a 
tree with the vertices $y_1,y_2,\ldots,y_k$, as we did in the proof 
of Lemma~\ref{lem:r}, let us take a Jordan curve~$J$ on the 
$y$-plane passing through the points $y_1,y_2,\ldots,y_k$, and let 
$C$ be its preimage $C=T^{-1}(J)$. Then $C$ is a tree-like map 
often called {\em cactus}\/: it does not contain any cycles except $n$ 
``copies'' of $J$ glued together at the vertices which are preimages of 
$y_i$\/; the number of copies of $J$ glued together at a vertex is equal to
the multiplicity of the corresponding root of the equation $T(x)=y_i$\/;
see Figure~\ref{fig:cactus}. Equation (\ref{eq:thom}) may then be 
interpreted as Euler's formula for the cactus since the cactus has 
$\sum_{i=1}^kp_i$ vertices, $kn$ edges, and $n+1$ faces ($n$ copies 
of $J$ and the outer face). 
We leave it to the reader to verify that another proof of the necessity 
of formulas \eqref{eq:thom} can be deduced from the fact that 
the sum $\sum_{i=1}^k(n-p_i)$ in the second equality in (\ref{eq:thom}) 
represents the degree of the derivative~$T'(x)$. 

These observations prove that conditions (\ref{eq:thom}) are necessary. 
Notice that the partition $\la=1^n$ may be eliminated from $\Lambda$ 
(if it belongs to it), and may also be added to it, and this does not 
invalidate equalities (\ref{eq:thom}).

\bsk

The proof that (\ref{eq:thom}) is also sufficient is divided into two parts. 
The first part is purely combinatorial and consists in constructing a 
cactus (at least one) with the vertex degrees corresponding to~$\Lambda$. 
The second part is just a reference to Riemann's existence theorem 
which relates combinatorial  data to the complex structure, as it was 
already the case in Proposition~\ref{prop:belyi-exist}.

\begin{figure}[htbp]
\begin{center}
\epsfig{file=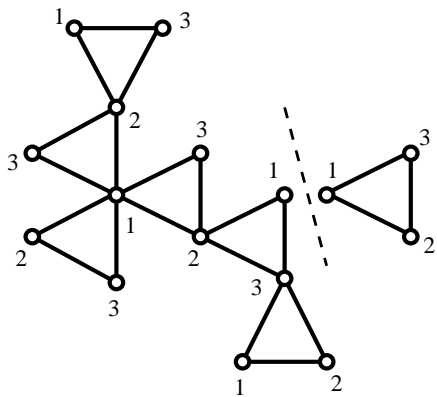,width=5.7cm}
\end{center}
\caption{\small Cutting off a leaf from a cactus. A leaf exists since
every partition in $\Lambda$, except maybe one, contains a part equal
to~1: this is a consequence of (\ref{eq:thom}).}
\label{fig:cut-leaf}
\end{figure}

The proof of the {\em existence of a cactus}\/ in question is similar to
that of Proposition~\ref{prop:forest-exist}; namely, it consists
in ``cutting a leaf''. Here a {\em leaf}\/ means a copy of~$J$ which 
is attached to $C$ at a single vertex (see Figure~\ref{fig:cut-leaf}). 
This cutting operation must be carried out not with the cactus itself
(since it is not yet constructed) but with its passport $\Lambda$:
it is easy to verify that \eqref{eq:thom} implies that all partitions 
$\la_i\in\Lambda$ except maybe one contain a part equal to~1. We eliminate 
these parts, and diminish by~1 a part in the remaining partition. In this 
way we obtain a valid passport $\Lambda'$  of degree $n-1$; then we 
construct inductively a smaller cactus; and then glue to it an $n$th 
copy of~$J$. We leave details to the reader. 
\hfill $\Box$

\bsk

Note that for rational functions, and even for Laurent polynomials, 
a similar statement is not valid, see \cite{Pakovich-09}: conditions 
based on the Euler formula remain necessary but they are no longer
sufficient. See also Example~\ref{ex:non-realize}.

\bigskip

Another approach to the proof of Proposition~\ref{prop:thom} is to use 
an  enumerative formula due to Goulden and Jackson~\cite{GouJac-92} 
which gives the number of cacti corresponding to a given list of 
partitions~$\Lambda$. Let us write a partition $\la\vdash n$ in the 
power notation: $\la=1^{d_1}2^{d_2}\ldots n^{d_n}$ where $d_i$ is 
the number of parts of $\la$ equal to $i$, so that 
$\sum_{i=1}^nd_i=p$ (here $p$ is the total number of parts in $\la$), 
and $\sum_{i=1}^nid_i=n$. Denote
$$
N(\la) \eq \frac{(p-1)!}{d_1!\,d_2!\,\ldots\, d_n!}\,.
$$
Then the following is true:

\begin{proposition}[Enumerative formula]\label{prop:GJ}
For a given $\Lambda=(\la_1,\la_2,\ldots,\la_k)$ satisfying conditions
{\rm (\ref{eq:thom})} we have
\begin{eqnarray}\label{eq:GJ}
\sum \frac{1}{|{\rm Aut}\,(C)|} \eq n^{k-2}\, \prod_{i=1}^k\,N(\la_i)
\end{eqnarray}
where the sum is taken over the cacti $C$ with the passport $\Lambda$,
and\/ $|{\rm Aut}\,(C)|$ is the size of the automorphism group of~$C$.
\end{proposition}

\noindent
Now in order to prove Proposition \ref{prop:thom} it suffices to remark
that the right-hand side of formula (\ref{eq:GJ}) is always positive.
\hfill $\Box$

\begin{remark}[Enumeration of ordinary trees]
Taking $k=2$ in Proposition~\ref{prop:thom} we may put the critical 
values $y_1$ and $y_2$ to 0 and 1, and replace the Jordan curve~$J$ 
passing through these points by the segment $[0,1]$. Then a cactus 
becomes an {\em ordinary}\/ bicolored plane tree with the passport
$\Lambda=(\la_1,\la_2)$. In this case the number of trees (with the 
weights $1/|{\rm Aut}\,(C)|$) is also given by formula (\ref{eq:GJ}). 
This fact will be useful in the future: in order to verify that a given 
ordinary tree is a unitree we can just compute the number given by 
(\ref{eq:GJ}). 

\end{remark}

\subsection{Proof of Theorem \ref{th:weak}}



Consider first the case $d={\rm gcd}(\al,\be)=1$, so that $p+q>n+1$.
In this case the inequality (\ref{eq:weak}) is trivial: it is
reduced to $\deg R\ge 0$. Thus, we only need to prove that this
bound is attained.

We have $n+1\le p+q\le 2n$, therefore $1\le(2n+1)-(p+q)\le n$. 
Let us take an arbitrary partition $\la_3\vdash n$ having 
$(2n+1)-(p+q)$ parts, and also take $\la_1=\al$ and $\la_2=\be$. 
Then for $\Lambda=(\la_1,\la_2,\la_3)$ conditions 
(\ref{eq:thom}) are satisfied. Hence, 
there exists a polynomial $T(x)$ satisfying all the conditions of 
Proposition~\ref{prop:thom}, with three critical values $y_1,y_2,y_3$ 
which may be chosen arbitrarily. Taking $P(x)=T(x)-y_1$ and 
$Q(x)=T(x)-y_2$ we obtain two polynomials which factorize as in 
(\ref{eq:P,Q}) and whose difference is
$$
R(x) \eq P(x)-Q(x) \eq y_2-y_1 \eq {\rm Const}\,.
$$ 
Thus, the obvious lower bound\/ $\deg R\ge 0$ is indeed attained.

\msk

Let us now consider the case ${\rm gcd}\,(\al,\be)=d>1$. In this case
we must prove both the bound~(\ref{eq:weak}) and its attainability.

We have $P=f^d$ and $Q=g^d$. Therefore, $R=f^d-g^d$ factors into $d$ 
factors $f-\zeta g$, where $\zeta$ runs over the $d$-th roots of unity. 
If one of the factors, which we may without loss of generality assume 
to be $f-g$, has degree $<n/d$, then the leading coefficients of 
$f$~and~$g$ coincide. Hence, the leading coefficients of $f$ and 
$\zeta g$ for $\zeta\ne 1$ do not coincide, and all the remaining 
$d-1$ factors have the degree exactly equal to $n/d$. This gives
us the inequality
$$
\deg R \eq \deg\, (f-g)+(d-1)\cdot\frac{n}{d}\, \ge\, \frac{(d-1)\,n}{d}\,.
$$
According to the first part of this proof, the bound $\deg\,(f-g)\ge 0$
is attained, and therefore the bound (\ref{eq:weak}) is also attained. 

This finishes the proof of Theorem \ref{th:weak}.
\hfill$\Box$

\bsk

Notice that the above reasoning may be used for deducing the attainabilty 
of the bound~\eqref{eq:main-bound} in the case when condition 
(\ref{eq:realize}) is satisfied and $d>1$, from the case $d=1$. 
Indeed, it follows from  the case of coprime $\al$ and $\be$ that
$$
\min\deg\, (f-g) \eq \left(\frac{n}{d}+1\right)-(p+q)\,,
$$
and hence 
$$
\min\deg\, (f^d-g^d) \eq \left[\left(\frac{n}{d}+1\right)-(p+q)\right]+
(d-1)\cdot\frac{n}{d} \eq (n+1)-(p+q).
$$

\begin{example}[Weak bound]
Let us take $n=6$, $\al=4^12^1$, and $\be=2^3$, so that $p=2$, $q=3$, 
and $d=2$. Then we have
$$
(n+1)-(p+q)=(6+1)-(2+3)=2
$$ 
but this bound cannot be attained. The correct answer is 
given by Theorem~\ref{th:weak}:
$$
\min\deg R \eq (d-1)\cdot\frac{n}{d} \eq (2-1)\cdot\frac{6}{2} \eq 3\,.
$$
And, indeed, there is only one plane map with two black vertices of 
degrees~4 and~2, respectively, and with three white vertices of degree~2, 
see Figure~\ref{fig:small}. This map has two finite faces, but one of 
them is not of degree~1. The sum of degrees of the finite faces is~3.


\begin{figure}[htbp]
\begin{center}
\epsfig{file=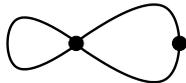,width=2.4cm}
\end{center}
\caption{\small This is the only plane map having the passport 
$(4^12^1,2^3)$. We see that one of the faces is of degree~2.}
\label{fig:small}
\end{figure}

\end{example}

\begin{remark}[Non-realizable planar data]\label{ex:non-realize}
Let us take $\al$ and $\be$ such that ${\rm gcd}\,(\al,\be)=d>1$
and ${\displaystyle \frac{n}{d}+1<p+q\le n+1}$. Let us also take
$r=(n+1)-(p+q)$ and $\ga=(n-r,1^r)$. 
Then the passport $\pi=(\al,\be,\ga)$ satisfies the Euler relation: 
there are $p+q$ vertices, $n$ edges, and $r+1$ faces, so that
$$
(p+q)-n+(r+1) \eq (p+q)-n+[(n+1)-(p+q)+1] \eq 2.
$$
However, a plane map with these data does not exist.

The principal vocation of this paper is to use combinatorics for the
study of polynomials. But here in this particular example we, essentially, 
deduce a non-trivial statement about plane maps from a trivial property of
polynomials. Namely, we deduce the non-existence of certain maps from the 
fact that the degree of a polynomial cannot be negative.
\end{remark}

\section{Classification of unitrees}
\label{sec:unitrees}

This section contains the main results of the paper: here we give a 
complete classification of the passports satisfying the following 
property: a weighted bicolored plane tree having this passport is unique. 
As we have explained before, in Proposition~\ref{prop:in-Q}, Belyi 
functions corresponding to such trees are defined over~$\Q$.

Ordinary unitrees were classified by Adrianov in 1989. However, his
initial proof was never published since it looked too cumbersome.
Then, in 1992, appeared the paper \cite{GouJac-92} by Goulden and Jackson 
with the enumerative formula (\ref{eq:GJ}), which opened a possibility
for another proof, by carefully analizing this formula and looking for 
the cases in which it gives a number~$\le 1$ (recall that (\ref{eq:GJ}) 
counts each tree $C$ with the weight $1/|{\rm Aut}(C)|$). 

Our situation is more difficult for two reasons. First, we deal with 
weighted trees; and, second, we don't have an enumerative formula at 
our disposal. For these (and many other) reasons such a formula would 
be very welcome. To be more specific, we need a formula which would 
give us, in an explicit way, the number of the weighted bicolored plane 
trees {\em corresponding to a given passport}. An additional difficulty 
here ensues from the fact that the same passport may correspond not only 
to trees but also to forests.

\begin{assumption}[Passports from now on] 
\label{ass:simple-passport}
In the remaining part of the paper we will consider only the passports 
$(\al,\be)$ such that ${\rm gcd}\,(\al,\be)=1$ and $p+q\le n+1$. 
\end{assumption}

According to Lemma~\ref{lem:d>1}, the case ${\rm gcd}\,(\al,\be)>1$ is
reduced to this one. Indeed, starting from a tree $\cal T$ with 
${\rm gcd}\,(\al,\be)=d>1$ we can obtain a tree $\widetilde {\cal T}$ 
with ${\rm gcd}\,(\widetilde \al,\widetilde \be)=1$ by dividing the 
weights of all edges of $\cal T$ by $d$, and it is easy to see that 
$\cal T$ is a unitree if and only if $\widetilde {\cal T}$  is a unitree.

Recall that the face partition $\ga$ is defined by $(\al,\be)$ and is 
always equal to $\ga=(n-r,1^r)$ where $r=(n+1)-(p+q)$.

\begin{definition}[Unitree]
A weighted bicolored plane tree such that there is no other tree with the
same passport is called a {\em unitree}.
\end{definition}

\subsection{Statement of the main result}

\begin{definition}[Diameter]
The {\em diameter}\/ of a tree is the length of the longest path in
this tree.
\end{definition}

\begin{figure}[!h]
\begin{center}
\epsfig{file=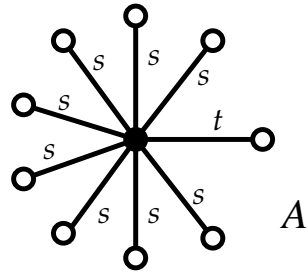,width=4cm}
\caption{\small Series $A$: stars. The edge of the weight $s$ is repeated
$k\ge 0$ times.}
\label{fig:stars-A}
\end{center}
\end{figure}

\begin{figure}[!h]
\begin{center}
\epsfig{file=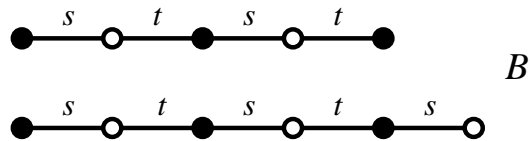,width=6.9cm}
\caption{\small Series $B$: periodic chains of an arbitrary length. We 
distinguish the chains of even and odd length since they have passports 
of different type.}
\label{fig:chains-B}
\end{center}
\end{figure}

\begin{figure}[!h]
\begin{center}
\epsfig{file=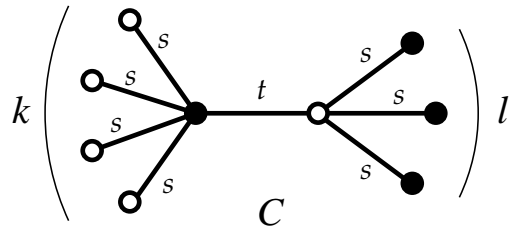,width=6.6cm}
\caption{\small Series $C$: brushes of diameter 3. Here $k,l\ge 1$.}
\label{fig:brushes-C}
\end{center}
\end{figure}

\begin{figure}[!h]
\begin{center}
\epsfig{file=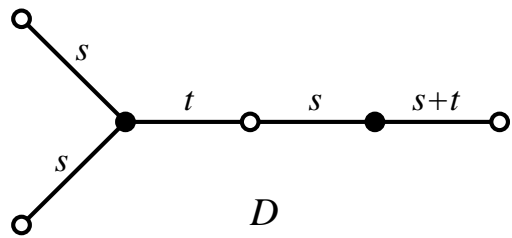,width=6.6cm}
\caption{\small Series $D$: brushes of diameter 4. There are exactly two 
leaves of weight~$s$ and exactly one leaf of weight~$s+t$.}
\label{fig:brushes-D}
\end{center}
\end{figure}

\begin{figure}[htbp]
\begin{center}
\epsfig{file=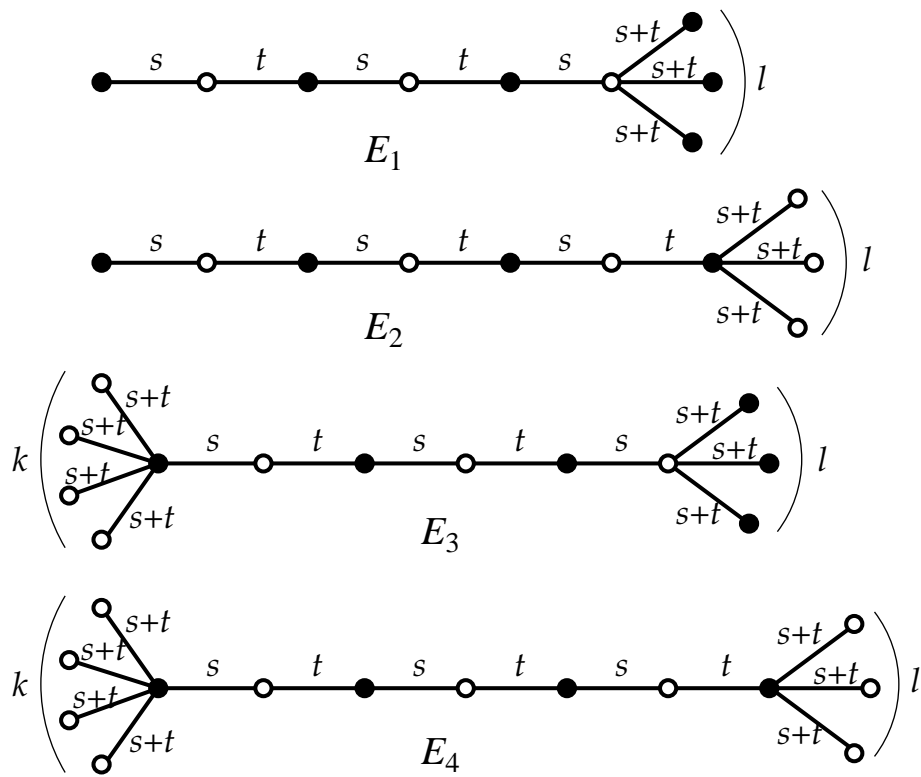,width=12cm}
\caption{\small Series $E$: brushes of an arbitrary length. If there 
is a leaf of weight $s$, it is ``solitary'' on one of the ends of the 
brush; otherwise, all the leaves are of the weight $s+t$. The parameters 
$k,l\ge 1$.} 
\label{fig:brushes-E}
\end{center}
\end{figure}

\begin{figure}[!h]
\begin{center}
\epsfig{file=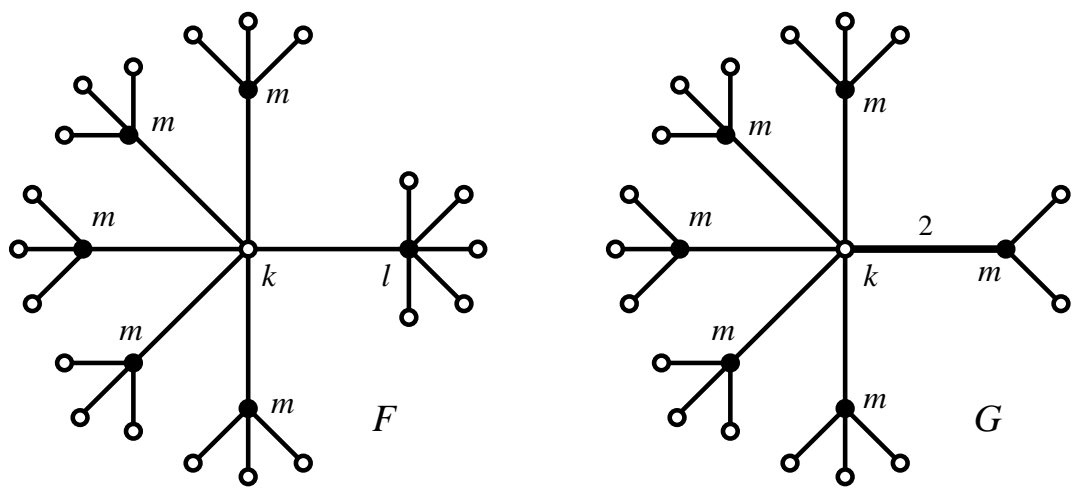,width=14cm}
\caption{\small Two series of unitrees of diameter 4. In the trees of
the series $F$ all the edges are of weight~1; the degrees of vertices 
(except the leaves) are indicated. In the trees of the series $G$, there 
is exactly one edge of weight~2, all the other edges being of weight~1; 
note that this time the degrees of the black vertices are all equal.}
\label{fig:diameter-4}
\end{center}
\end{figure}

\begin{figure}[htbp]
\begin{center}
\epsfig{file=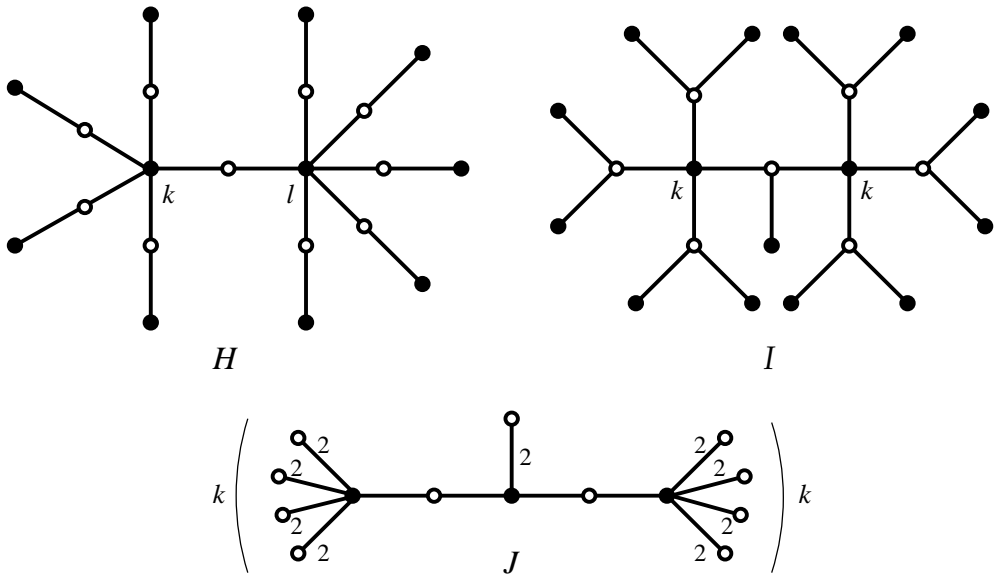,width=13cm}
\caption{\small Three series of unitrees of diameter 6. In $H$ and $I$,
all edges are of weight~1. In $H$ the black vertices which are not 
leaves are of degrees $k$ and $l$ which may be non-equal; in~$I$ they 
are of the same degree~$k$. In $H$ all white vertices are of degree~2; 
in~$I$, they are all of degree~3. In $J$, the number of leaves of the 
weight~2 on the left and on the right is the same.}
\label{fig:diameter-6}
\end{center}
\end{figure}

\begin{figure}[htbp]
\begin{center}
\epsfig{file=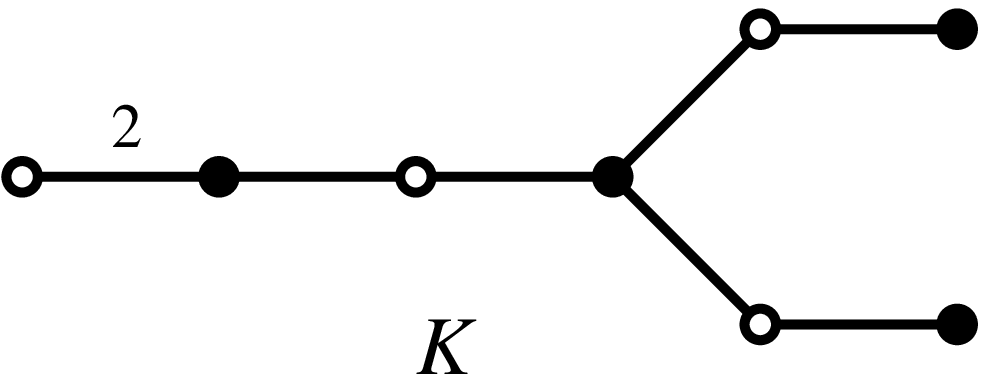,width=5.4cm}
\caption{\small A sporadic unitree of diameter 5.}
\label{fig:sporadic-5}
\end{center}
\end{figure}

\begin{figure}[htbp]
\begin{center}
\epsfig{file=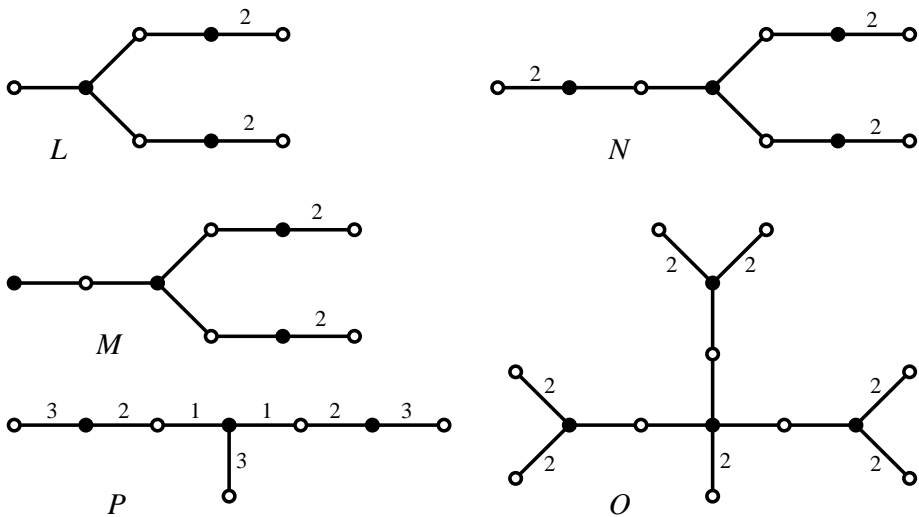,width=12cm}
\caption{\small Five sporadic unitrees of diameter 6.}
\label{fig:sporadic-6}
\end{center}
\end{figure}

\begin{figure}[htbp]
\begin{center}
\epsfig{file=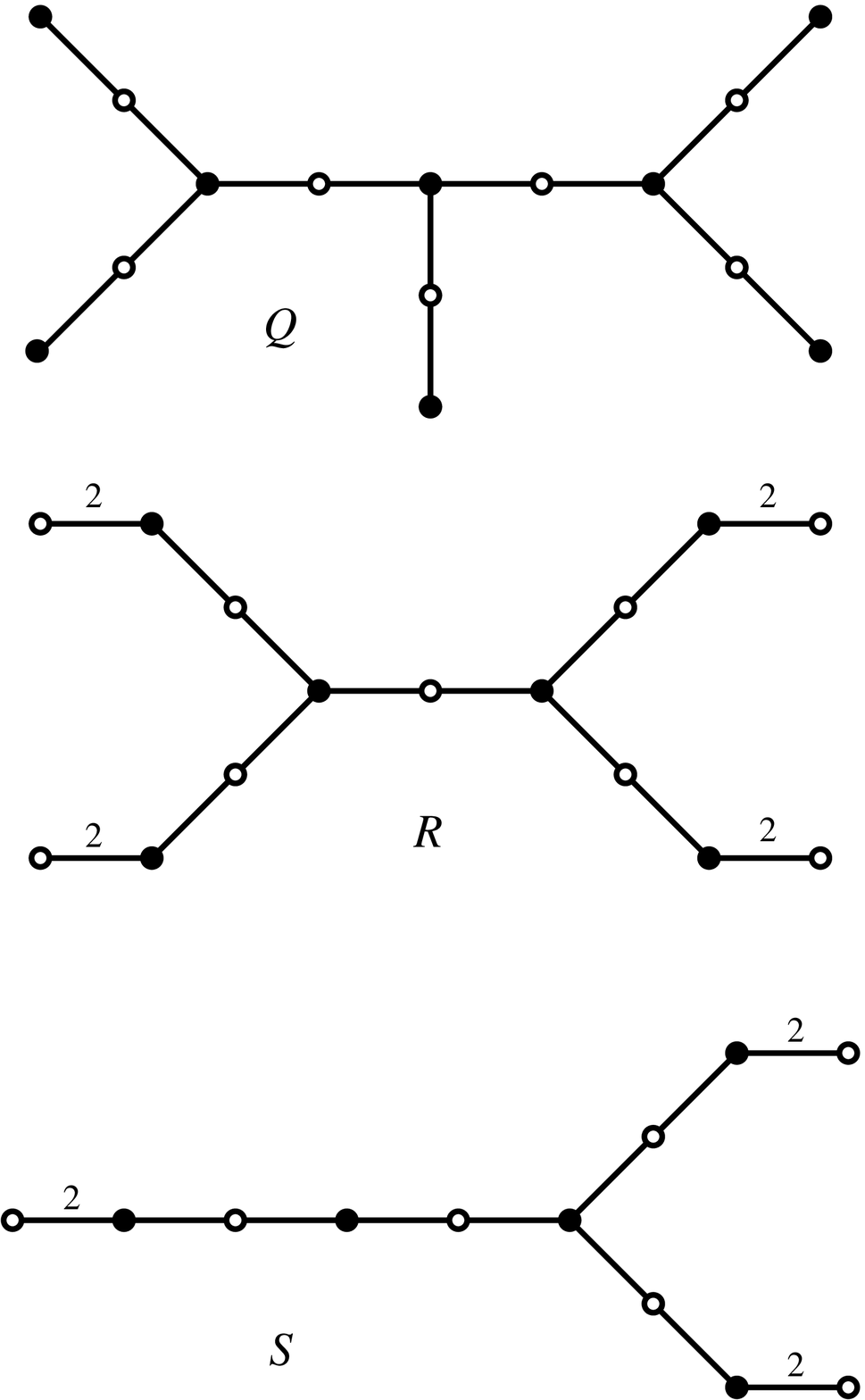,width=8.5cm}
\caption{\small Three sporadic unitrees of diameter 8.}
\label{fig:sporadic-8}
\end{center}
\end{figure}

\begin{figure}[htbp]
\begin{center}
\epsfig{file=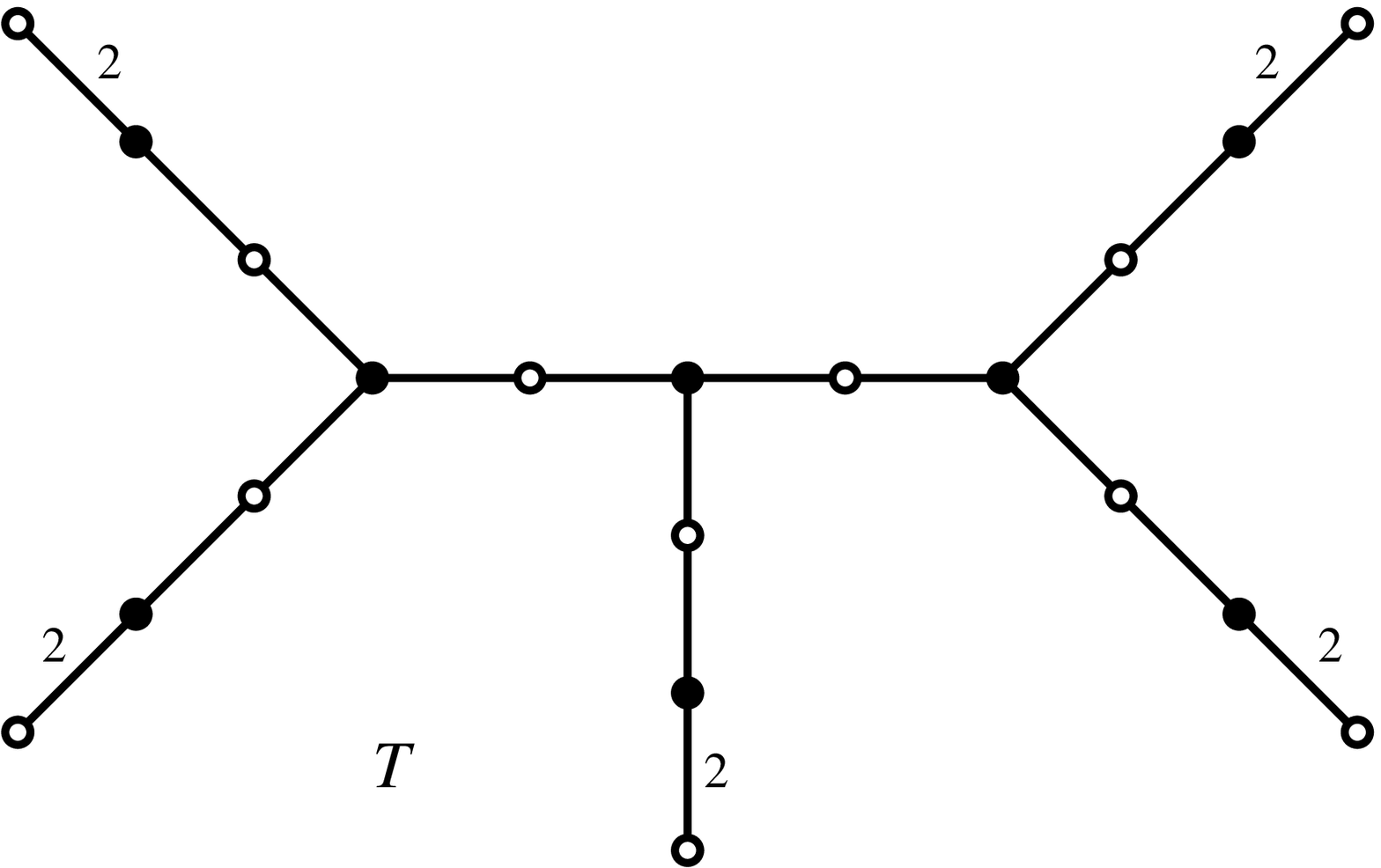,width=10cm}
\caption{\small A sporadic unitree of diameter 10.}
\label{fig:sporadic-10}
\end{center}
\end{figure}

The classification of unitrees is summarized in the following theorem:

\begin{theorem}[Complete list of unitrees]\label{th:main}
Up to an exchange of black and white and to a multiplication of all the
weights by $d>1$, the complete list of unitrees consists of the following\/ 
$20$ cases:
\begin{itemize}
\item   Five infinite series A, B, C, D, E of trees shown in 
        {\rm Figures \ref{fig:stars-A}, \ref{fig:chains-B}, 
        \ref{fig:brushes-C}, \ref{fig:brushes-D}, and \ref{fig:brushes-E}},
        involving two integral weight parameters $s$ and $t$ which are 
        supposed to be coprime $($thus either $s\ne t$, or $s=t=1)$. 
        Note that
        \begin{itemize}
        \item   for the diameter $\ge 5$, only the trees of the types~$B$ 
                and~$E$ exist;
        \item   for the diameter\/~$4$, the trees of types $B$, $D$, and $E$ 
                exist;
        \item   for the diameter\/~$3$, the trees of types $B$, $C$, and $E$ 
                exist.
        \end{itemize}
\item   Five infinite series F, G, H, I, J shown 
        in\/ {\rm Figures \ref{fig:diameter-4}} 
        and\/~{\rm \ref{fig:diameter-6}}.
\item   Ten sporadic trees K, L, M, N, O, P, Q, R, S, T shown in\/ 
        {\rm Figures \ref{fig:sporadic-5}, \ref{fig:sporadic-6},
        \ref{fig:sporadic-8}}, and\/ {\rm \ref{fig:sporadic-10}}.
\end{itemize}
\end{theorem}

\begin{remark}[Non-disjoint]
The above series are not disjoint. For example, the trees of the 
series~$C$ with $k=l=1$ also belong to the series~$B$. If $s>t$,
$C$~becomes a part of $E_3$, up to a renaming of variables; etc.
\end{remark}

\begin{remark}[Adrianov's list]
The list of ordinary unitrees compiled by Adrianov in 1989 consists of
the following cases: the series $A$, $B$, and $C$ with $s=t=1$; 
the series $F$, $H$, and~$I$; and the sporadic tree $Q$. 
\end{remark}

\begin{figure}[htbp]
\begin{center}
\epsfig{file=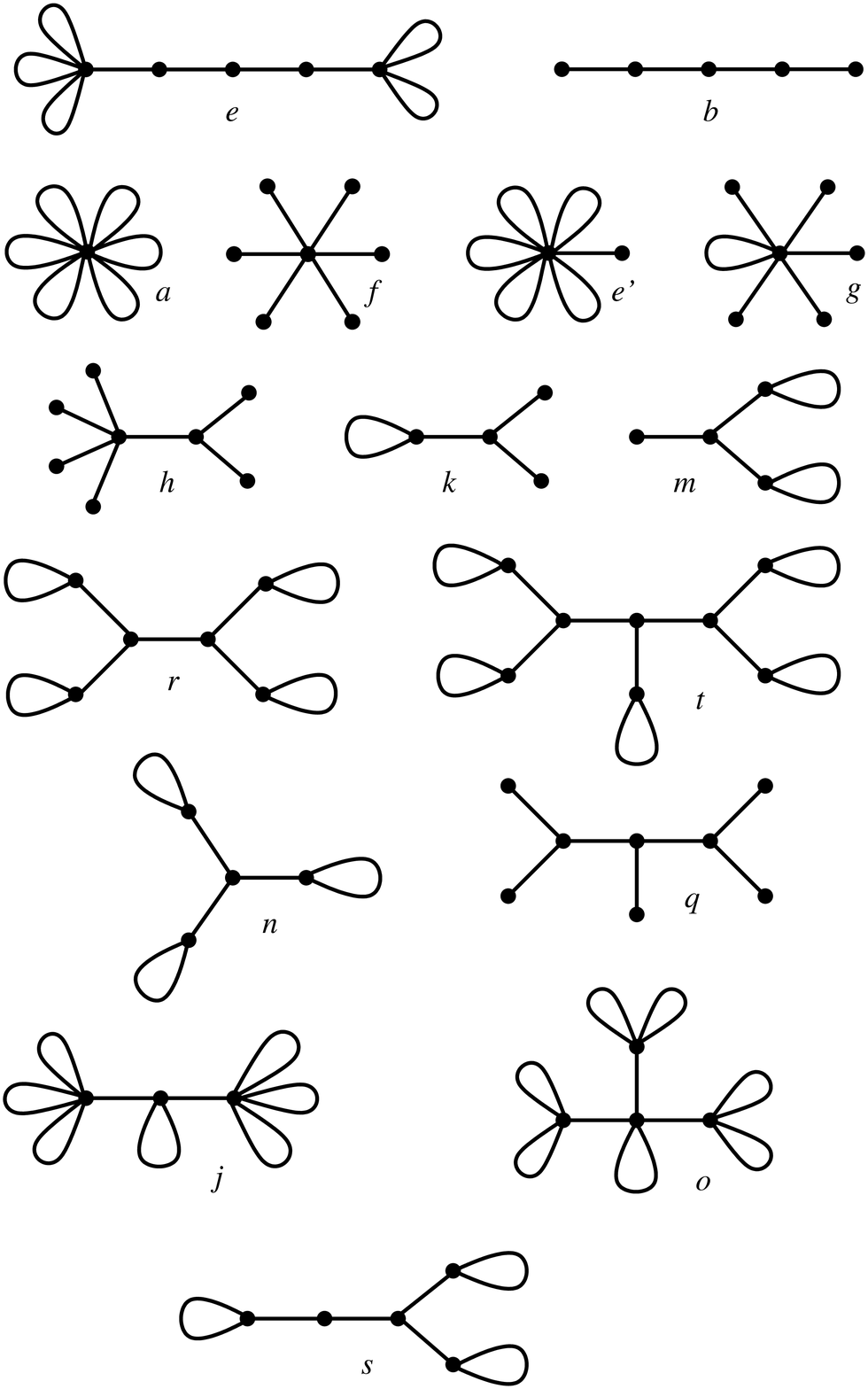,width=12cm}
\caption{\small Unimaps. Small letters correspond to the capital letters
by which we have previously denoted the unitrees; for example, the series
$e$ here is a particular case of the series $E$ (when $s=t=1$). Note that 
$a$, $b$, $e'$, $f$, and $g$ are particular cases of $e$. Note also that 
the series $h$ and $j$ and the sporadic unimaps $k$, $m$, $n$, $o$, $q$, 
$r$, $s$ are not ``particular cases'' but just coincide with $H$, $J$, $K$, 
etc., respectively.}
\label{fig:trees2maps}
\end{center}
\end{figure}

\begin{remark}[White vertices of degree 2]
Notice that quite a few of our trees have all their white vertices being
of degree~2, and thus, according to Convention~\ref{con:without-white},
we can make these vertices implicit and draw the pictures as usual plane
maps. The corresponding maps are shown in Figure~\ref{fig:trees2maps}. 
\end{remark}

The strategy of the proof of Theorem~\ref{th:main} is as follows.
We propose various transformations of trees changing the trees 
themselves while preserving their passports: this is a way to show
that the combinatorial orbit of a given tree consists of more that
one element. The trees which survive such a surgery are (a)~those
to which the transformation in question cannot be applied, and (b)~those 
for which the transformed tree turns out to be isomorphic to the 
initial one. In this way we gradually eliminate all the trees which 
are not unitrees. Then, at the final stage, we show that all the trees 
which have passed through all the sieves are indeed unitrees.

The proof ends on page~\pageref{end-proof}.

\vspace{\stretch{1}}

\vspace*{\stretch{1}}

\vspace{\stretch{1}}

\newpage

\newpage

\begin{definition}[Rooted tree]\label{def:rooted}
A tree with a distinguished leaf edge is called {\em rooted tree},
and the distinguished edge itself is called its {\em root}. Two
rooted trees are isomorphic if there exists an isomorphism which
sends the root of one of the trees to the root of the other one.
\end{definition}

\begin{definition}[Branch]\label{def:branch}
Let a vertex $v$ of a tree $\cal T$ be given. Then a {\em branch}\/ 
of~$\cal T$ attached to~$v$ is a rooted tree which is a subtree 
of~$\cal T$ containing a single edge incident to $v$. This edge is 
the root of the branch.
\end{definition}

The following characterization of unitrees eliminates a vast amount 
of possibilities.

\begin{lemma}[Branches of a unitree]\label{lem:branches}
All branches going out of a vertex of a unitree, except maybe one, 
are isomorphic as rooted trees. This property must be true for every 
vertex of a unitree.
\end{lemma}

\paragraph{Proof.} Let us call a vertex {\em central}\/ if it is obtained 
by the following procedure. We cut off all the leaves of the tree; 
then do the same with the remaining smaller tree, then again, etc. 
In the end, what remains is either a single vertex, or an edge. 
In the first case, there is a single central vertex; in the second case, 
there are are two of them, one black and one white. Obviously, 
{\em any isomorphism of trees sends a central vertex to a central one}, 
and if there are two central vertices, it sends the black central vertex 
to the black one, and the white, to the white one. On the other hand, 
according to Definition~\ref{def:iso}, any isomorphism which sends
a vertex~$v$ to itself must preserve the cyclic order of the branches
attached to~$v$. Thus, the property affirmed in this lemma, namely,
that all the branches except maybe one are isomorphic, is valid for a 
central vertex or vertices, since otherwise an operation of 
{\em exchanging}\/ of two non-isomorphic branches attached to the central 
vertex would change the cyclic order of branches around this vertex.

Further, observe that the operation of exchanging of two non-isomorphic
branches attached to a vertex of a {\em rooted}\/ tree always changes 
this tree unless  all the 
branches attached to this veretx, except maybe the 
branch containing the root, are isomorphic.  Indeed, the introduction of 
a root makes a cyclic order on branches around any vertex into a 
{\em linear order}\/ on the branches incident to it and not containing 
the root. The only possibility to make this linear order invariant under 
the operation of exchanging the branches is to make them all equal.

Now, if a vertex $v$ of a tree $\cal T$ is not central then it belongs 
to a branch $\cal V$ attached to a central vertex. This branch itself is 
a rooted tree and, in case if the condition of the lemma is not satisfied, 
the operation of exchanging of the branches changes  $\cal V$.
However, changing the branch $\cal V$ would mean changing a single branch of
$\cal T$ attached to its central vertex. This would make $\cal T$ not
isomorphic to itself. Thus, in this case $\cal T$ would not be a unitree. 
\hfill$\Box$

\subsection{Weight distribution}
\label{sec:weight-exchange}

Sometimes, we can change not the topology of the tree but the
distribution of its weights, while remaining in the same combinatorial
orbit. Let us first formulate a statement which is entirely obvious:

\begin{lemma}[Weight distribution]\label{lem:wd}
If a passport $(\al,\be)$ corresponds to a unitree then the corresponding
weight distribution $\mu$ {\rm (see Definition~\ref{def:weighted})} 
is determined by $\al$ and $\be$ in a unique way.
\end{lemma}

\begin{figure}[htbp]
\begin{center}
\epsfig{file=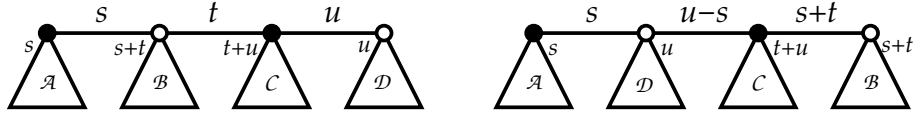,width=12cm}
\caption{\small Weight exchange. If $s<u$ but $u\ne s+t$ then exactly 
two parts of the weight distribution $\mu$ have changed.}
\label{fig:weight-exchange}
\end{center}
\end{figure}

\begin{lemma}[Condition on weights]\label{lem:weight-exchange}
Let $s,t,u$ be the weights of three successive edges of a unitree, as
in\/ {\rm Figure~\ref{fig:weight-exchange}, left}. If $s\le u$, then 
either $u=s$, or $u=s+t$. Similarly, if $s\geq u$, then either $u=s$, 
or $u=s-t$. 
\end{lemma}

\vspace{-5mm}

\paragraph{Proof.} 
Rotating if necessary the tree under consideration, without loss of 
generality we may assume that $s\le u$. If $s<u$ then we can construct 
the tree shown in Figure~\ref{fig:weight-exchange}, right, replacing 
the weight~$u$ with~$s+t$, the weight~$t$ with~$u-s$, and exchanging 
the places of the subtrees $\cal B$~and~$\cal D$. We see that the vertex 
degrees of the new tree are the same as in the initial one, while the 
weights of two edges have changed, unless $u=s+t$. Thus, if $s<u$ but 
$u\ne s+t$, then exactly two parts of the weight distribution $\mu$ 
have changed, which contradicts Lemma~\ref{lem:wd}.
\hfill $\Box$

\begin{corollary}[Adjacent edges]\label{cor:adj-edges}
If in a unitree there are two adjacent edges of the same weight $s$, 
and at least one of them is not a leaf, then $s=1$.
\end{corollary}

\vspace*{-5mm}

\begin{figure}[htbp]
\begin{center}
\epsfig{file=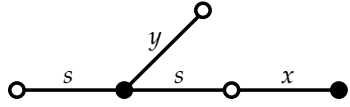,width=4.5cm}
\caption{\small Two adjacent edges of the same weight $s$; one of them is
not a leaf.}
\label{fig:adj-edges}
\end{center}
\end{figure}

\paragraph{Proof.} An edge which is not a leaf must be the middle edge of
a path of length~3, see Figure~\ref{fig:adj-edges}. According to 
Lemma~\ref{lem:weight-exchange} we have either $x=s$ or $x=2s$. 
If there is an edge of weight~$y$ attached to the middle edge of
the path, like in Figure~\ref{fig:adj-edges}, then we have either $y=x$ 
or $y=x+s$, so the possible values for $y$ are $s$, $2s$, or $3s$. 
Dealing in the same way with the other edges of the tree we see that 
the weights of all of them are multiples of $s$. According to 
Assumption \ref{ass:simple-passport} this means that $s=1$. 
\hfill$\Box$

\vspace{-5mm}

\begin{proposition}[Path $s$, $t$, $s$]\label{prop:chain-sts}
Suppose that a unitree contains a path of three successive edges having
the weights $s$, $t$, $s$. Then the only possible weights for all the
edges of this tree are $s$, $t$, or $s+t$.
\end{proposition}

\begin{figure}[htbp]
\begin{center}
\epsfig{file=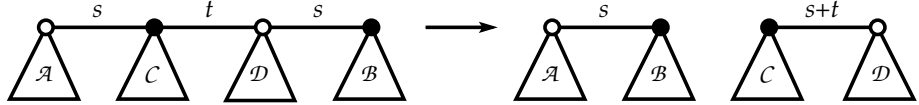,width=12cm}
\caption{\small ``Ripping'' a tree: an operation inverse to that of 
the proof of Theorem~\ref{th:existence}.}
\label{fig:unglue}
\end{center}
\end{figure}

\paragraph{Proof.} Let us make an operation inverse to the one used in
the proof of Theorem~\ref{th:existence}, that is, ``rip'' the tree 
along the edge of the weight $t$, as in Figure~\ref{fig:unglue}.

\noindent
Now suppose that in one of the subtrees ${\cal A,B,C,D}$ there exists
an edge of a weight $x\ne s,t,s+t$, and try to stitch the two trees 
together in a different way.
\begin{enumerate}
\item   Suppose that in one of the subtrees $\cal A$ or $\cal B$ there 
        exists an edge of a weight $x\ne s+t$. Stitch it to the edge 
        of the weight $s+t$ by the procedure explained in the proof
        of Theorem~\ref{th:existence}.
        \begin{enumerate}
        \item   If $x<s+t$ then the weights of the four edges participating
                in the operations of ripping and stitching, instead of being 
                $s,t,s,x$ become
                $s+t-x,x,x,s$. Removing from the two sets the coinciding
                elements $s$ and $x$, we get, on the one hand, $s,t$,
                and, on the other hand, $s+t-x,x$. These sets coincide 
                only when $x=s$ or $x=t$.
        \item   If $x>s+t$ then, instead of $s,t,s,x$, we obtain
                $x-s-t,s+t,s+t,s$. These two sets cannot coincide at all
                since $x$ is greater than every term in the second set.
        \end{enumerate}
\item   Suppose that in one of the subtrees $\cal C$ or $\cal D$ there 
        exists an edge of a weight $x\ne s$. Stitch it to the edge of 
        the weight $s$ by the same procedure as above.
        \begin{enumerate}
        \item   If $x>s$ then, instead of $s,t,s,x$, we get $x-s,s,s,s+t$. 
                Removing $s,s$ from both sets we get, on the one hand, 
                $t,x$, and on the other hand, $x-s,s+t$. These sets coincide 
                only when $x=s+t$.
        \item   If $x<s$ then the new set of weights is $s-x,x,x,s+t$;
                this set cannot coincide with $s,t,s,x$ since $s+t$ is
                greater than every term in the second set.
        \end{enumerate}
\end{enumerate}
Thus, the hypothesis that there exists an edge of a weight $x\ne s,t,s+t$ 
leads to a contradiction. The proposition is proved.
\hfill$\Box$

\begin{remark} 
Notice that the operation of ripping and stitching introduced in the proof
of Proposition \ref{prop:chain-sts} often leads to another tree even in 
the case when the weights of all the edges of the tree under consideration 
are $s$, $t$, or $s+t$. Below we will often use this operation and call 
it $sts$-operation. 
\end{remark}

\begin{proposition}[Path $s$, $t$, $s+t$\,, I]
\label{prop:chain-st-sum}
Suppose that a unitree contains a path of three successive edges having
the weights $s$, $t$, $s+t$, and suppose also that $s\ne t$. Then
the edge of the weight $s+t$ is a leaf.
\end{proposition}

\begin{figure}[htbp]
\begin{center}
\epsfig{file=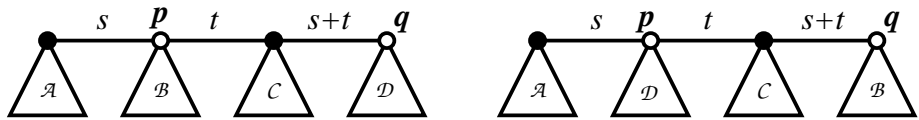,width=12cm}
\caption{\small Suppose that the subtree $\cal D$ is not empty.}
\label{fig:sts-I}
\end{center}
\end{figure}

\paragraph{Proof of Proposition \ref{prop:chain-st-sum}.} 
Take the tree shown in Figure~\ref{fig:sts-I}, left, and exchange the
subtrees $\cal B$ and $\cal D$. Obviously, both trees in the figure have 
the same passport. Suppose that the edge of the weight $s+t$ is not 
a leaf, so that the subtree $\cal D$ of the left tree is non-empty. 
According to  Lemma~\ref{lem:branches}, the edges  of ${\cal D}$ which 
are adjacent to the vertex $q$ have the same weight. Denote this weight
by $x$. By Lemma~\ref{lem:weight-exchange}, the possible values of $x$ 
are either~$t$ or $s+2t$. In the first case, we get a sub-path containing 
three edges of the weights $t$, $s+t$, $t$, and, according to 
Proposition~\ref{prop:chain-sts}, the only possible edge weights in such 
a tree could be $t$, $s+t$, and $s+2t$. But this contradicts the 
supposition that we have already an edge of the weight $s$ with $s\ne t$.

If, on the other hand, $x=s+2t$, then there are at least three 
{\em non-isomorphic}\/ subtrees attached to the vertex $p$ in the right
tree, since there are edges of three different weight $s$, $t$, and
$s+2t$ incident to this vertex. This situation violates 
Lemma~\ref{lem:branches}. 
\hfill$\Box$

\begin{remark}[$s=t$]
If $s=t$ then both arguments in the above proof are no longer valid,
though it is still difficult to make left and right trees of
Figure~\ref{fig:sts-I} isomorphic. (Recall that if $s=t$ then $s=t=1$,
see Corollary~\ref{cor:adj-edges}.) However, it is possible to construct
unitrees having the paths of length~4 with the weights $1,1,2,1$ 
(see the series~$G$, Figure~\ref{fig:diameter-4}), and the paths with 
the weights $1,1,2,3$ (see the tree~$P$, Figure~\ref{fig:sporadic-6}).
\end{remark}

\begin{proposition}[Path $s$, $t$, $s+t$\,, II]
\label{prop:chain-st-sum1} 
Suppose that a unitree contains a path of three 
successive edges having the weights $s$, $t$, $s+t$, and suppose also that 
the vertex adjacent to the edges of weights $s$ and $t$ has the valency $s+t$.
Then the edge of the weight $s+t$ is a leaf.
\end{proposition}

\paragraph{Proof.} 
Keeping notation of Figure~\ref{fig:sts-I}, it is enough to observe 
that if the subtree $\cal B$ is empty, then the two trees cannot be
isomorphic since the right one has more leafs than the left one. This
statement remains valid also for $s=t$.
\hfill$\Box$

\subsection{Brushes}
\label{sec:brushes}

A brush is a chain with two bunches of leaves attached to its ends: 
see formal definition below. Typical representatives of brushes are 
the trees shown in Figure~\ref{fig:brushes-E}. 

In this section we classify all brush unitrees.

\begin{definition}[Crossroad]\label{def:crossroad}
A vertex of a tree is {\em profound}\/ if, after having removed all the 
leaves from the tree, this vertex does not become a leaf. A vertex of a 
tree is a {\em crossroad}\/ if it is profound and has at least three 
branches going out of it.
\end{definition}

\begin{definition}[Brush]
A tree is called a {\em brush}\/ if it does not contain crossroads.
\end{definition}

\begin{proposition}[Brush unitrees]\label{prop:uni-brush}
A brush unitree belongs either to one of the series $A$, $B$, $C$, $D$, $E$ 
{\rm (Figures \ref{fig:stars-A}, \ref{fig:chains-B}, \ref{fig:brushes-C}, 
\ref{fig:brushes-D}, and \ref{fig:brushes-E})}, or to the series~$F$ with 
the degree of the central vertex $k=2$, or to the series $G$ with the 
degree of the central vertex $k=3$ {\rm (Figure~\ref{fig:diameter-4})}.
\end{proposition}

\paragraph{Proof.} In this part of the proof we only eliminate brush trees 
which are not unitrees. The uniqueness of the remaining brush trees will 
be proved later, in Section~\ref{sec:uniqueness}.

When all the edges of a tree are leaves we get the series~$A$ consisting 
of stars (Figure~\ref{fig:stars-A}). According Lemma~\ref{lem:branches}, 
only one of the leaves may have a weight different from the others. 

\begin{figure}[htbp]
\begin{center}
\epsfig{file=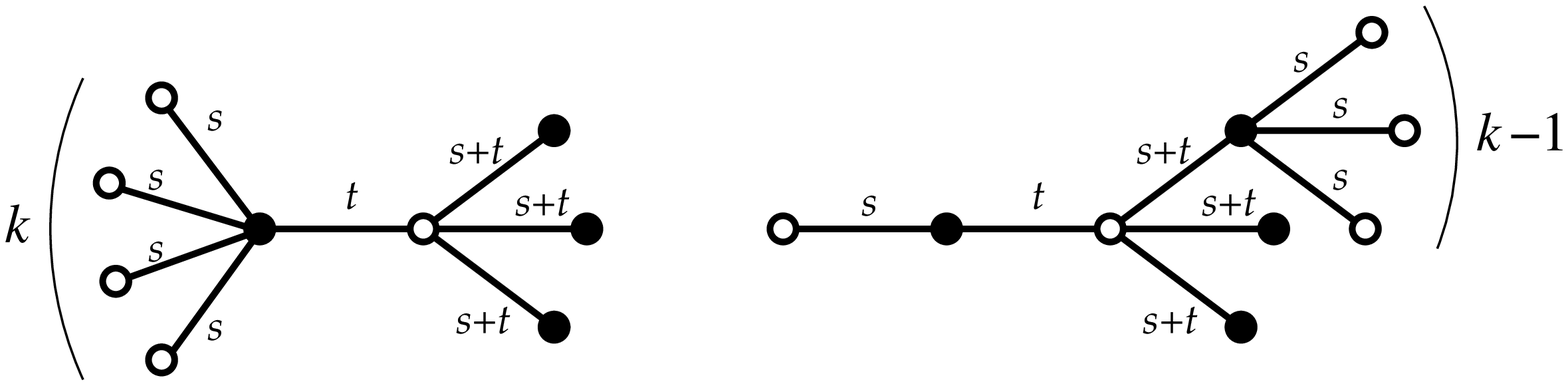,width=12cm}
\caption{\small If $k>1$, the tree on the left can be reconstructed 
into the one on the right. The new tree is different from the initial 
one since it has bigger diameter.}
\label{fig:not-C2}
\end{center}
\end{figure}

For the trees of diameter three, Lemma~\ref{lem:weight-exchange} leads 
to two possible patterns. One of them corresponds to the series $C$
(Figure~\ref{fig:brushes-C}); the other one is shown in
Figure~\ref{fig:not-C2}, left. We see that when $k>1$ we can transform
the left tree of Figure~\ref{fig:not-C2} into the right one. The new 
tree has the same passport but is not isomorphic to the initial one. 
Thus, the pattern shown in Figure~\ref{fig:not-C2}, left, is not a unitree. 

If $k=1$ then this pattern is a particular case of the series~$E_1$,
so it is a unitree.

\begin{figure}[htbp]
\begin{center}
\epsfig{file=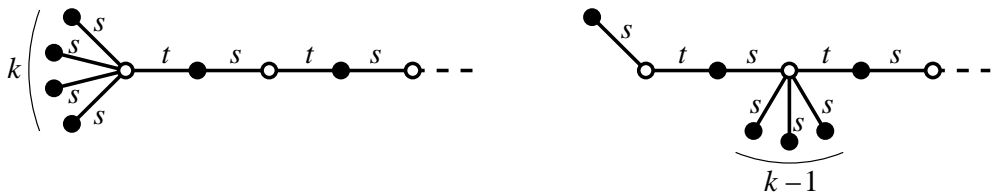,width=13cm}
\caption{\small These trees have the same passport. They are not isomorphic
since the right tree contains a crossroad while the left one does not. 
Therefore, if $k\ge 2$ then the weight of the leaves must be $s+t$ and 
not $s$ (or $t$).}
\label{fig:brush-pere}
\end{center}
\end{figure}

Now let us consider first the trees of diameter $\geq 5$, and after that 
return to the diameter~4 case. Suppose that a tree contains two adjacent
edges of weights $s$ and $t$ which are not leaves. It follows from 
Lemma \ref{lem:weight-exchange} and Proposition \ref{prop:chain-st-sum1} 
that the weights $s$ and $t$ alternate along all the path connecting 
vertices from which leafs grow. Furthermore, since this path contains 
at least three edges, it follows from Proposition~\ref{prop:chain-sts} 
that the only possible weight of a leaf which is not obtained by the 
further alternance of $s$ and $t$ is $s+t$. Now look at 
Figure~\ref{fig:brush-pere}, where an $sts$-operation is applied
to a brush tree having $k\ge 2$ leaves of weight $s$ on one of its ends.
The tree thus obtained, shown on the right, is distinct from the initial 
one since it contains a crossroad. A similar surgery can be made if there
are $l\ge 2$ leaves of the weight $s$ or $t$ (depending on the parity
of the diameter) on the right end of the tree. Thus, for the diameters 
$\ge 5$ only the types $B$ and $E$ survive. Namely, if a bunch of leaves 
at an end of the tree contains two or more leaves then the weight of these 
leaves is $s+t$.

\begin{figure}[htbp]
\begin{center}
\epsfig{file=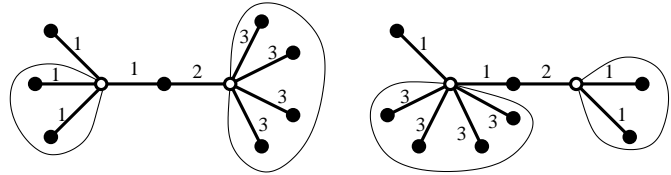,width=8.7cm}
\caption{\small The bunch of leaves of weight~1 which is transplanted
from left to right can be empty, if there is only one leaf on the left.}
\label{fig:no-3}
\end{center}
\end{figure}

The above argument fails for the brush trees of diameter~4: indeed, this 
time the operation shown in Figure~\ref{fig:brush-pere} does not create 
a crossroad. Therefore, the diameter~4 case needs a special consideration. 
Let us take a diameter of the tree, that is, a chain of length~4 going 
from one of its ends to the other. By Lemma~\ref{lem:weight-exchange} 
the sequence of the weights of its three first edges can be either 
$s$, $t$, $s$, or $s$, $t$, $s+t$, or $s+t$, $t$, $s$. Consider first the 
case $s$, $t$, $s+t$, so that the edge of the weight $s+t$ is not a leaf. 
In this case by Proposition~\ref{prop:chain-st-sum} we necessarily have
$s=t=1$ and a tree either belongs to the series~$G$, where the degree of 
the central vertex is $k=3$, or has the form shown in Figure~\ref{fig:no-3} 
on the left. In the last case, however, the tree under consideration is not 
a unitree, which can be seen by a transformation shown on the right. 

\begin{figure}[htbp]
\begin{center}
\epsfig{file=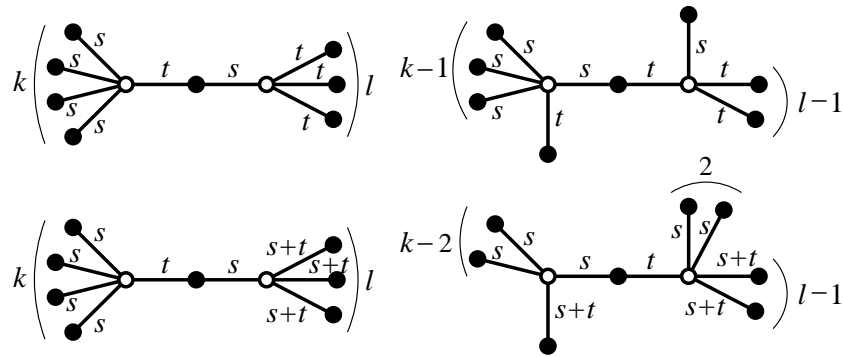,width=11cm}
\caption{\small The trees on the upper level have the same passport; 
they are different if $s\ne t$ and at least one of the parameters $k$, 
$l$ is greater than~1. The trees on the lower level also have the same 
passport; they are different if either $k>2$, or $l>1$, or both.}
\label{fig:diam4-pere}
\end{center}
\end{figure}

Assume now that the sequence of weights starts with $s$, $t$, $s$. Using 
Lemma \ref{lem:weight-exchange} again we see that it must be a part of 
one of the following three possible sequences: either $s,t,s,t$, or 
$s,t,s,s+t$, or $s,t,s,t-s$. For the latter one, taking $s'=t-s$, $t'=s$, 
we find the already considered above case $s',t',s'+t'$ read from right to 
left. Two other forms are shown in Figure~\ref{fig:diam4-pere} on the left. 

The first of these forms (above, left) can be transformed in a way shown 
on the right. The new tree is not equal to the initial one, unless either 
it belongs to the type $B$\/ (that is, $k=l=1$), or $s=t=1$, which is a 
particular case of the series~$F$, where the degree of the central vertex
is $k=2$. For the second form (below, left), the operation shown on the
right cannot be applied when $k=1$, that is, when the tree belongs to
the series $E_1$; and it does not change the tree when $k=2$ and $l=1$, 
which corresponds to the series $D$.

Finally, if the sequence of weights of edges of a diameter starts as $s+t$, 
$s$, $t$, then either a tree belongs to the series~$E_4$ (or $E_2$), or 
the sequence of weights of the edges of a diameter is $s+t$, $s$, $t$, 
$s-t$. In the latter case, however, Proposition~\ref{prop:chain-st-sum} 
yields that $s-t=t=1$, implying that a tree is the one shown in 
Figure~\ref{fig:no-3} on the left.
\hfill$\Box$




\subsection{Trees with repeating branches of height 2} 

\begin{figure}[htbp]
\begin{center}
\epsfig{file=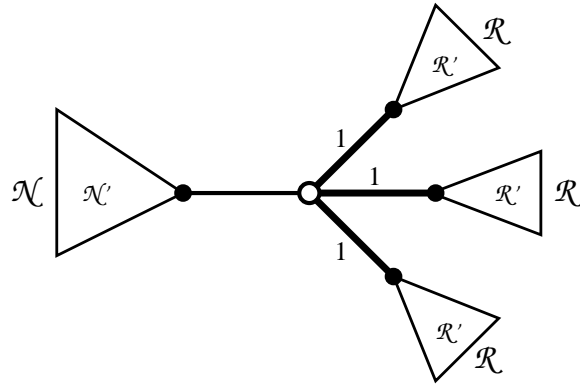,width=7.5cm}
\caption{\small A typical tree with a crossroad. The subtrees ${\cal R}'$
are all non-empty. The branch $\cal N$ is also non-empty. It may, or may 
not be isomorphic to $\cal R$.}
\label{fig:crossroad}
\end{center}
\end{figure}

From now on we will assume that the trees we consider are not brushes, 
that is, they contain at least one crossroad. Recall that a crossroad 
is a profound vertex at which three or more branches meet, see 
Definition~\ref{def:crossroad}. A~typical tree with a crossroad is shown 
in Figure~\ref{fig:crossroad}. According to Lemma~\ref{lem:branches}, 
all the branches attached to the crossroad, except maybe one, are 
isomorphic as rooted trees, with their root edges (shown by thick lines 
in the figure) being incident to the crossroad. We call these branches 
{\em repeating}\/; in the figure they are denoted by the same 
letter~$\cal R$; the subtrees of $\cal R$ attached to the root are denoted 
by ${\cal R'}$. The subtrees ${\cal R'}$, are non-empty since otherwise 
the vertex to which $\cal R$ and $\cal N$ are attached would not be 
profound. The number of branches of the type $\cal R$ can be two or more, 
but the majority of the transformations given below involve only two 
branches; therefore, in the majority of pictures we will draw only two 
repeating branches. 

By convention, we suppose that {\em the branch $\cal N$ is always non-empty}.
If all the branches meeting at the crossroad are isomorphic to~$\cal R$
and might therefore be all considered as repeating, we take an arbitrary 
one of them and, somewhat artificially, declare it to be the ``non-repeating''
branch~$\cal N$. The subtree ${\cal N'}$ has a right to be empty.

The roots of repeating branches are adjacent edges which are not leaves. 
Therefore, according to Corollary~\ref{cor:adj-edges}, their weights must 
be equal to~1. Finally, the {\em height}\/ of a repeating or non-repeating
branch is the distance from its root vertex (that is, the crossroad) to 
its farthest leaf. 

In the previous section we classified the brush unitrees, which are 
by definition unitrees without crossroads. In this section we establish 
a complete list of all possible unitrees whose repeating branches all have 
the height~2. More precisely, we assume that for {\em any}\/ crossroad 
of a unitree under consideration the repeating branches are of height~2.

\begin{proposition}[Repeating branches of height 2]\label{prop:repdiam2}
A unitree whose all repeating branches are of height\/~$2$ belongs 
to one of the types $F$, $G$, $H$, or $K$.
\end{proposition}

\begin{figure}[htbp]
\begin{center}
\epsfig{file=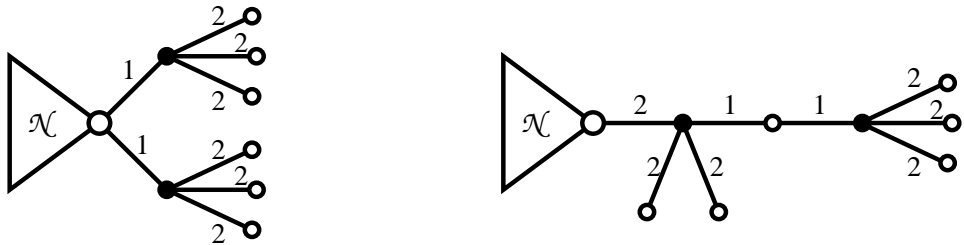,width=12.6cm}
\caption{\small A transformation of repeating branches of height~2 with 
leaves of weight 2.}
\label{fig:proof-rep-brush-2}
\end{center}
\end{figure}

\paragraph{Proof.} First of all observe that weights of leaves of repeating 
branches cannot be equal to 2 since otherwise the transformation shown in 
Figure~\ref{fig:proof-rep-brush-2} leads to a non-isomorphic tree (the 
tree on the right has fewer leaves than the one on the left). Thus, the 
weights of these leaves are equal to~1. Therefore, according to 
Lemma~\ref{lem:weight-exchange} the root edge of $\cal N$ can only be 
of weight~1 or~2. The case ${\cal N'}=\pusto$ (that is, when $\cal N$ 
consists of a single edge, see Figure~\ref{fig:crossroad}) is a particular 
case of the series $F$ and $G$. Suppose then that ${\cal N'}\ne\pusto$.

\begin{figure}[htbp]
\begin{center}
\epsfig{file=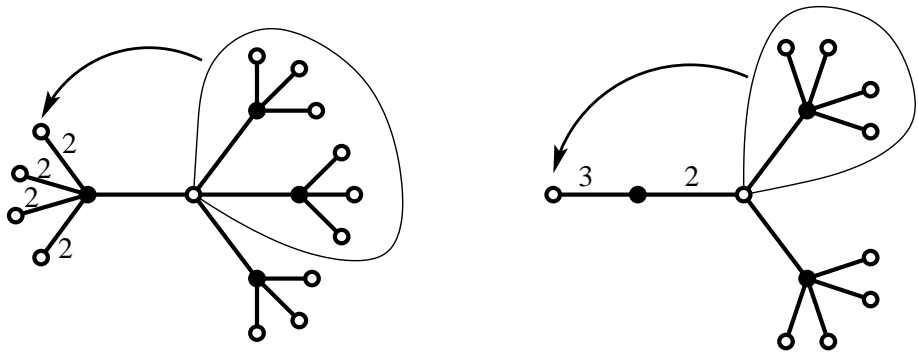,width=12cm}
\caption{\small Non-repeating branch of height 2. These transformations
show that the weight of its leaves cannot be 2 or 3.}
\label{fig:nonrep-12}
\end{center}
\end{figure}

If the non-repeating branch is of the height~$2$, that is, if it is
a root edge with a bunch of leaves attached to it, then these leaves
could be of weight~1 or~2 when the root edge is of weight~1, and they
could be of weight~1 or~3 when the root edge is of weight~2. However, 
the leaves of the weights 2 and 3 are impossible, as two transformations
of Figure~\ref{fig:nonrep-12} show. In this figure, we take all the 
repeating branches but one and re-attach them to one of the leaves of the 
non-repeating branch. These transformations always change the trees:
on the left, there appears a non-leaf of weight~2, and on the right, there 
appears a non-leaf of weight~3. Thus, all the leaves of the non-repeating 
branch must be of weight~1.

In addition, when the root edge of the non-repeating branch is of weight~2 
this branch cannot be isomorphic to the repeating branches. This observation 
implies that the degree of the black vertex lying on this branch must be 
equal to the degrees of the black vertices of repeating branches since 
otherwise an exchange of leaves between repeating and non-repeating branches 
could be possible. Thus, the only remaining possibilities are the trees of 
the types $F$~and~$G$, see Figure~\ref{fig:diameter-4}.

\medskip

\begin{figure}[htbp]
\begin{center}
\epsfig{file=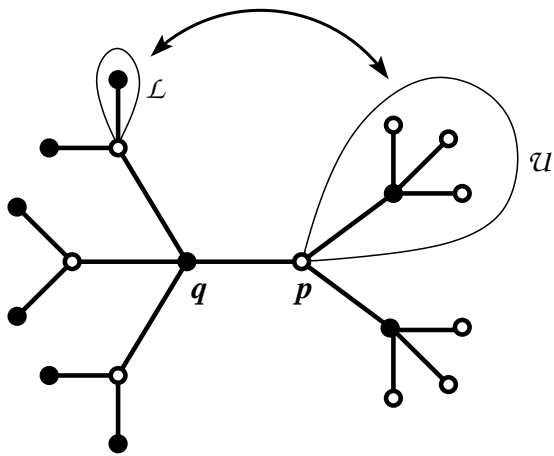,width=7.2cm}
\caption{\small A leaf $\cal L$ can be interchanged with a repeating
branch $\cal U$.}
\label{fig:rep-to-rep}
\end{center}
\end{figure}

Consider now the case of a non-repeating branch of height $\ge 3$.
First suppose that the vertex~$q$, which is the nearest neighbor 
of the crossroad vertex~$p$ when we move
along the non-repeating branch, is itself a crossroad. According
to our supposition, the tree does not have repeating branches of the
height greater than~2. Hence, the repeating branches growing out of~$q$
are of height~2. But then a leaf $\cal L$ of such a branch can be
interchanged with a repeating branch $\cal U$ growing out of~$p$,
see Figure~\ref{fig:rep-to-rep}. This operation would create at least
three different trees attached to~$p$: one of them would be of height~1,
another one of height~2, and a third one of height~4. Thus, this
possibility is ruled out.

\begin{figure}[h!]
\begin{center}
\epsfig{file=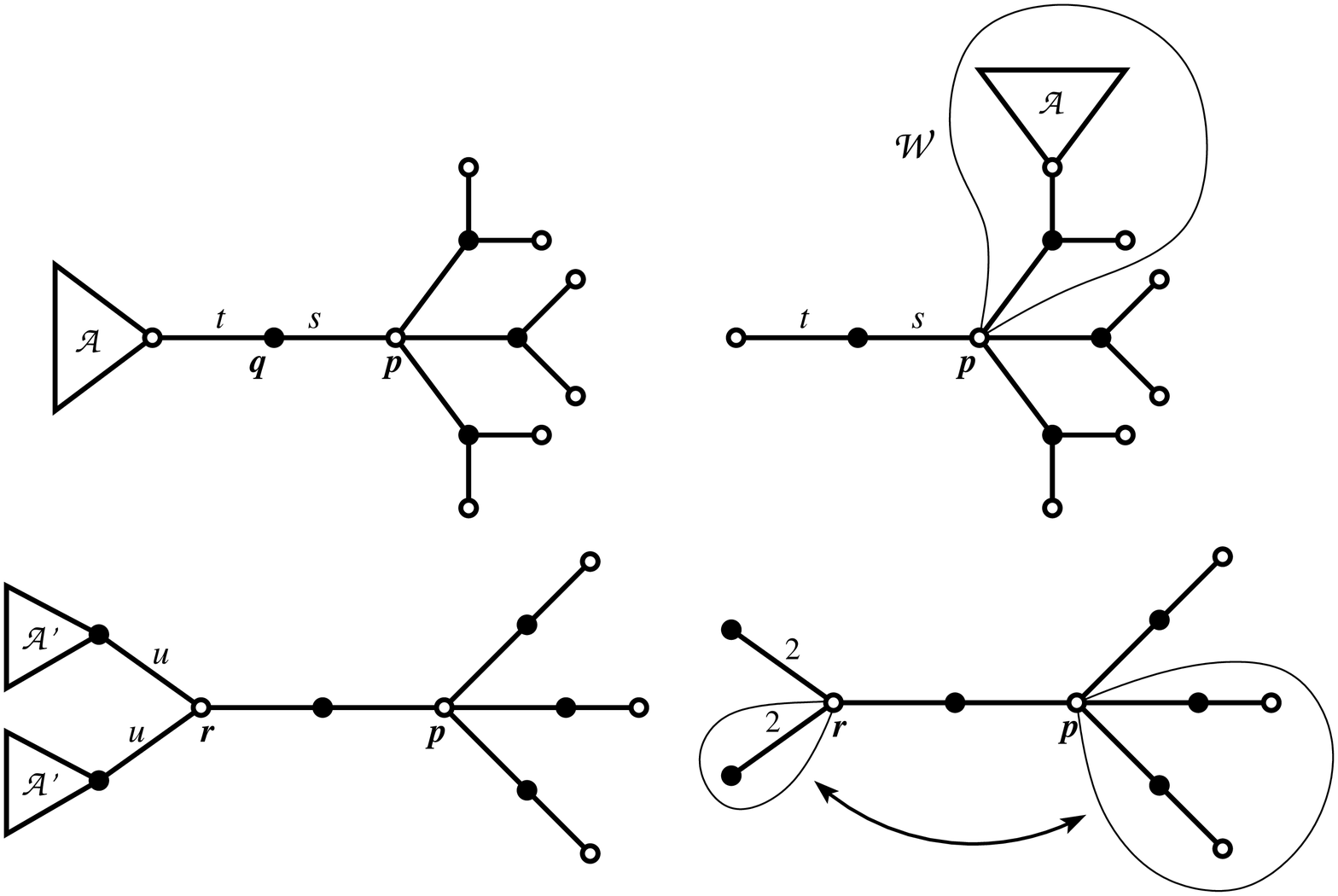,width=14.4cm}
\caption{\small Illustration to the proof of Proposition~\ref{prop:repdiam2}.}
\label{fig:rep-2-nonrep-3}
\end{center}
\end{figure}

Suppose next that the vertex $q$ is not a crossroad. Then the tree looks
like the one in Figure~\ref{fig:rep-2-nonrep-3}, top left, where $\cal A$
is non-empty. A priori, there are four possibilities for the values
$(s,t)$, namely, $(1,1)$, $(1,2)$, $(2,1)$, and $(2,3)$. The case $(2,3)$
can be immediately ruled out since the edge of weight~3 should be
a leaf by Proposition~\ref{prop:chain-st-sum}, but we have supposed 
that ${\cal A}\ne\pusto$.

The cases $(s,t)=(1,1)$ or $(2,1)$ can be treated together. When $t=1$ we 
can re-attach $\cal A$ to one of the leaves of the repeating branches, as 
is shown in the same figure on the top right. Among the branches attached to 
the vertex $p$ of the tree thus obtained there is only one branch of a 
height greater than~2: it is $\cal W$. Therefore, all the remaining branches 
are repeating, so we may conclude that $s=1$ (the case $s=2$ is impossible), 
and all the repeating branches have only one leaf. Thus, the tree looks like 
the one on the bottom left in Figure~\ref{fig:rep-2-nonrep-3}, where two
possibilities may occur: either $u=1$; or $u=2$, and then, according to 
Proposition~\ref{prop:chain-st-sum1}, ${\cal A'}=\pusto$, so the edges 
of the weight $u=2$ are leaves.

In the first case, we can exchange the repeating branches attached
to the vertices $p$~and~$r$. Therefore, they all must be equal, and we
get a tree of the type~$H$, see Figure~\ref{fig:diameter-6}.
In the second case, we can interchange one of the leaves of weight~2 with 
two repeating branches attached to~$p$. The only tree which does not change 
after this transformation is the one which has exactly one leaf of weight~2 
and exactly two repeating branches attached to~$p$, that is, the tree~$K$, 
see Figure~\ref{fig:sporadic-5}.

\medskip

\begin{figure}[htbp]
\begin{center}
\epsfig{file=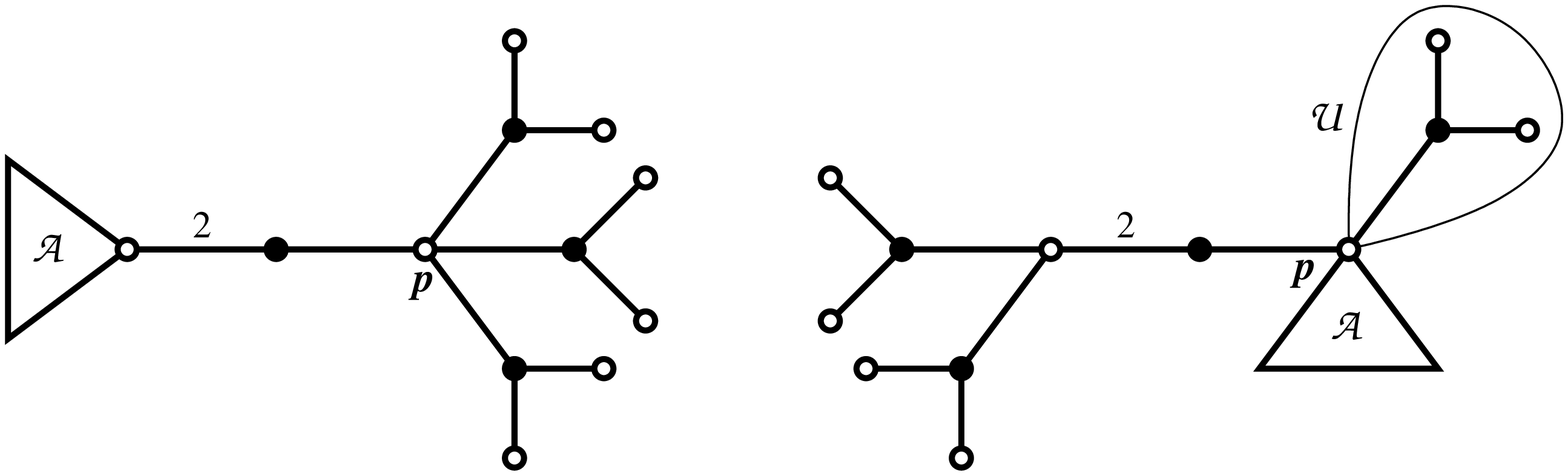,width=15cm}
\caption{\small Illustration to the proof of Proposition~\ref{prop:repdiam2}.}
\label{fig:outrule}
\end{center}
\end{figure}

There remains the last case to be ruled out: when the tree shown in
Figure~\ref{fig:rep-2-nonrep-3}, top left, has $s=1$ and $t=2$; see
also Figure~\ref{fig:outrule}, left. In this case we can interchange
the subtree~$\cal A$ with all but one repeating branches, see the tree
on the right of Figure~\ref{fig:outrule}. We see that $\cal A$ must
consist of several copies of the branch $\cal U$ since otherwise
$\cal A$ should consist of copies of the longer branch at~$p$, and
we would get repeating branches of the height greater than~2. 
Then, we may take the left tree of  Figure~\ref{fig:outrule}, cut all 
the repeating branches form $p$ and re-attach them to one of the leaves 
of $\cal A$. This operation will necessarily produce a different tree 
since the only edge of the weight~2 will be now at distance~1 from the 
leaf while it was at distance~2 in the initial tree. Therefore, this 
case is impossible.

\medskip
 
Proposition~\ref{prop:repdiam2} is proved.
\hfill $\Box$

\subsection{Trees with repeating branches of the type $(1,s,s+1)$}

\begin{figure}[!h]
\begin{center}
\epsfig{file=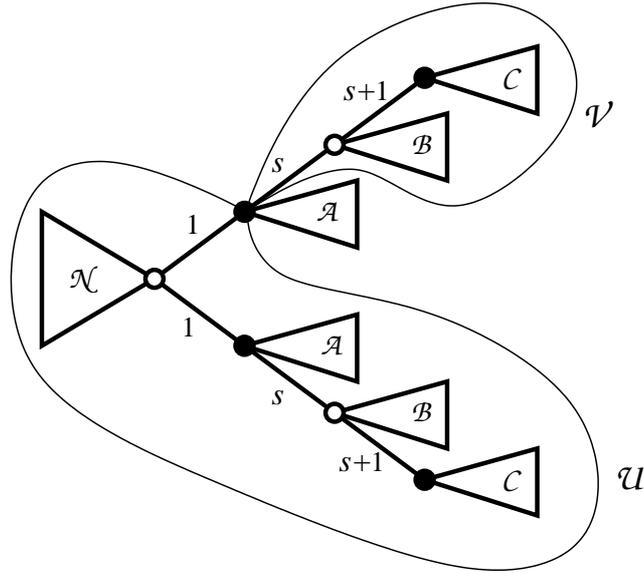,width=8.4cm}
\caption{\small Illustration to the proof of Lemma~\ref{lem:xxx}: the
subtrees $\cal A$ and $\cal C$ are empty; the subtree $\cal B$ is a
bunch of leaves of weight $s+1$.}
\label{fig:1.t.t+1}
\end{center}
\end{figure}

If a unitree has a crossroad from which grow repeating branches of height
$>2$, then these branches ``start'' either with a path $1$, $s$, $s+1$,
or with a path $1$, $t$, $1$, where $s$ and $t$ here may be equal to
either~1 or~2 (see Figures \ref{fig:1.t.t+1} and \ref{fig:1.t.1}).   
In this subsection we classify unitrees which have no crossroads
of the second type. We start with the following lemma.

\begin{figure}[htbp]
\begin{center}
\epsfig{file=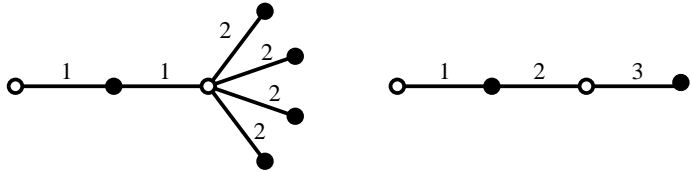,width=9cm}
\caption{\small Illustration to the proof of Lemma~\ref{lem:xxx}.}
\label{fig:xxx}
\end{center}
\end{figure}

\begin{lemma}\label{lem:xxx} 
If a unitree has a repeating branch of the type $(1,s,s+1)$, then this 
branch has one of the two forms shown in\/ {\rm Figure~\ref{fig:xxx}.} 
Furthermore, in the second case the unitree is necessarily the tree $P$.
\end{lemma} 

\vspace{-5mm}

\paragraph{Proof.} First of all observe that the subtrees $\cal A$ and 
$\cal C$ in Figure~\ref{fig:1.t.t+1} can be interchanged, and if one of
them was empty while the other was not, this operation would change 
the number of leaves, so that the tree in question could not be a unitree.
We will show now that the assumption that both trees $\cal A$ and $\cal C$ 
are not empty also leads to a contradiction (so that, in fact, both of
them are empty).

Since the tree $\cal V$ is isomorphic to a subtree of $\cal U$ and is
therefore distinct from $\cal U$, if $\cal A$ is not empty, then it consists 
of a certain number of copies of $\cal U$ or of $\cal V$. The first case 
is impossible since $\cal A$ is a subtree of $\cal U$. Therefore, $\cal A$ 
consists of a certain number of copies of $\cal V$ implying that $\cal C$ 
is a proper subtree of $\cal A.$ Now, interchanging $\cal A$ and $\cal C$
in every repeating branch we may prove in the same way that that  $\cal A$ 
is a proper subtree of $\cal C$, implying the contradiction that we need. 
Thus, $\cal A$~and~$\cal C$ are empty. In particular, $\cal B$ is merely 
a collection of leaves of weight $s+1$.

\begin{figure}[!h]
\begin{center}
\epsfig{file=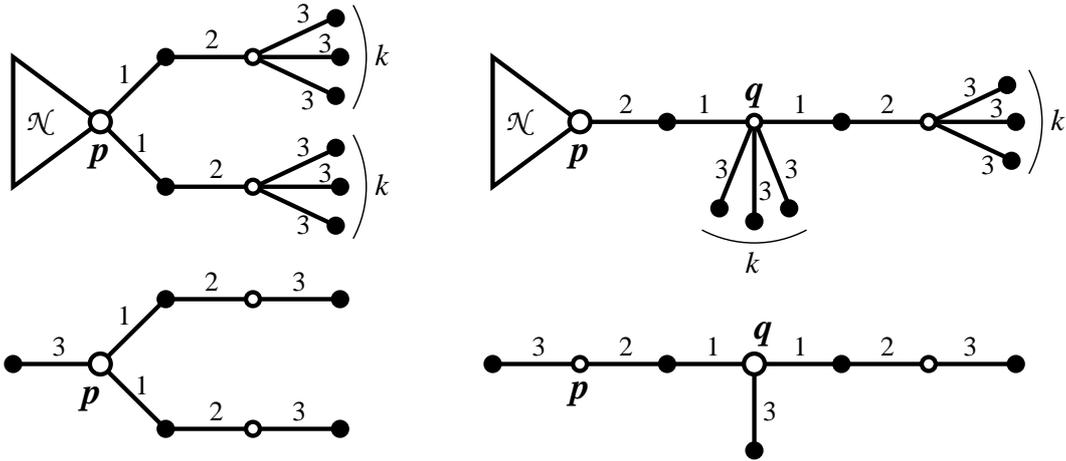,width=14cm}
\caption{\small Illustration to the proof of Lemma~\ref{lem:xxx}: 
transformations of repeating branches of height~3 with the weight 
sequence $1,2,3$.}
\label{fig:proof-rep-brush-12}
\end{center}
\end{figure}

Assume now that $s=2$. Then our tree must look as in 
Figure~\ref{fig:proof-rep-brush-12}, top left,
where the number of repeating branches at the vertex $p$ might be two 
or more. Let us take two of these branches, apply the transformation
shown on top right, and see what takes place at the vertex $q$.
According to Lemma~\ref{lem:branches},
all the subtrees growing out of this vertex, except maybe one, must
be isomorphic. This can only happen when $k=1$, and the subtree growing
from $q$ to the left is isomorphic to the one growing from $q$ to the
right. Therefore, before the transformation there were exactly two 
(and not more) repeating branches at $p$, and the subtree $\cal N$ was
reduced to a single leaf of weight~3. The resulting situation is shown 
in Figure~\ref{fig:proof-rep-brush-12}, bottom left. In this case, our
transformation can still be applied, but it leads to a tree isomorphic
to the initial one. The unitree thus obtained is~$P$, see 
Figure~\ref{fig:sporadic-6}.
\hfill $\Box$

\begin{proposition}[Branches of the type $(1,s,s+1)$]\label{prop:repsts+1}
A unitree which has at least one crossroad of type $(1,s,s+1)$ 
but no crossroads of type $(1,t,1)$ belongs to one of the types 
$J$, $L$, $N$, $M$, $O$, $P$, $R$, or $S$.
\end{proposition}

\begin{figure}[htbp]
\begin{center}
\epsfig{file=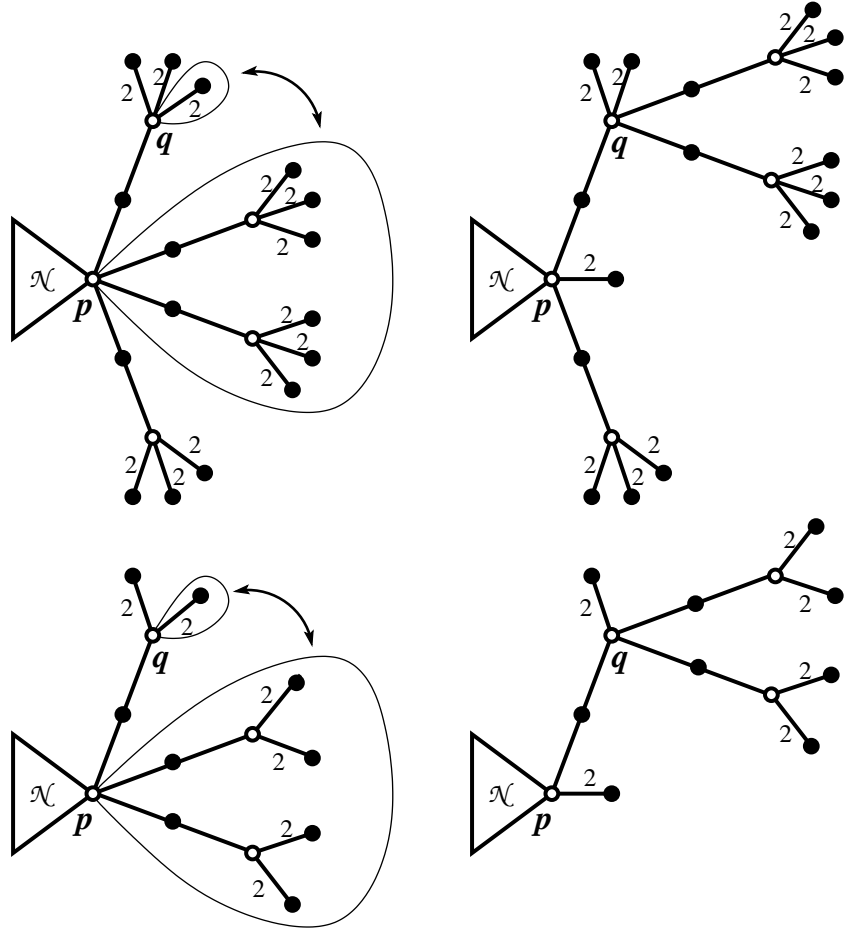,width=11cm}
\caption{\small Illustration to the proof of 
Proposition~\ref{prop:repsts+1}: transformations of repeating branches
of height~3 with the weight sequence $1,1,2$.}
\label{fig:rep-branch-112}
\end{center}
\end{figure}

\paragraph{Proof.}
In view of Lemma \ref{lem:xxx} we may assume that the repeating branches 
have the form shown in Figure~\ref{fig:xxx} on the left. Suppose first
that the number of the repeating branches is three or more, and apply the 
transformation shown in Figure~\ref{fig:rep-branch-112}, that is, 
interchange the positions of a leaf of weight~2 and of a pair of repeating 
branches. If the number of the repeating branches was more than three 
then the principle ``all branches except maybe one are isomorphic'' would 
be violated at the vertex~$p$. The same principle would be violated at the
vertex $q$ if the number of leaves in a repeating branch was more than two. 
Therefore, the number of repeating branches is three, and our transformation 
looks as is shown in Figure~\ref{fig:rep-branch-112}, bottom.
If the number of leaves in a repeating branch is two, then  applying 
once again the same principle at the point $q$, we arrive at the tree~$O$. 
Assume now that this number is equal to one. Then the height of~$\cal N$ 
is less than two, since otherwise the new tree would have more crossroads 
than the initial one. Furthermore, if~$\cal N$ is a bunch of leaves of weight 
two, then  we could transfer all these leaves to the vertex $q$ changing 
the tree. Therefore,~$\cal N$ is empty, and we arrive at the tree $N.$

\begin{figure}[htbp]
\begin{center}
\epsfig{file=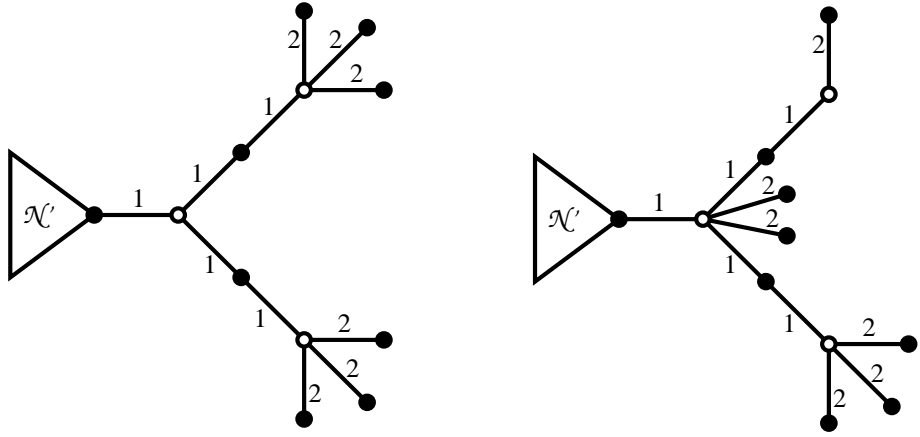,width=12cm}
\caption{\small Illustration to the proof of 
Proposition~\ref{prop:repsts+1}: two repeating branches.}
\label{fig:rep-br-R1}
\end{center}
\end{figure}

Suppose next that the number of repeating branches is two. The starting
edge of the non-repeating branch $\cal N$ is either of weight~2, and then,
according to Proposition~\ref{prop:chain-st-sum1},
it is a leaf, and we get the tree~$J$ (Figure~\ref{fig:sporadic-6}), or
it is of weight~1. In the latter case it does not have to be a leaf, though
this situation imposes another constraint: the repeating branches must
have only one leaf, otherwise the transformation shown in
Figure~\ref{fig:rep-br-R1} can be applied, producing three non-isomorphic 
branches growing from the crossroad.

What remains is to study more attentively the structure of the
non-repeating branch $\cal N$. Here we consider the following
cases:
\begin{itemize}
\item   The height of $\cal N$ is 1, that is, ${\cal N'}=\pusto$.
\item   The height of $\cal N$ is 2.
\item   The height of $\cal N$ is 3 or more, and $\cal N$ starts with a
        path having the weights $1$, $s$, $1$ where $s$ is equal to
        1~or~2.
\item   The height of $\cal N$ is 3 or more, and $\cal N$ starts with a
        path having the weights $1$, $s$, $s+1$ where $s$ is equal to
        1~or~2.
\end{itemize}

The height of $\cal N$ equal to 1 case is trivial: we get the tree $L$
(Figure~\ref{fig:sporadic-6}).

\ssk

\begin{figure}[htbp]
\begin{center}
\epsfig{file=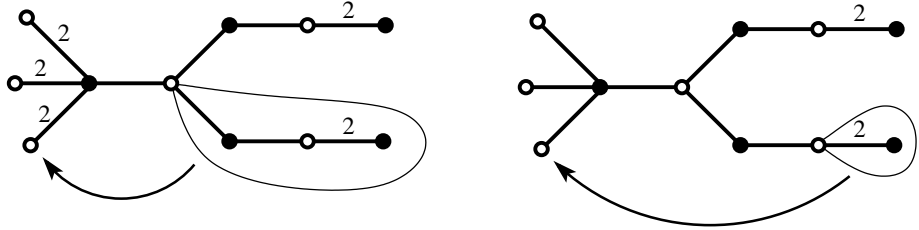,width=12cm}
\caption{\small Illustration to the proof of 
Proposition~\ref{prop:repsts+1}: non-repeating branch of height 2.}
\label{fig:nonrep-height-2}
\end{center}
\end{figure}

The height of $\cal N$ equal to 2 case is illustrated in
Figure~\ref{fig:nonrep-height-2}. If the weight of the leaves of the
non-repeating branch is equal to~2 then, whatever is their number, the
reattachment shown on the left changes the tree since the new
tree has one leaf less than the initial one. If the weight of the
leaves of the non-repeating branch is equal to~1 then the reattachment
shown on the right also changes the tree unless there is only one leaf
in the non-repeating branch. The latter case gives us the tree~$M$
(Figure~\ref{fig:sporadic-6}).

\begin{figure}[!h]
\begin{center}
\epsfig{file=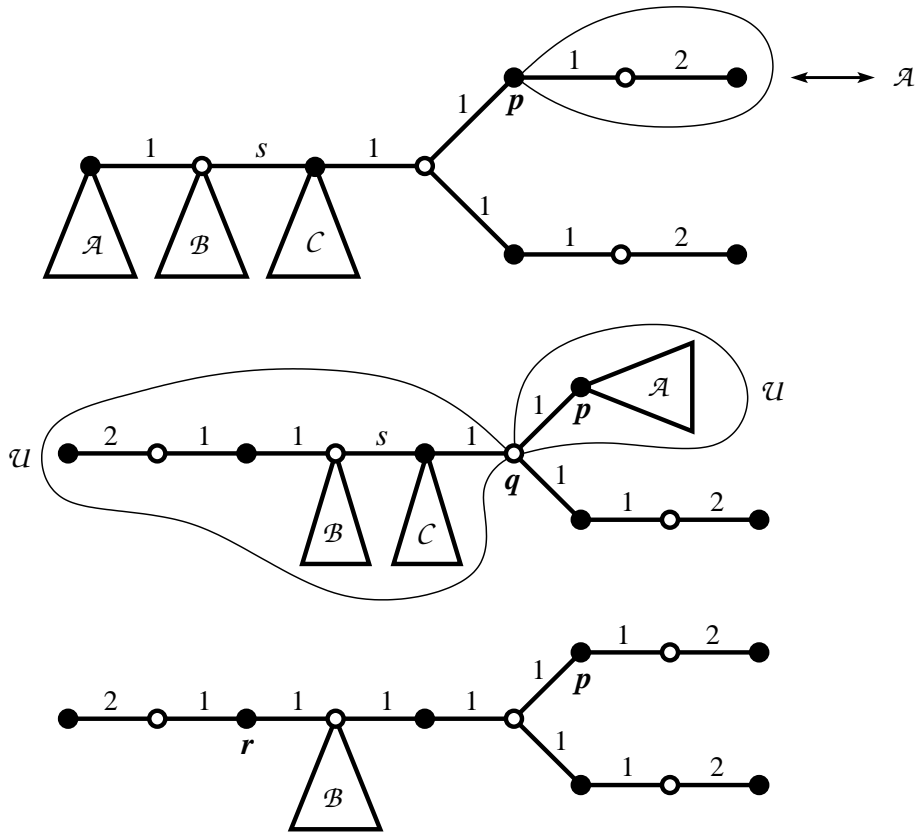,width=12cm}
\caption{\small Illustration to the proof of 
Proposition~\ref{prop:repsts+1}: a non-repeating branch containing a path
with the weights 1, $s$, 1.}
\label{fig:nonrep.1.s.1}
\end{center}
\end{figure}

Assume now that the height of the non-repeating branch is $\ge 3$, and
this branch contains a path with the weights 1,~$s$,~1, see the upper
tree in Figure~\ref{fig:nonrep.1.s.1}. First of all, we remark that the
subtree $\cal A$ may be interchanged with a chain of length~2 attached
to the vertex~$p$. Then, according to the principle ``all branches except 
maybe one are isomorphic'', two situations may occur. First, we could
thus create two repeating branches $\cal U$ attached to the vertex $q$, 
see the tree in the middle. But such a tree would contain repeating 
branches of the type $(1,s,1)$ which contradicts our supposition. 
The other possibility is that $\cal A$ is equal to the chain which was 
attached to $p$. Then we get a vertex~$r$ (see the lower tree) which 
is of degree~2 and is incident to two edges of weights 1~and~1. 
Therefore, according to Proposition~\ref{prop:chain-st-sum1}, the edge 
of weight~$s$, which is not a leaf, cannot have weight~2; hence, $s=1$. 
Finally, we affirm that ${\cal C}=\pusto$, otherwise it could be 
reattached to the vertex~$p$ and we would get three different trees 
attached to $q$. The resulting tree is shown in Figure~\ref{fig:nonrep.1.s.1}, 
bottom. If ${\cal B}=\pusto$ we get the tree~$S$ (Figure~\ref{fig:sporadic-8}).
If ${\cal B}\ne\pusto$ then, again according to the principle ``all branches
except maybe one are isomorphic'', $\cal B$ must be equal to a chain
of weights 1,~1,~2, and we get the tree~$R$ (Figure~\ref{fig:sporadic-8}),
since otherwise a tree would contain repeating branches of the type $(1,s,1)$. 
(Note that in the last case we obtain the tree $T$ which is considered in 
Proposition~\ref{prop:repsts+2} which treats the
case of repeating branches of the type $(1,t,1)$.)

\begin{figure}[!h]
\begin{center}
\epsfig{file=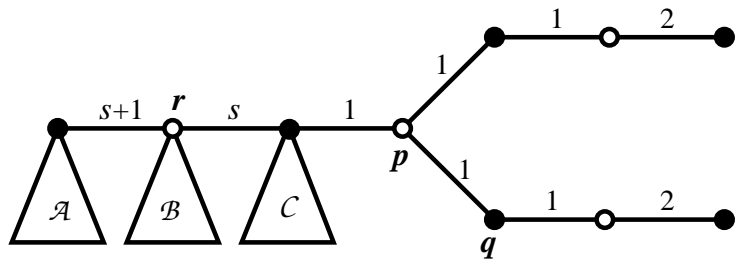,width=9.6cm}
\caption{\small Illustration to the proof of 
Proposition~\ref{prop:repsts+1}: a non-repeating branch containing a path
with the weights 1, $s$, $s+1$.}
\label{fig:nonrep.1.s.s+1}
\end{center}
\end{figure}

Finally, consider the case when the non-repeating branch is of height
$\ge 3$ and contains a path with the weights 1,~$s$,~$s+1$, see 
Figure~\ref{fig:nonrep.1.s.s+1}. We affirm that in this case $s=1$
and all the three subtrees $\cal A$, $\cal B$, $\cal C$ are empty,
so that we get the tree $N$ of Figure~\ref{fig:sporadic-6} (and what
we call ``non-repeating branch'' is in this case equal to the repeating
branches). Indeed, the tree contains a path with the weights~1,~1,~1; 
therefore, according to Proposition \ref{prop:chain-sts}, the only 
possible weights are 1~and~2, so $s=1$. Now, if  ${\cal A}\ne\pusto$
then it could be reattached to the vertex $q$ thus producing a tree
with one more leaf. Therefore, ${\cal A}=\pusto$. Then, $\cal C$ is
also empty since $\cal A$ and $\cal C$ can be interchanged. Finally,
if ${\cal B}\ne\pusto$ then there are two possibilities. Either $\cal B$
is a bunch of leaves of weight~2; but then it can be reattached to 
the vertex~$p$. Or $\cal B$ is a number of copies of the long branch
growing out of the vertex $r$; but then, once again, we would create
repeating branches of the type $(1,t,1)$.

\ssk

Proposition~\ref{prop:repsts+1} is proved.
\hfill $\Box$

\newpage

\subsection{Trees with repeating branches of the type $(1,t,1)$}

In this subsection we classify unitrees which have crossroads
of the type $(1,t,1)$.

\begin{proposition}[Branches of type $(1,t,1)$]\label{prop:repsts+2}
A unitree which has at least one crossroad of type $(1,t,1)$
belongs to one of the types $I$, $Q$, or $T$.
\end{proposition}

\begin{figure}[htbp]
\begin{center}
\epsfig{file=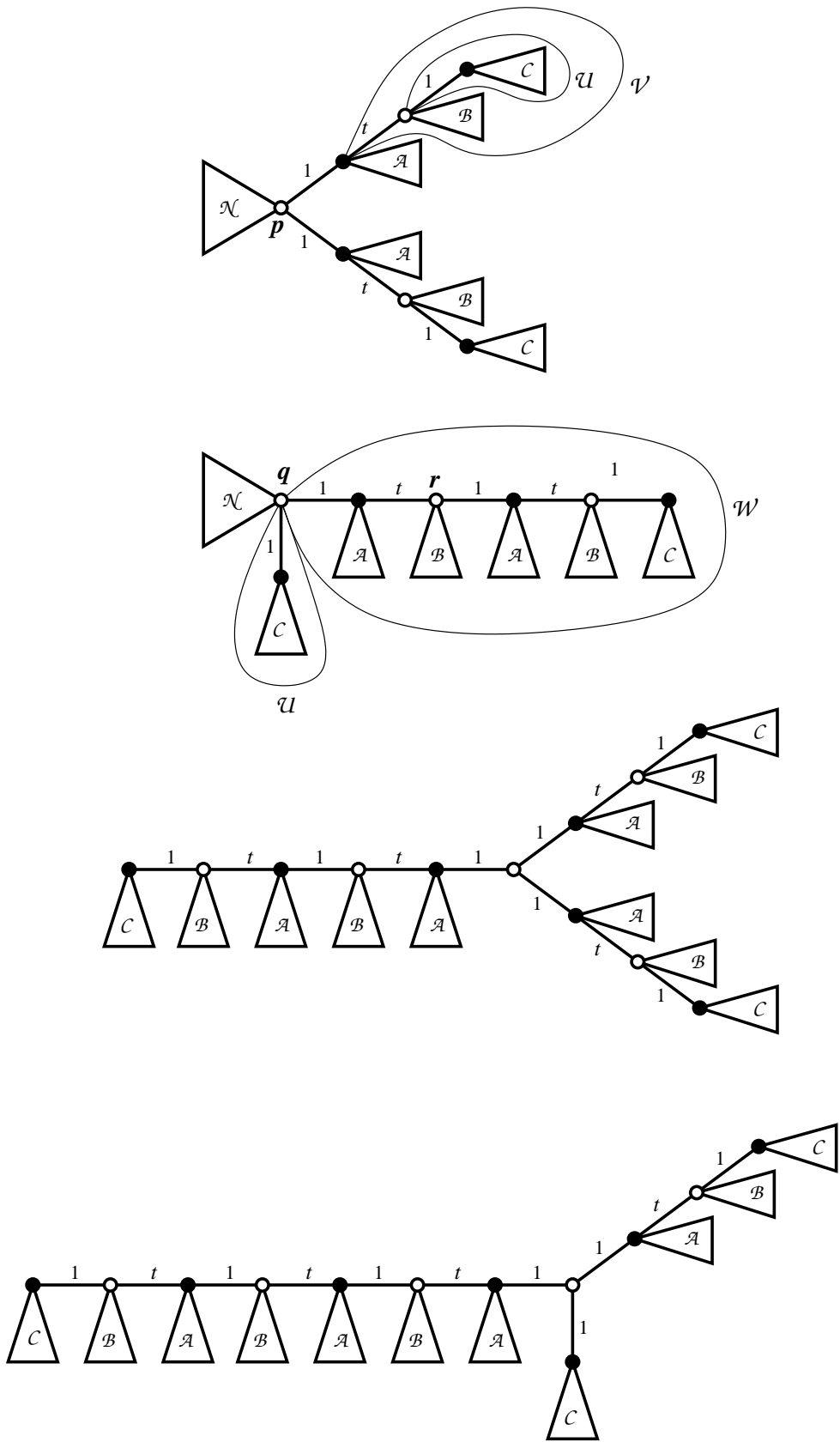,width=12.8cm}
\caption{\small Illustration to the proof of Proposition~\ref{prop:repsts+2}.}
\label{fig:1.t.1}
\end{center}
\end{figure}

\paragraph{Proof.} First of all, observe that by Lemma~\ref{lem:branches} 
the subtree $\cal B$ is a collection of 
copies of the subtree $\cal U$, and the subtree $\cal A$ is a collection 
of copies of the subtree $\cal V$ (see Figure~\ref{fig:1.t.1}). 
Further, an $sts$-operation, applied to the first tree of 
Figure~\ref{fig:1.t.1} gives the second tree shown in this figure.
This image implies that there are only two repeating branches growing 
from the vertex~$p$, otherwise the tree would certainly change. 
Now, looking at the vertex $q$ of the second tree of Figure~\ref{fig:1.t.1} 
we see that either ${\cal N}={\cal U}$ or ${\cal N}={\cal W}$. If 
${\cal N}={\cal W}$ then the initial tree would look like the third tree
of the same figure. Then we could once again apply an $sts$-transformation 
and make the long branch even longer, and one of the repeating branches,
shorter (see the fourth tree of the figure), which would give us three 
different branches attached to~$p$. Hence, ${\cal N}={\cal U}$. 
In particular, we have proved that whenever a unitree has a crossroads 
of type $(1,t,1)$ the corresponding non-repeating branch is a subtree 
of the repeating branch.

\begin{figure}[htbp]
\begin{center}
\epsfig{file=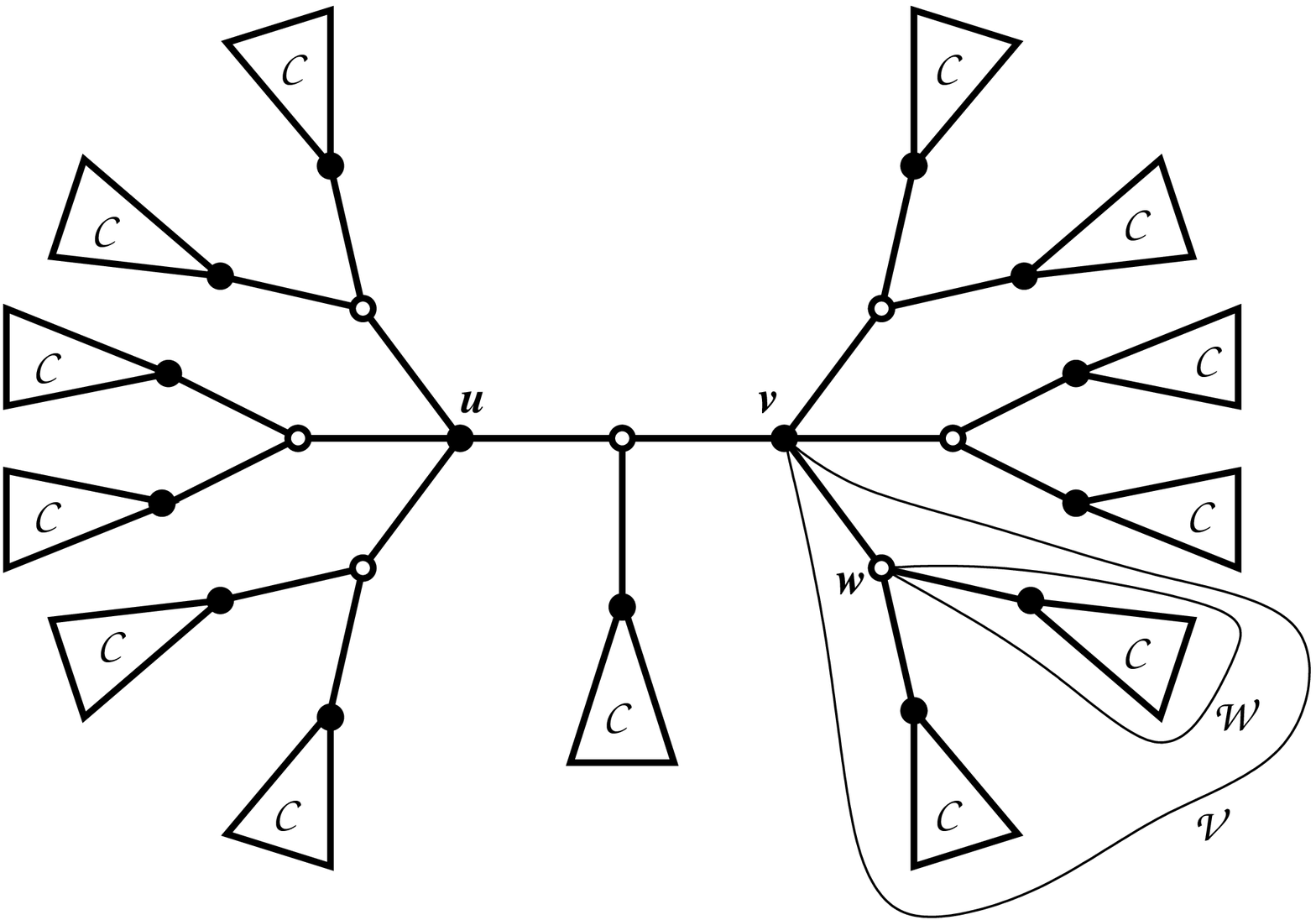,width=10cm}
\caption{\small Illustration to the proof of Proposition~\ref{prop:repsts+2}.}
\label{fig:1.t.1-bis}
\end{center}
\end{figure}

Now, it follows from ${\cal N}={\cal U}$ that the first tree of 
Figure~\ref{fig:1.t.1} has a (unique) center at the vertex $p$ while the 
second one has a (unique) center at the vertex $r$. Hence, the vertex~$p$ 
of the first tree must correspond to the vertex $r$ of the second one, 
and thus we must have $t=1$ and ${\cal B}={\cal N}={\cal U}$.
Therefore, the tree has the form shown in Figure~\ref{fig:1.t.1-bis},
with the same number $l\geq 1$ of branches growing out of the vertices
$u$ and $v$. 

If ${\cal C}=\pusto$ then we get a tree of the type $I$. Thus we may 
assume that ${\cal C}$ is non-empty. Observe first that $l=1$. Indeed, 
if $l>1$ then $\cal V$ is a repeating branch. Furthermore, since 
${\cal C}$ is non-empty, $\cal V$ is either $(1,t,t+1)$-branch or 
$(1,t,1)$-branch. The first case is impossible by Lemma~\ref{lem:xxx}, 
while the second case is impossible since, as we have shown in the previous 
paragraph, here the corresponding non-repeating branch should be 
a subtree of~$\cal V$, and this is not so. 

If $\cal C$ is a collection of leaves, then the transformation of 
Figure~\ref{fig:proof-rep-brush-2} shows that all leaves are of weight~1. 
Moreover, $\cal C$ contains not more than one leaf since 
otherwise we could transfer all the other leaves to the vertex $v$, 
which would change the tree. Therefore, in this case we get the tree~$Q$.
Finally, if $\cal C$ is not a collection of leaves, then $\cal W$ is a 
repeating branch of height at least~3, which, as above, is necessarily 
of type $(1,t,t+1)$, and Lemma~\ref{lem:xxx} implies that $\cal W$ has 
the form shown in Figure~\ref{fig:xxx} on the left, where the number 
of leaves is equal to one since otherwise we could transport all the 
leaves but one to the vertex $w$. Therefore, in this case we get the 
tree~$T$.

\ssk

Proposition~\ref{prop:repsts+2} is proved.
\hfill$\Box$

\subsection{Proof of the uniqueness of unitrees}\label{sec:uniqueness}
Our main tool will be ``cutting and gluing leaves'', though
these operations will be carried out not with the trees themselves
but with their passports; the trees must be kept in mind for an
intuitive understanding of the proof. We do not repeat every time
that ``the same reasoning remains valid if we interchange black and
white''.

\medskip

(A) There is only one black vertex (see Figure~\ref{fig:stars-A}); 
therefore, all white vertices must be adjacent to it, which means that 
they are leaves. The uniqueness is evident.

\msk

The uniqueness proofs for the cases from (B) to (E) all follow the same 
lines. If the structure of a passport implies the existence of a leaf of 
degree, say, $s$, then the only way to construct a corresponding tree 
is to glue an edge of the weight $s$ to a vertex of the opposite color 
of a tree having one less edge. Furthermore, in the initial (bigger) 
tree this edge can only be attached to a vertex of degree bigger than~$s$. 
If the smaller tree is a unitree (and usually it is by induction), and 
if there is essentially one way to attach the new edge to it, then the 
bigger tree is also a unitree. In certain cases more than one way of 
attaching a new edge may exist, but they all lead to isomorphic trees.

\msk

(B) Let us consider, for example, the case of an odd length, and examine
not the tree itself but its passport $(\al,\be)$. For this tree, $\al$
and $\be$ are the same: $\al=\be=((s+t)^k,s)$. The passport implies that 
there are two vertices of degree~$s$, a~black one and a white one, while 
all the other vertices, both black and white, are of degree $s+t$. 
These latter vertices cannot be leaves since otherwise there should exist 
vertices of a bigger degree to which such leaves would be attached. 
A tree must have at least two leaves. We conclude that there are exactly 
two leaves, and they are vertices of degree~$s$. They are connected to the 
tree by edges of the weight~$s$. The degree of the vertex to which such a 
leaf is attached is~$s+t$, and its color is opposite to the color of the 
leaf.

Now let us cut off one of these leaves, for example, the white one.
Then we get a tree with one less edge and with the passport 
$\al=((s+t)^{k-1},s,t)$ and $\be=(s+t)^k$. This passport corresponds to 
the chain tree of smaller (and even) length. We may inductively suppose 
that this tree is a unitree. The vertex of degree $t$ of this tree is a 
leaf. Now, we must make an operation which would simultaneously fulfill 
the following three goals:
\begin{itemize}
\item   it re-attaches back a white leaf of weight~$s$ to the smaller tree;
\item   it makes the black vertex of degree~$t$ in the smaller tree
        to disappear;
\item   it makes to appear an additional, $k$th black vertex of degree
        $s+t$, to the already existing $k-1$ ones.
\end{itemize}

It is clear that the only way to do all that is to attach this white leaf 
of weight~$s$ to the black vertex of degree~$t$.  This operation re-creates 
the initial chain-tree.

\ssk

The proof for an even length repeats the previous one almost word
to word, only the leaves are now of the same color and of degrees
$s$ and $t$. The base of induction is a tree consisting of a single
edge, which is obviously unique.

\msk

(C) The passport of a tree of the type $C$ is $\al=(ks+t,s^l)$, 
$\be=(ls+t,s^k)$. We affirm that there exists a leaf of degree~$s$.
Indeed, a tree must have at least two leaves, and the vertex of
the biggest degree cannot be a leaf. The biggest degree is either 
$ks+t$, or $ls+t$, or both. 

Suppose that we have a black leaf of degree $s$. Then it has to be attached 
to the only white vertex of degree bigger than $s$, which is the white 
vertex of degree $ls+t$. Cut this leaf off. We get a smaller tree, with 
one less edge, with $l$ being replaced with $l-1$. This tree is a smaller
instance of $C$ which may supposed to be a unitree by induction. 
(Note that in particular cases it can also be of type $B$, or even of $A$, 
the latter one when $l$ was equal to~1.)

Now we no longer work with the passports but with the trees. We know 
the smaller tree since it is unique, and we must re-attach the previously 
cut-off black leaf to a white vertex of this smaller tree. Here two cases 
may take place.
\begin{enumerate}
\item   If $s\ne t$, or even if $s=t=1$ but $l\ne 1$, the initial (bigger) 
        tree did not have a white vertex of degree~$2s$. Therefore, we 
        cannot attach the cut-off leaf to a white vertex of degree $s$. 
        Hence, the only vertex to which it can be attached is the white 
        vertex of degree~$(l-1)s+t$.
\item   If $s=t=1$ and $l=1$ then the smaller tree is the star with all
        its leaves being of degree~1. Then we may re-attach the leaf
        to any one of them, the resulting tree will be the same.
\end{enumerate}

There is an additional subtlety here. The planar structure of our trees
means that we must choose not only a vertex to which we attach a new
edge. We must also choose an angle between neighboring edges incident 
to the vertex of attachment, and to insert the leaf into the angle between
these edges. If there are $m$ edges incident to a vertex, there also are 
$m$ angles between them, and therefore $m$ ways of placing the new edge. 
But, obviously, in our case all these ways give the same plane tree, see 
Figure~\ref{fig:brushes-C}.

\msk

(D) Two black vertices of degree $2s+t$ cannot be leaves. Therefore,
there exists a leaf of degree $s$ or $s+t$. Cut it off, and we get
either $C$ or $E_1$. Indeed, if the cut-off (white) leaf was of degree~$s$, 
and was attached to one of the black vertices of degree $2s+t$ (there
are no other black vertices), then the passport of the smaller tree
becomes $\al=(2s+t,s+t)$, $\be=((s+t)^2,s)$. This passport corresponds
to the pattern $E_1$, with $l=1$ and the length of the chain equal to~4. 
The uniqueness of the corresponding tree will be proved in a moment. 
The only way to glue a leaf of degree~$s$ to this tree and to create a 
vertex of degree~$2s+t$ instead of vertex of degree~$s+t$ is to glue 
this leaf to the vertex of degree~$s+t$.

If, on the other hand, the cut-off (white) leaf was of degree~$s+t$, and 
was attached to one of the black vertices of degree~$2s+t$, then the
passport of the smaller tree becomes $\al=(2s+t,s)$, $\be=(s+t,s^2)$.
This passport corresponds to a tree of type $C$, with $k=2$ and $l=1$.
The uniqueness of such a tree was proved above. Then the only way to glue 
back the leaf of degree~$s+t$ is to glue it to the black vertex of
degree~$s$ of the smaller tree.

It is easy to see that in both cases we get the same tree $D$.


\msk

(E) The proof is similar to the cases considered above, so we will
shorten our presentation. Consider first the cases $E_3$ and $E_4$. 
All the vertices except
two are of degree $s+t$; the two remaining ones are of degrees
$(k+1)s+kt$ and $(l+1)s+lt$ for $E_3$, and $(k+1)s+kt$ and $ls+(l+1)t$
for $E_4$. Without loss of generality we may suppose that $(k+1)s+kt$
is the bigger of the two; therefore, it cannot be a leaf. For $E_4$,
the ``second best'' vertex cannot be a leaf either since it has the
same color. For $E_3$, if $k>l$, the vertex of degree $(l+1)s+lt$ might 
in principle be a leaf. Whatever is the case, there exists a leaf of 
degree $s+t$. Cut it off, and we obtain a smaller tree, with the possible 
pattern transitions as follows: $E_4\to E_4$; $E_3\to E_3$; $E_4\to E_2$; 
or $E_3\to E_1$, the latter two maybe with renaming the variables.

Now, for the cases $E_1$ and $E_2$ the situation is similar. All the
vertices except two are of degree $s+t$. The vertex of the biggest degree
cannot be a leaf. Therefore, there exists a leaf of degree $s$ or $s+t$.
Cut it off, and we get a smaller tree,  with the possible pattern 
transitions as follows: $E_1\to E_1$; $E_1\to E_2$; $E_2\to E_1$;
$E_2\to E_2$, or we may arrive to the patterns $A$ or~$B$.

What remains now is to see that there is only one way to re-attach
the cut-off leaf to the smaller unitree.

\msk

(F,\,H,\,I,\,Q) The trees $F$, $H$, $I$, $Q$ are ordinary; therefore, the 
enumerative formula (\ref{eq:GJ}) can be applied. 

If $m\ne l$, a tree of the series $F$ is asymmetric, and therefore its 
contribution to (\ref{eq:GJ}) is~1. Now, formula (\ref{eq:GJ}) in this 
case gives 1; therefore, there is no other tree with this passport. 

When $m=l$, a tree of the series $F$ is symmetric, with the rotational 
symmetry of order $k$. Therefore, its contribution to (\ref{eq:GJ}) 
is~$1/k$. Now, the formula itself gives $1/k$; therefore, there is no 
other tree in this case either.

For the trees of the series $H$, formula (\ref{eq:GJ}) gives 1 when 
$k\ne l$, and gives $1/2$ when $k=l$. This corresponds to the symmetry
order of these trees: they are asymmetric when $k\ne l$, and symmetric
of order~2 when $k=l$.

The trees of the series $I$ are asymmetric, and formula (\ref{eq:GJ}) 
gives~1.

The tree $Q$ is asymmetric, and formula (\ref{eq:GJ}) gives 1.

\msk

(G) The tree has $km$ vertices and hence $km-1$ edges. Since the total 
weight is $km$, there exists exactly one edge of weight 2 while all the 
other edges are of weight 1. The only white vertex to which the edge of 
weight 2 can be attached is the vertex of degree $k$, since all the other 
white vertices are of degree 1. The rest is obvious.

\msk

(J) All white vertices are of degree 2; therefore, a weight of an 
edge can only be 1 or 2. There are only three black vertices, their degrees 
being 4, $2k+1$, $2k+1$. Therefore, the black vertices cannot be leaves 
since such leaves could not be attached to a white vertex of degree~2; 
thus, all the leaves are white. A white vertex which is not a leaf must 
have two black neighbors; therefore, there are exactly two white vertices 
which are not leaves: they are ``intermediate'' white vertices between 
the black ones. The edges incident to them are both of weight 1. The black 
vertex of degree 4 cannot have two incident edges of weight 2 since these
edges should be leaves, and such a tree would not be connected; it cannot 
have four incident edges of weight 1 either since such a tree should need 
more than three black vertices. Therefore, the weights of the edges attached 
to this vertex must be 2, 1, 1, and the edge of weight~2 is a leaf. The 
rest is obvious.

\msk

(P) All black vertices are of degree 5, all white ones are of degree 3. 
Therefore, black vertices cannot be leaves. Arguing now as in the case (J) 
we conclude that there are exactly three white leaves, implying easily 
that there exists only one tree with this passport.  


\msk

(K,\,L,\,M,\,N,\,O,\,R,\,S,\,T) The proof of all these cases follows 
the same pattern.

Let us take, for example, the tree $O$. Its passport is $(5^4,2^{10})$.
Therefore, the number of vertices is $4+10=14$ and the number of edges
is 13, while the total weight is $5\cdot 4 = 2\cdot 10 = 20$. Thus,
the overweight is 7, and it must be distributed among the edges.

Now, no edge can have a weight greater than 2 since the degrees of
white vertices are all equal to 2. Therefore, the tree $O$ has seven
edges of weight 2. Moreover, {\em all of them are leaves}\/; indeed, 
if something were attached to the white end of such an edge then this 
white end would have a degree greater than~2.

The same reasoning may be carried out for all the above cases, with 
their respective overweights and numbers of leaves of weight 2.

Now, let us cut off all the leaves of weight 2. What remains is an
ordinary tree, and we must verify that it is a unitree. Usually it
is immediately obvious since the ordinary tree in question is very
small; otherwise, we may apply formula (\ref{eq:GJ}), or else we may
remark that such a tree belongs to one of the previously established 
cases. For example, for the tree $T$ what remains after cutting off 
the leaves of weight~2 is the tree $Q$.

The last step consists in proving that there is only one way to glue
back to this ordinary unitree the leaves of weight 2 that were previously
cut off. For example, in the case $O$ the ordinary unitree has black
vertices of degrees 3,~1,~1,~1, and we have, by gluing to them seven
edges of weight 2, made these degrees equal to 5, 5, 5, 5. Obviously,
there is only one way to do that. In fact, in certain cases there are
several ways of gluing but they give the same result because of a
symmetry of the underlying ordinary tree. For the tree $L$, there is
an additional condition we must satisfy: white leaves can only be attached
to black vertices.

\msk

Theorem \ref{th:main} is proved. 
\hfill $\Box$ \label{end-proof}

\section{Other combinatorial Galois invariants}
\label{sec:galois}

The theory of dessins d'enfants studies combinatorial invariants
of the Galois action on dessins. These invariants have various
levels of generality. The most general (and the most simple) one is 
the passport; the subject of this paper is precisely the case when 
the passport alone guarantees the definability of a dessin over $\Q$. 
But our exposition would be incomplete if we did not mention several
other Galois invariants which lead to further examples of dessins and
DZ-triples defined over $\Q$.

\newpage

\subsection{Composition}

The following proposition is obvious.

\begin{proposition}[Composition]\label{prop:composition}
Let $f=f(x)$ and $h=h(t)$ be two rational functions such that:
\begin{itemize}
\item   $f$ is a Belyi function, with the corresponding dessin $D_f$;
\item   $h$ is a function $($not necessarily a Belyi one\/$)$, all of 
        whose critical values are either vertices or face centers of $D_f$.
\end{itemize}
Then the function $F(t)$ obtained as a composition
$$
F(t)=f(h(t)), \qquad \mbox{\rm that is}, \qquad
F\,:\,\Cbar\,\stackrel{h}{\fleche}\,\Cbar\,\stackrel{f}{\fleche}\,\Cbar,
$$ 
is a Belyi function. If, furthermore, both $f$ and $h$ are defined 
over~$\Q$, then, obviously, the same is true for $F$.
\end{proposition}

The above proposition gives us a very general method of constructing 
Belyi functions with all its finite poles being simple, or, in other
words, Belyi functions corresponding to weighted trees.

\begin{corollary}[Decomposable weighted trees]\label{cor:composition}
Suppose that the functions $f$ and $h$ of the above proposition
satisfy the following properties:
\begin{itemize}
\item   the dessin~$D_f$ corresponds to a weighted tree, that is, all its
        finite faces are of degree~$1$;
\item   $h$ is a polynomial all of whose critical values except 
        infinity are vertices of $D_f$.
\end{itemize}
Then all the finite faces of the dessin $D_F$ corresponding to the Belyi 
function $F(t)=f(h(t))$ are of degree~$1$. If, furthermore, both $f$ and $h$ 
are defined over\/~$\Q$, then, obviously, the same is true for~$F$.
\end{corollary}

\paragraph{Proof.} Since $h$ is a polynomial, the only poles of $F=f\circ h$, 
except infinity, are the preimages of the simple poles of~$f$, i.\,e.,
the preimages of the centers of the small faces of~$D_f$. Since $h$ is 
not ramified over these simple poles, they remain simple for $F$, and 
each of them is ``repeated'' $\deg h$ times.
\hfill $\Box$

\begin{example}[Composition 1]\label{ex:composition-1}
Consider the following functions:
$$
f \eq -\,\frac{64\,x^3(x-1)}{8\,x+1}\,, \qquad
u \eq \frac{1}{5^5}\cdot(t^2+4)^3(3\,t+8)^2\,.
$$
Here $f$ is a Belyi function corresponding to the upper left dessin 
in~Figure~\ref{fig:composition}, and $u$ is a Belyi function corresponding
to the lower left dessin. 

Substituting $x=u(t)$ in $f$ we obtain a Belyi function $F$ corresponding
to the dessin shown on the right of Figure~\ref{fig:composition}.
It is obvious that the {\em combinatorial}\/ orbit of the dessin $D_F$ 
consists of more than one element: for example, the petals attached to 
the vertices of degrees 9 and 6 can be cyclically arranged in many 
different ways. Still, $F\in\Q(x)$ by construction.

Note that the dessins $D_f$ and $D_u$ serving as building blocks
for the above example both belong to the classification we have
established in Section~\ref{sec:unitrees}: they both correspond to
unitrees, and it is their passports that guarantee that they are defined 
over~$\Q$. We don't have a simple way of constructing more general
examples over~$\Q$ (when, e.\,g., the polynomial $u$ has more than two 
critical values with prescribed positions at vertices of $D_f$).

\begin{figure}[htbp]
\begin{center}
\epsfig{file=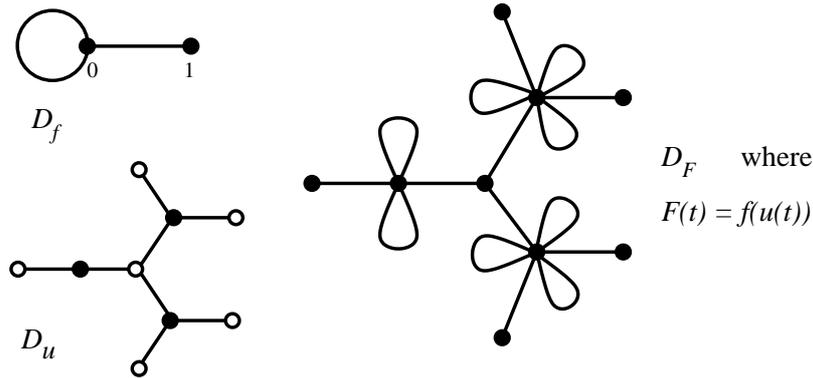,width=10.8cm}
\caption{\small The pictures corresponding to $f(x)$ and $F(t)$ are drawn 
according to Convention \ref{con:without-white}: only their black vertices
are shown explicitly. In the picture corresponding to $u$ the black 
vertices are those sent to~0, and the white ones are those sent to 1.}
\label{fig:composition}
\end{center}
\end{figure}

\end{example}

\begin{example}[Composition 2]\label{ex:composition-2}
Another example, based on the same function $f$, is as follows. We have
$$
f-1 \eq -\,\frac{(8\,x^2-4\,x-1)^2}{8\,x+1}\,,
$$
so the white vertices of the dessin $D_f$ (which are not shown 
explicitly in Figure~\ref{fig:composition}) are the 
roots of $8\,x^2-4\,x-1$, that is, they are equal to $(1\pm\sqrt{3})/4$. 
Now, the critical values of the polynomial
$$
v \eq \frac{1}{3}\,t^3-\frac{3}{4}\,t+\frac{1}{4}\,,
$$
that is, the values of $v$ at the roots of $v'=t^2-3/4$, are equal to 
exactly $(1\pm\sqrt{3})/4$. Therefore, the composition $G(t)=f(v(t))$ 
is once again a Belyi function, and all its poles except infinity are 
simple. The corresponding dessin is shown in Figure~\ref{fig:composition-2}.

\begin{figure}[htbp]
\begin{center}
\epsfig{file=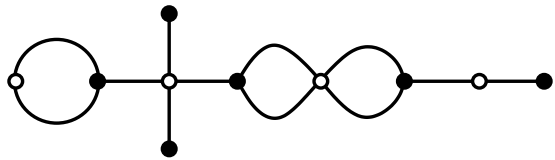,width=7.2cm}
\caption{\small The dessin $D_G$ corresponding to the function 
$G(t)=f(v(t))$; this time not all white vertices are of degree 2, 
therefore we show them explicitly.}
\label{fig:composition-2}
\end{center}
\end{figure}

It is obvious that the dessin $D_G$ is not the only one having the 
passport $(3^31^3,4^22^2)$. For example, the two ``vertical'' edges
can both be put above, or both can be put below the horizontal axis, 
or they can be cut off and attached to the leftmost white vertex of 
degree~2 (the one on the loop), or to the rightmost one (the one on the
horizontal segment). Nevertheless, the dessin thus obtained is defined 
over $\Q$ by construction.
\end{example}

Now the dessins $D_F$ and $D_G$, being defined over $\Q$, may themselves 
serve for a similar construction: if, for example, $w$ is a polynomial
with coefficients in $\Q$ whose critical values are vertices of $D_G$,
then the function $H=G\circ w=f\circ v\circ w$ is a Belyi function
corresponding to a dessin $D_H$, all of whose finite faces are small.
In such a composition, only $f$ has to be a Belyi function while the 
subsequent terms may have more than three critical values.

\begin{remark}[Symmetric trees]
The group of the orientation preserving  automorphisms of a plane tree 
is always cyclic. If it is $\Z_k$
then the Belyi function for the corresponding map is $F(x)=f(x^k)$
where $f$ is the Belyi function for the map corresponding to a
single branch of the tree (the vertex of this branch, which will become
the center of the symmetric tree, must be put to the origin). Among the
unitrees classified in Section~\ref{sec:unitrees}, the trees $N$ and
and $R$ are symmetric (of order 3 and 2 respectively). Some elements
of the infinite series may also be symmetric (for special values of 
parameters).
\end{remark}

We leave it to the reader to see that the series $H$ and $I$ are
compositions (although the composition in this case is not reduced to 
a rotational symmetry), and that a multiplication of all the weights 
of the edges of a tree by a factor $d$ can be represented as a 
composition with the following Belyi function:
$$
f(x)\eq\frac{x^d}{x^d-(x-1)^d}\,, \qquad \mbox{hence} \qquad 
f(x)-1\eq\frac{(x-1)^d}{x^d-(x-1)^d}\,.
$$

\subsection{Primitive monodromy groups}
\label{sec:primitive}

\begin{definition}[Primitive and special groups]\label{def:primitive}
A permutation group of degree $n$ acting on a set $X$, $|X|=n$, is
called {\em imprimitive}\/ if the set $X$ can be subdivided into
$m$ disjoint {\em blocks}\/ $X_1,\ldots,X_m$ of equal size $|X_i|=n/m$,
where $1<m<n$, such that an image of a block under the action of any 
element of the group is once again a block. A permutation group which 
is not imprimitive is called {\em primitive}. A primitive permutation 
group not equal to ${\rm S}_n$ or ${\rm A}_n$ is called {\em special}.
\end{definition}

It is known (see \cite{LanZvo-04}) that a covering is a composition of 
two or more coverings of smaller degrees {\em if and only if}\/ its 
monodromy group is imprimitive. Thus, a covering which is not a composition 
has a primitive monodromy group. However, in a vast majority of cases this
group is equal to either ${\rm S}_n$ or ${\rm A}_n$, for a very simple
reason: a permutation group generated by a randomly chosen pair of
permutations is  ${\rm S}_n$ or ${\rm A}_n$ with a probability close
to~1. This is why special groups are of a particular interest: since 
{\em the monodromy group is a Galois invariant}, such a group gives an 
additional invariant in the absence of the composition.

In the case of weighted trees, that is, in the case of coverings realized
by Belyi functions with all poles except one being simple, the monodromy 
group must contain a permutation of the cycle structure $(n-r,1^r)$: 
it is the monodromy permutation corresponding to a loop around infinity
on the Riemann sphere. Motivated by our study of weighted
trees, Gareth A. Jones classified all special permutation groups 
containing such a permutation, see \cite{Jones-12} (2012). 
In particular, it is shown that in all such cases $r\le 2$.
This property is based on two results. The first is an old theorem 
by Jordan~(1871)~\cite{Jordan-1871} stating that a primitive group
containing a permutation with the cycle structure $(n-r,1^r)$ is
$(r+1)$-transitive. The second is the complete list of multiply
transitive groups: it is based on the classification theorem of
finite simple groups.

The classification due to Jones looks as follows (we use standard
notation for projective, cyclic and affine groups and for the Mathieu
groups):

\begin{theorem}[G.\,Jones's classification]\label{th:jones}
Let $G$ be a primitive permutation group of degree $n$ not equal to\/
${\rm S}_n$ or ${\rm A}_n$. Suppose that $G$ contains a permutation
with cycle structure $(n-r,1^r)$. Then $r\le 2$, and one of the
following holds:
\begin{itemize}
  \item[\rm 1.] $r=0$ and either
    \begin{itemize}
      \item[\rm (a)] ${\rm C}_p \le G \le {\rm AGL}_1(p)$ 
        with $n=p$ prime, or
      \item[\rm (b)] ${\rm PGL}_d(q) \le G \le {\rm P{\Gamma}L}_d(q)$
        with $n=(q^d-1)/(q-1)$ and $d\ge 2$ for some prime power $q=p^e$, or
      \item[\rm (c)] $G={\rm L}_2(11)$, ${\rm M}_{11}$ or\/ ${\rm M}_{23}$
        with $n=11$, $11$ or\/ $23$ respectively.  
    \end{itemize}
  \item[\rm 2.] $r=1$ and either
    \begin{itemize}
      \item[\rm (a)] ${\rm AGL}_d(q) \le G \le {\rm A{\Gamma}L}_d(q)$
        with $n=q^d$ and $d\ge 1$ for some prime power $q=p^e$, or
      \item[\rm (b)] $G={\rm L}_2(p)$ or\/ ${\rm PGL}_2(p)$ with
        $n=p+1$ for some prime $p\ge 5$, or
      \item[\rm (c)] $G={\rm M}_{11}$, ${\rm M}_{12}$ or\/ ${\rm M}_{24}$
        with $n=12$, $12$ or\/ $24$ respectively.
    \end{itemize}
  \item[\rm 3.] $r=2$ and ${\rm PGL}_2(q) \le G \le {\rm P{\Gamma}L}_2(q)$
    with $n=q+1$ for some prime power $q=p^e$.
\end{itemize}
\end{theorem}

\pagebreak[4]

\begin{example}[A tree with a special monodromy group]\label{ex:special}
There are six weighted trees of degree $n=8$ with the passport 
$(7^11^1,2^31^2,6^11^2)$. One may expect that their common moduli field
would be an extension of $\Q$ of degree~6. However, this is not the case.
Five trees out of six have the monodromy group~${\rm S}_8$, while the 
remaining one, shown in Figure~\ref{fig:primitive-1} (left) has
the monodromy group~${\rm PGL}_2(7)$. Therefore, this tree is defined
over~$\Q$. The other trees form a single Galois orbit over a field 
of degree~5.

\begin{figure}[htbp]
\begin{center}
\epsfig{file=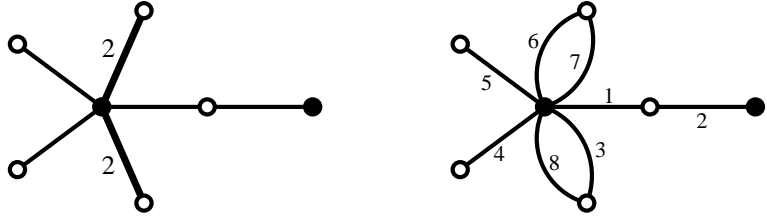,width=10cm}
\caption{\small The tree on the left (equivalently, the map on the
right) has the monodromy group ${\rm PGL}_2(7)$. All the other trees
with the same passport have the monodromy group ${\rm S}_8$.}
\label{fig:primitive-1}
\end{center}
\end{figure}

In order to establish that the group in question is indeed ${\rm PGL}_2(7)$
we may proceed as follows. First, we draw the bicolored map represented
by this tree, as it is done in  Figure~\ref{fig:primitive-1}, right.
Then, we label the edges of the map and write down two permutations:
the first one represents the cyclic order of the edges around its black 
vertices in the counterclockwise direction, while the second one
represents the cyclic order of the edges around their white vertices in 
the same direction. In our case these permutations are
$$
a \eq (1,7,6,5,4,8,3), \qquad b \eq (1,2)(3,8)(6,7).
$$
The permutation corresponding to the faces is
$$
c \eq (ab)^{-1} \eq (1,2,3,4,5,6),
$$
so that $abc=1$. The cycle $c$ can be read in the picture by going around 
the outer face. Note that the cycle structures of $a$, $b$ and $c$ are 
$7^11^1$, $2^31^2$ and $6^11^2$, respectively. Then, the monodromy group is
$$
G \eq \langle a,b \rangle \eq \langle a,b,c \rangle.
$$
Using Maple it is easy to find out that $|G|=336$, and the only transitive
subgroup of ${\rm S}_8$ of order~336 is ${\rm PGL}_2(7)$ (see, for
example, \cite{ButMcK-83}). For the other five trees with the same
passport the same Maple package shows that the size of their monodromy
group is $40\,320=8!$, so the group in question is~${\rm S}_8$.
\end{example}

\begin{example}[${\rm PGL}_2(7)$ once again]\label{ex:special-bis}
One more example is shown in Figure~\ref{fig:primitive-2}. There are
five trees with the passport $(6^11^2,3^21^2,6^11^2)$. One of them is
symmetric and therefore forms a Galois orbit containing a single element
and thus defined over $\Q$. Three trees have the monodromy group~${\rm S}_8$;
they form a cubic Galois orbit. Finally, the remaining tree shown in  
Figure~\ref{fig:primitive-2} has the monodromy group ${\rm PGL}_2(7)$. 
Therefore, it forms a Galois orbit in itself and is thus defined over~$\Q$.

\begin{figure}[htbp]
\begin{center}
\epsfig{file=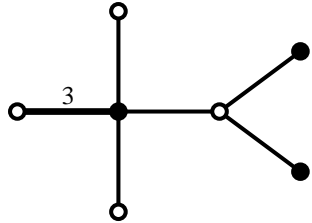,width=4cm}
\caption{\small One more tree with the monodromy group ${\rm PGL}_2(7)$.}
\label{fig:primitive-2}
\end{center}
\end{figure}

\end{example}

\subsection{Duality and self-duality}
\label{sec:duality}

A dual to a map is usually constructed as follows. First, one puts a new 
vertex inside every face of the initial map: this vertex is called ``center''
of the face. Then, the centers of the adjacent faces are connected by
edges in such a way that every edge of the initial map is crossed in
its ``middle point'' by a new edge. A dual of a dual is the initial map.

\begin{figure}[htbp]
\begin{center}
\epsfig{file=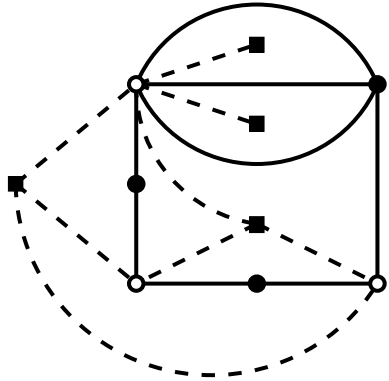,width=5cm}
\caption{\small A bicolored map (in unbroken line), and its dual (in 
dashed line). The black vertices of the dual map are designated by 
squares. The white vertices belong to both maps.}
\label{fig:dual}
\end{center}
\end{figure}

For the {\em bicolored maps}\/ a specific variant of the above construction
is used, when only black vertices are considered as vertices, while the
white vertices play the role of the edge midpoints. An association is 
thus made between the faces of the initial map and {\em black vertices}\/ 
of the dual map. The white vertices belong to both maps. More exactly, a
center of a face is connected by edges with all the white verices lying
on the border of the face: see an example in Figure~\ref{fig:dual} where
the initial map is shown in an unbroken line, and its dual, in a dashed 
line; the black vertices of the dual map are designated by the little 
squares.

From the point of view of Belyi functions, if $f(x)$ is a Belyi function
for the initial map then $1/f(x)$ is a Belyi function for its dual.
Indeed, $1/y$ interchanges 0 and $\infty$ while leaving 1 untouched.
Therefore, the former poles become roots (i.\,e., black vertices),
and vice versa.

\begin{definition}[Self-dual map]\label{def:autodual}
A bicolored map is called {\em self-dual}\/ if it is isomorphic to its
dual map. 
\end{definition}

Of course, the fact that a map is self-dual does not mean that $f=1/f$
where $f$ is its Belyi function. It means that $1/f(x)=f(w(x))$ where
$w(x)$ a linear fractional transformation of the variable~$x$. 
The self-duality is an invariant of the Galois action: if a function 
satisfies an algebraic relation while the other function does not satisfy 
the same relation, they cannot belong to the same Galois orbit.

\begin{figure}[htbp]
\begin{center}
\epsfig{file=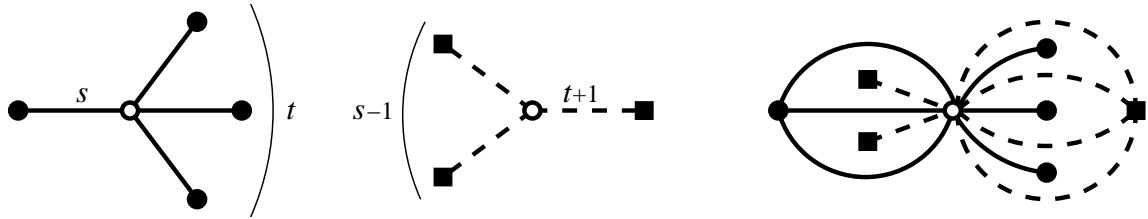,width=15cm}
\caption{\small Two branches of weighted trees of diameter 4 dual to each 
other. The figure on the right shows how these branches, represented as 
maps, fit to one another.}
\label{fig:dual-branches}
\end{center}
\end{figure}

A weighted tree represents a map whose all faces except one are of
degree~1. Therefore, its dual map must have all its black vertices except
one being of degree~1. This can only happen if the dual map corresponds
to a weighted tree of diameter~4: it has a black vertex of a degree greater
than~1 (its central vertex), while all its black leaves are of degree~1. 
Therefore, if we are interested in self-dual maps which correspond to 
weighted trees then we must consider only the trees of diameter~4. 
The condition on the branches of such trees in order for them to be dual 
to each other is shown in Figure~\ref{fig:dual-branches}.

Now we are ready to give an example where the self-duality plays the
role of a Galois invariant.

\begin{example}[Self-duality as a Galois invariant]
Let us take two integers $p$ and $q$, $p<q$, and consider the
following passport of degree $n=2p+2q-2$:
\begin{itemize}
\item   there is a black vertex of degree $p+q$ (the center), and $p+q-2$ 
        black vertices of degree~1 (the leaves);
\item   there are two white vertices, of degrees $2p-1$ and $2q-1$
        respectively;
\item   the above data imply that the trees have $p+q$ edges, and therefore 
        the outer face is of degree $p+q$, 
        the same as the degree of the central black vertex.
\end{itemize}

There are exactly $2p-1$ trees with this passport. Their general appearance
is shown in Figure~\ref{fig:autodual}. Here the parameters take the
following values: $s=1,2,\ldots,2p-1$ while 
$$
t=(p+q)-s, \quad k=(2p-1)-s, \quad l=(2q-1)-t.
$$
Among all these trees, only one is self-dual: it corresponds to the
values $s=p$, $t=q$, $k=p-1$, and $l=q-1$. Therefore, this tree is 
defined over $\Q$.

(In this example both branches are dual to themselves. An attempt to
make one branch dual to the other leads to the equality $p=q$, but
we have supposed that $p<q$.)

\begin{figure}[htbp]
\begin{center}
\epsfig{file=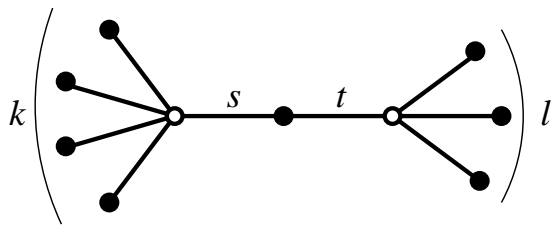,width=7.2cm}
\caption{\small Here $s+t=p+q$, $s+k=2p-1$, $t+l=2q-1$, where $1<p<q$.
The combinatorial orbit consists of $2p-1$ trees but it splits into
at least two Galois orbits since exactly one of these trees is self-dual,
the one with $s=p$ and $t=q$.}
\label{fig:autodual}
\end{center}
\end{figure}

\end{example}

\begin{remark}[Example~\ref{ex:special-bis} revisited]
All the five trees with the passport $(6^11^2,3^21^2,6^11^2)$ considered 
in Example~\ref{ex:special-bis} are self-dual. Therefore, for them the
self-duality cannot serve as an additional Galois invariant leading to
a splitting of the combinatorial orbit into two (or more) Galois orbits.
\end{remark}

\subsection{A sporadic example}\label{sec:galois-split}

\begin{figure}[htbp]
\begin{center}
\epsfig{file=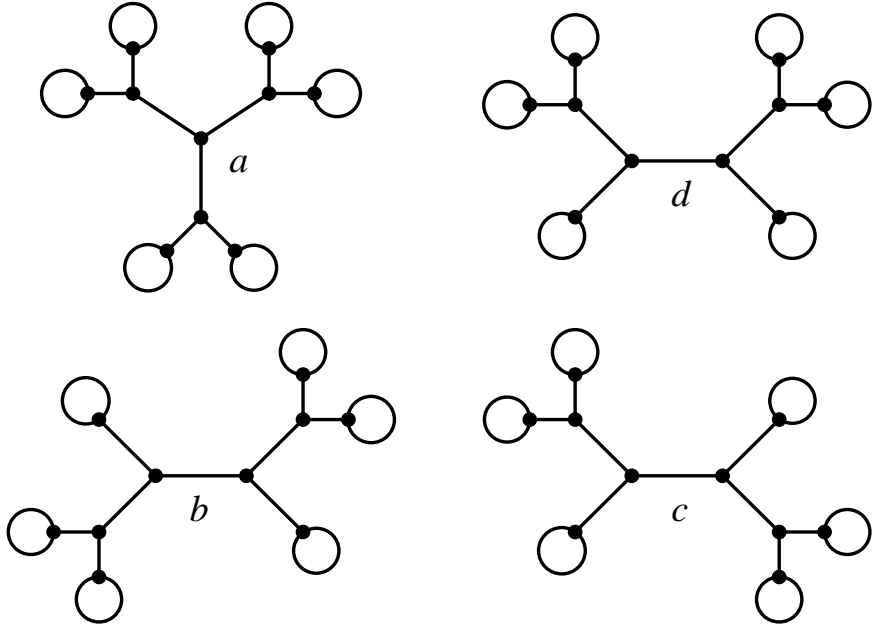,width=11.4cm}
\caption{\small This combinatorial orbit, corresponding to the passport
$(\al,\be,\ga)$ where $\al=3^{10}$, $\be=2^{15}$, $\ga=24^11^6$, splits 
into three Galois orbits: $\{a\}$, $\{b,c\}$, and $\{d\}$. The dessins
$a$ and $d$ are defined over $\Q$.}
\label{fig:zannier-5}
\end{center}
\end{figure}

The world of dessins d'enfants is rich with various specific cases. 
Let us consider, 
for example, the set of dessins shown in Figure~\ref{fig:zannier-5}. 
They constitute a combinatorial orbit for the passport $(\al,\be,\ga)$
where $\al=3^{10}$, $\be=2^{15}$, and $\ga=24^11^6$. We might na\"{\i}vely
suppose that this combinatorial orbit also constitutes a Galois orbit; 
if this were the case, this orbit would be defined over a field of 
degree~4 (since it has four elements). However, the reality is more 
complicated and, in fact, more exciting. 

Namely: the dessin $a$ is the only one having a rotational symmetry 
of order~3 around a black vertex. Therefore, the singleton $\{a\}$ 
constitutes a Galois orbit. Two dessins $b$ and $c$ are the only ones 
which have rotational symmetry of order~2, the center being a white 
vertex (we recall that the white vertices, being all of degree~2, 
are not shown explicitly in the picture). Therefore, the set $\{b,c\}$ 
must also be taken apart from the combinatorial orbit. There are two
{\em a priori}\/ possibilities: $b$ and $c$ may make two Galois orbits, 
both defined over $\Q$, or they may make a single Galois orbit defined 
over a quadratic field. But any map whose Belyi function is defined over
a real field must be axially symmetric since it remains invariant under 
the complex conjugation. This observation excludes the possibility of
two orbits over $\Q$, and it also excludes a real quadratic field.
We may conclude that the set $\{b,c\}$ constitutes a single Galois orbit
defined over an imaginary quadratic field.

The dessin $d$ is not symmetric and does not have any other specific 
combinatorial features. But it remains solitary, and therefore it 
constitutes a Galois orbit all by itself. Since the orbits $\{a\}$ 
and $\{d\}$ consist of a single element, their Belyi functions are 
defined over $\Q$. Thus, the dessin~$d$ is defined over $\Q$ for no 
other reason than the fact that {\em it remains alone after all the 
other Galois orbits are taken away}\/.

This combinatorial orbit is markworthy for the reason that different authors 
returned to it, or to some of its elements, many times. The Belyi function 
for $a$ was computed by Birch and already appeared in~\cite{BCHS-65} (1965); 
the one for $d$ was computed 35~years later by Elkies~\cite{Elkies-00} 
(2000). All the four Belyi functions were independently computed by Shioda 
\cite{Shioda-04} (2004). In particular, he found out that the orbit
$\{a,b\}$ is defined over the field $\Q(\sqrt{-3})$. Shioda had already
used as a starting point the above combinatorial orbit; the other
authors have apparently made a ``blind'' search.

Our combinatorial approach does not make the computational part of the 
work any easier. Its advantage is elsewhere. It consists in the fact that, 
before any computation, we may be sure of the following.
\begin{itemize}
\item   There exist exactly four non-equivalent Belyi functions with
        the passport $(3^{10},2^{15},24^11^6)$; here ``non-equivalent''
        means that they cannot be obtained from one another by a linear
        fractional change of variables.
\item   Belyi functions corresponding to $a$ and $d$ are defined over~$\Q$.
\item   Belyi function corresponding to $a$ is a rational function in $x^3$
        (because of the threefold symmetry of the dessin $a$).
\item   Belyi functions for the orbit $\{b,c\}$ are defined over an
        imaginary quadratic field.
\end{itemize}

More examples similar to this one are given below.

\subsection{Sporadic examples of Beukers and Stewart \cite{BeuSte-10}}
\label{sec:BeuSte}

All the examples in this section are borrowed from the 
paper~\cite{BeuSte-10} by Beukers and Stewart, which was one of the 
sources of inspiration for our study. In their paper, the authors
consider only the case of powers of polynomials. Namely, they look for 
polynomials $A$ and $B$, {\em defined over}\/ $\Q$, for which the degree 
$\deg\,(A^p-B^q)$ attains its minimum. The degrees of polynomials
in question are $\deg A=qr$, $\deg B=pr$ where the parameter $r$ may
be greater than~1. The passport of the corresponding tree is
$(p^{qr},q^{pr})$.

The authors find, as we do, several infinite series of DZ-triples 
(which they call Davenport pairs), and several sporadic examples.
The first such example, for which $(p,q,r)=(5,2,2)$, corresponds to
our sporadic tree $O$. The next one, $(p,q,r)=(5,3,1)$, corresponds
to the sporadic tree $P$. However, the subsequent examples do not
correspond to anything we have found up to now. What is going on?

It turns out that here we encounter once again the phenomenon that we
already explained in Section~\ref{sec:galois-split}.

\begin{example}[$(p,q,r)=(7,3,1)$]
There exist two trees corresponding to the passport $(7^3,3^7)$: they
are shown in Figure~\ref{fig:split-7-3}. We see that one of the trees
is symmetric, with the symmetry of order~3, while the other one is not.
Therefore, this combinatorial orbit splits into two Galois orbits,
and hence both trees are defined over $\Q$. The left-hand one corresponds
to the example given in~\cite{BeuSte-10}.

\begin{figure}[htbp]
\begin{center}
\epsfig{file=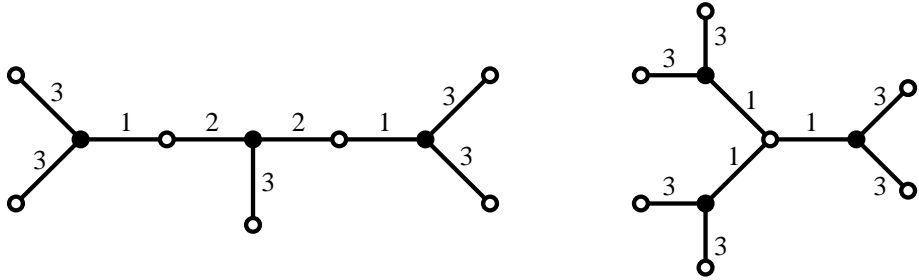,width=12cm}
\caption{\small Two trees corresponding to the passport $(7^3,3^7)$; one
of them is symmetric, the other one is not.}
\label{fig:split-7-3}
\end{center}
\end{figure}

\end{example}

Note that an axial symmetry is not a Galois invariant.

\begin{example}[$(p,q,r)=(8,3,1)$ and $(10,3,1)$]\,
The\, situation\, for\, the\, passports\, $(8^3,3^8)$\, and\, 
$(10^3,3^{10})$ is similar to the previous one.
For the first passport there are two trees, and one of them is symmetric 
while the other is not (see Figure~\ref{fig:split-8-3}); therefore, both 
are defined over~$\Q$. For the second passport there are three trees (see 
Figure~\ref{fig:split-10-3}). One of them is symmetric with the symmetry 
of order~2; one is symmetric with the symmetry of order~3; and one is 
asymmetric. Therefore, all the three trees are defined over $\Q$. In both 
cases ``sporadic'' polynomials given in~\cite{BeuSte-10} correspond to 
asymmetric trees.

\begin{figure}[htbp]
\begin{center}
\epsfig{file=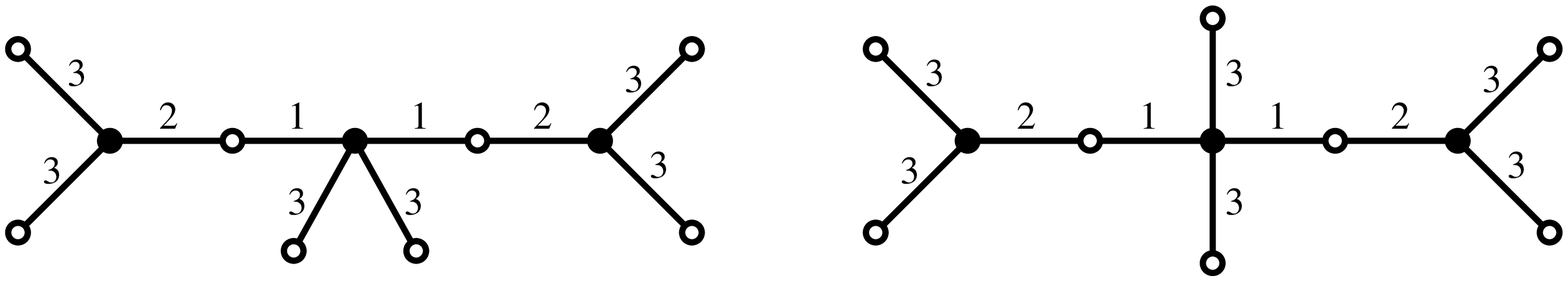,width=13.5cm}
\caption{\small Two trees corresponding to the passport $(8^3,3^8)$.} 
\label{fig:split-8-3}
\end{center}
\end{figure}

\pagebreak[4]

\begin{figure}[htbp]
\begin{center}
\epsfig{file=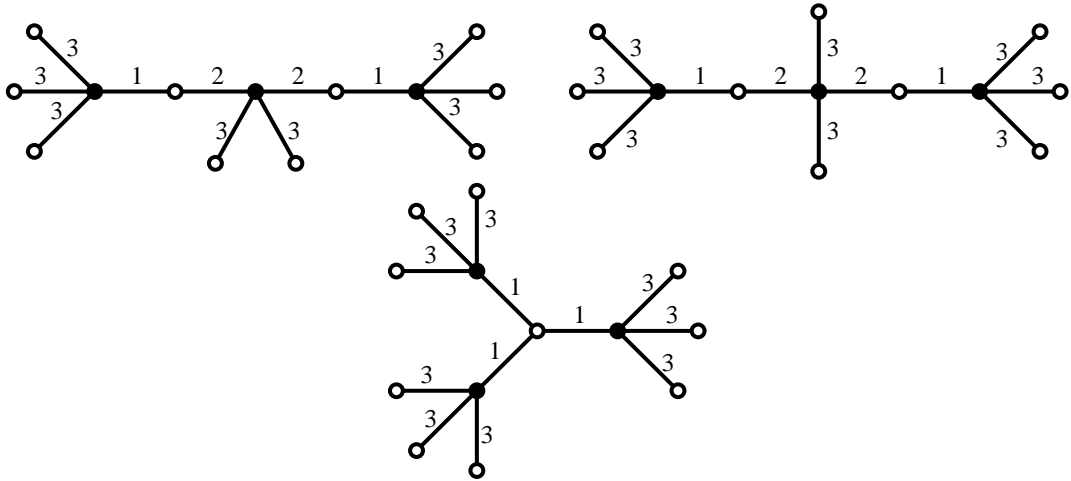,width=14cm}
\caption{\small Three trees corresponding to the passport $(10^3,3^{10})$.} 
\label{fig:split-10-3}
\end{center}
\end{figure}

\end{example}

\begin{example}[Further sporadic DZ-triples]
The next example given in \cite{BeuSte-10} corresponds to the passport
$(5^4,4^5)$. This time, there are three trees: one of them is symmetric
with the symmetry of order~2; another one is symmetric with the symmetry
of order~4; the third one is asymmetric. All the three are therefore
defined over $\Q$.

\ssk

For the passport $(6^5,5^6)$ there are four trees. One of them is symmetric
with the symmetry of order~5; two are symmetric with the symmetry of
order~2; the remaining tree is asymmetric. Therefore, the combinatorial
orbit containing four trees splits into three Galois orbits. The asymmetric
tree corresponds to the sporadic example given in \cite{BeuSte-10}.

We leave it to the reader to draw the trees in question.
\end{example}

\begin{example}[When nothing works]\label{ex:why}
All known combinatorial invariants fail to explain why the tree with the 
passport $(9^5,5^9)$ shown in Figure~\ref{fig:why} is defined over $\Q$.

\begin{figure}[htbp]
\begin{center}
\epsfig{file=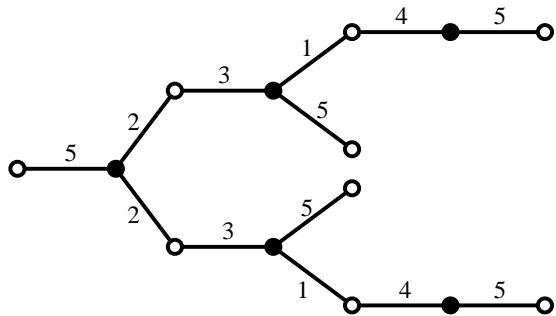,width=7.2cm}
\caption{\small This tree, corresponding to the passport $(9^5,5^9)$,
is defined over $\Q$. All known combinatorial invariants of Galois
action fail to explain this phenomenon.} 
\label{fig:why}
\end{center}
\end{figure}

The trees corresponding to this passport have 13 edges, so the outer face
is of degree 13, which is prime. This implies that this tree cannot be a
composition. Indeed, the outer face of a tree $D_F$ corresponding 
to a composition $F=f\circ h$ can only be ramified over the outer face of 
the tree $D_f$, so the degree of the outer face of $D_F$ should be the
product of the degree of the outer face of $D_f$ and of $\deg h$. But
13 cannot be a product of two integers.

Now, the monodromy group cannot be special because of Theorem~\ref{th:jones}. 
The tree cannot be self-dual either since its diameter is greater than~4, etc. 

For the moment, this example is the only one of its kind. 
However, one cannot hope to reduce the whole body of Galois theory to 
combinatorics. Note nevertheless that the direction of ``Diophantine 
invariants'' (see the next section) for the weighted trees remains 
entirely unexplored.

This example is also borrowed from the paper by Beukers and Stewart
\cite{BeuSte-10}. There are no trees in their paper; the DZ-triple 
corresponding to this example, as well as several other sporadic triples,
are found by a brute force computation using Gr\"obner bases.
\end{example}

\section{Further questions}
\label{sec:further}

Here we discuss briefly some possibilities for further research. 
To begin with, there are quite a few results known for ordinary trees, 
which might eventually be generalized to weighted trees.

\paragraph{Enumeration of weighted trees.} It would be very interesting to
find an enumerative formula which would count the number of weighted
trees having a given passport. However, this problem may turn out to be
very difficult because of the fact that the same passport can be realized
by a tree and by a forest. Therefore, an inclusion-exclusion procedure
might be necessary, preventing a nice closed formula of Goulden--Jackson's 
type (see formula (\ref{eq:GJ})). 

Right now we can prove only much more modest results, namely, enumerate
the weighted trees by their weight and number of edges. We formulate
these results without proof. Let us call a tree with a distinguished edge 
{\em edge-rooted}.

\begin{theorem}[Some enumerative results]\label{th:generating}
Let $a_n$ be the number of edge-rooted weighted bicolored plane trees of
weight~$n$. Then the generating function $f(t)=\sum_{n\ge 0}a_nt^n$ is
equal to
\begin{eqnarray*}
f(t) & = & \frac{1-t-\sqrt{1-6\,t+5\,t^2}}{2\,t} \\
     & = & 1+t+3\,t^2+10\,t^3+36\,t^4+137\,t^5+543\,t^6+2219\,t^7+9285\,t^8
           +\ldots 
\end{eqnarray*}
The asymptotic formula for the numbers $a_n$ is
\begin{eqnarray*}
a_n \sim \frac{1}{2}\sqrt{\frac{5}{\pi}}\cdot 5^n\,n^{-3/2}.
\end{eqnarray*}

Let $b_{m,n}$ be the number of edge-rooted weighted bicolored plane trees
of weight~$n$ with $m$ edges. Then the generating function
$h(s,t)=\sum_{m,n\ge 0}b_{m,n}s^mt^n$ is equal to
\begin{eqnarray*}
h(s,t) &=& \frac{1-t-\sqrt{1-(2+4s)\,t+(1+4s)\,t^2}}{2st} \\
       &=& 1 + st + (s+2s^2)\,t^2 + (s+4s^2+5s^3)\,t^3 + 
             (s+6s^2+15s^3+14s^4)\,t^4 + \ldots
\end{eqnarray*}
The following is an explicit formula for the numbers $b_{m,n}$:
\begin{eqnarray*}
b_{m,n} \eq \binom{n-1}{m-1}\cdot{\rm Cat}_m 
       \eq \binom{n-1}{m-1}\cdot\frac{1}{m+1}\binom{2m}{m},
\end{eqnarray*}
where ${\rm Cat}_m$ is the $m$th Catalan number.

Let $c_n$ be the number of non-isomorphic non-rooted trees of weight $n$,
each counted with the factor $1/|{\rm Aut}|$. Then
$$
c_n \eq \sum_{m=1}^n\frac{b_{m,n}}{m}\,.
$$
\end{theorem}

The sequence $a_n$ is listed in the On-Line Encyclopedia of Integer 
Sequences \cite{OEIS} as the entry A002212. It has many interpretations; 
the one corresponding to the weighted trees is submitted by Roland Bacher.

\paragraph{Inverse enumeration problem.} The problem is formulated as follows:
{\em For a given $m\ge 1$, classify all passports and corresponding weighted 
trees such that there exist exactly $m$ trees having this passport.}\/ 
In our paper, we have solved this problem for $m=1$. The following
result for ordinary trees was proved in \cite{Adrianov-07}. It does not
provide a classification, but nevertheless gives some important information
concerning a general pattern for the eventual classifications.

\begin{theorem}[Combinatorial orbits of a given size]
For any $m\ge 1$ the combinatorial orbits of ordinary trees containing 
exactly $m$ elements are classified as follows:
\begin{itemize}
\item   the series of chain trees\/ {\rm (}only for\/ $m=1${\rm )};
\item   a finite number of series of diameter\/ $4$;
\item   a finite number of series of diameter\/ $6$;
\item   a finite number of sporadic orbits whose elements have at most 
        $12m+1$ edges.
\end{itemize}
\end{theorem}

Our results for weighted trees and for $m=1$ fall into line with this 
pattern, only the chains must be replaced by brushes, and the bound 
$12m+1$ must be increased. It would be interesting to see if a similar 
theorem is valid for the general case.

\paragraph{Generic-sporadic splitting.} For the following three passports
for ordinary trees:
\begin{itemize}
\item   $(4^11^{n-4},p^2q^2)$: a series of trees of diameter 4; here 
        $n=2p+2q$ and $p\ne q$;
\item   $(4^p,q^21^{n-2q})$: a series of trees of diameter 6; here $n=4p$;
\item   $(4^31^8,2^{10})$: sporadic trees of diameter 8; here $n=20$
\end{itemize}
the combinatorial orbits consist of two (ordinary) trees, but one of these
trees is symmetric while the other one is not. Therefore, both trees are
defined over~$\Q$.

\begin{figure}[htbp]
\begin{center}
\epsfig{file=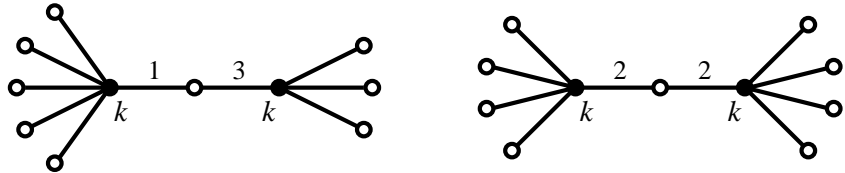,width=11cm}
\caption{\small The combinatorial orbit consists of two trees, but it
splits into two Galois orbits since one of the trees is symmetric while
the other one is not. The degrees of the black vertices are both equal
to $k\ge 3$; all leaves are of weight 1.} 
\label{fig:sym-and-not}
\end{center}
\end{figure}

For the weighted trees, we have seen similar examples in the previous
section: they correspond to the passports were $(7^3,3^7)$ and $(8^33^8)$. 
An infinite series of weighted trees with the same property is shown in 
Figure~\ref{fig:sym-and-not}. It would be interesting to produce a complete 
classification of such cases.

\vspace{-5mm}

\paragraph{Special weighted trees.} 
The condition of {\em planarity}, being added to Theorem~\ref{th:jones},
imposes strong constraints on the numbers of vertices and faces of the
corresponding maps. In~\cite{AdKoSu-97}, the complete list of ordinary
special trees is compiled. This list is finite: it contains 48~trees, the
biggest degree being~31. We are pretty sure that the list of special
weighted trees is also finite and can be drawn up.

\paragraph{Diophantine invariants.}
The following two examples may be found in \cite{LanZvo-04}. The first one
is due to Adrianov. Consider the following passport for ordinary trees:
$(5^11^{n-5},p^2q^3)$ (here $n=2p+3q$ and $p\ne q$). It is easy to see that 
there exist exactly two trees having this passport, and neither of 
them is symmetric. A simple computation shows that these trees are defined
over the field $\Q\,(\sqrt{\Delta})$ where $\Delta=3\,(p+2q)\,(2p+3q)$.
Now, if we take, for example, $p=6k^2-3l^2$ and $q=2l^2-3k^2$, choosing
$k$ and $l$ in such a way that $p$ and $q$ become positive and not equal, 
we get $\Delta$ to be a perfect square. Therefore, both trees become 
defined over $\Q$, and this splitting of the combinatorial orbit into 
two Galois orbits does not have any specific combinatorial reason: 
it is due to certain Diophantine relations between vertex degrees.
Once again, it would be interesting to extend this scheme to weighted trees.

The next example is maybe
the most spectacular one. We consider the (ordinary) trees corresponding 
to the passport $(7^11^{n-7},p^2q^5)$ (here $n=2p+5q$ and $p\ne q$). It is 
easy to see that there exist exactly three ordinary trees having this 
passport. Therefore, they are defined over a cubic extension of $\Q$; 
the cubic polynomial generating this field may be written explicitly. 
Now, we ask the following question: is it possible that this polynomial 
has a rational root? If yes, then the combinatorial orbit in question 
will split into two Galois orbits, one defined over $\Q$, and the other 
one quadratic.

It turns out that the search for polynomials having a rational root
can be reduced to the search for rational points on a particular elliptic 
curve. The curve in question contains infinitely many rational points.
We have computed the first 11 solutions. The smallest one corresponds
to trees having $n=686$ edges ($p=33$, $q=124$); the 11th solution 
corresponds to trees having $n\approx 3.45\cdot 10^{134}$ edges. 
Similar constructions certainly must also exist for weighted trees.

\begin{figure}[htbp]
\begin{center}
\epsfig{file=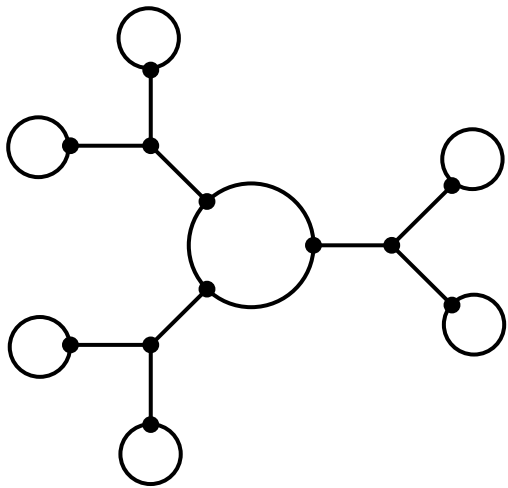,width=6.6cm}
\caption{\small This map represents two polynomials $A$ and $B$, of degrees
$2k=12$ and $3k=18$ respectively, such that $\deg\,(A^3-B^2)=9$. Thus, 
the degree of the difference does not attain its minimal value $k+1=7$, 
but in return both $A$ and $B$ are defined over $\Q$.} 
\label{fig:deg-relaxed}
\end{center}
\end{figure}

\begin{figure}[htbp]
\begin{center}
\epsfig{file=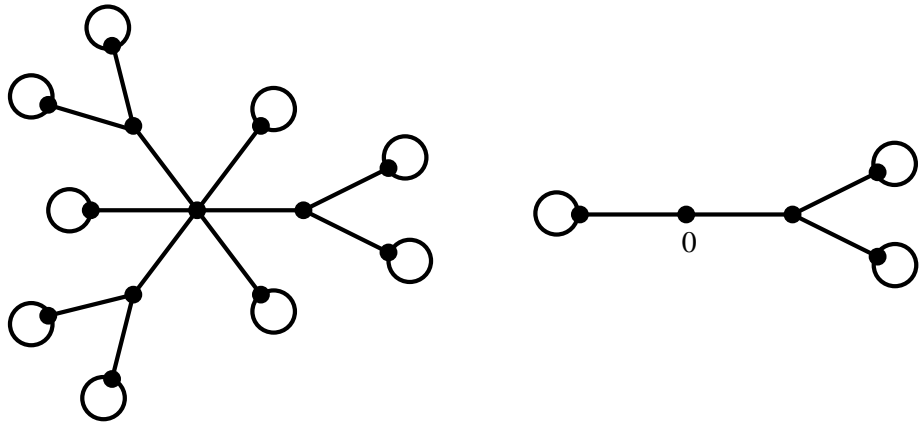,width=12cm}
\caption{\small This map represents two polynomials $A$ and $B$, of degrees
$2k=14$ and $3k=21$ respectively, such that $\deg\,(A^3-B^2)=9$. Thus, 
the degree of the difference does not attain its minimal value $k+1=8$, 
but in return both $A$ and $B$ are defined over $\Q$.} 
\label{fig:deg9-instead-8}
\end{center}
\end{figure}

\paragraph{Relaxing the minimum degree condition.} Let us revisit the initial
problem about the minimum degree of the polynomial $A^3-B^2$, see 
page~\pageref{init-problem}. When there are no DZ-triples defined
over~$\Q$, we may relax the condition of the $\deg R$ being the least
possible and thus obtain solutions with bigger $\deg R$ but, in return, 
defined over $\Q$. Two example of this kind are shown in Figures
\ref{fig:deg-relaxed} and~\ref{fig:deg9-instead-8}. In the first one,
$k=6$ but $\deg R=9$ instead of~7 since one of the faces is of degree~3
instead of~1. In the second example, $k=7$ but $\deg R=9$ instead
of~8 since, instead of two black vertices of degree~3 we have here
one black vertex of degree~6.

Now let us look at the second example. Though a computation
of the Belyi function in this case is not difficult, it is still 
interesting to analyze this example in purely combinatorial terms.
The map shown on the right in Figure~\ref{fig:deg9-instead-8}
is a unimap $s$ of Figure~\ref{fig:trees2maps}, which is also equal to
the unitree $S$ in Figure~\ref{fig:sporadic-8}. Therefore, it is defined
over $\Q$. Its black vertex of degree~2 is a bachelor 
(Definition~\ref{def:bachelor}); therefore, it can be placed at any
rational position (Remark~\ref{rem:bachelor}), for example, at the
point $x=0$. Then it remains to insert $x^3$ instead of $x$ in its
Belyi function, and we get a Belyi function for the bigger ``triple'' 
dessin. This example shows that the possibilities of the combinatorial 
approach to this problem are far from being exhausted.

In general, it would be interesting to establish an upper bound on the 
difference between the minimum degree attainable in $\C\,[x]$, and the 
one attainable in~$\Q\,[x]$.

\pagebreak[4]

\paragraph{Unimaps.} Besides the maps with all finite faces of degree~1,
there are other classes of maps for which the question of uniqueness
is interesting. One such example is the class of maps with exactly
two faces: their Belyi functions are Laurent polynomials. The existence
questions for such maps were completely settled in \cite{Pakovich-09};
the uniqueness remains to be studied. Some new Galois phenomena related
to the non-existence of bachelors appear in this case, see
\cite{Couveignes-94}.

\begin{center}
{\bf Acknowledgements}
\end{center}

Fedor Pakovich is grateful to the Bordeaux-I University, France, and 
Alexander Zvonkin is grateful to the Center for Advanced Studies in Mathematics
of the Ben-Gurion University of the Negev,
Israel, for their mutual hospitality. Fedor Pakovich is also grateful 
to the Max-Planck-Institut f\"ur Mathematik, Bonn, where the most part 
of this paper was written.


\end{document}